\documentstyle{amsppt} 

\input epsf.tex
\input MosherMacros.sty 
\epsfverbosetrue

\hsize = 6.5 true in

\topmatter

\title\nofrills A User's Guide to the Mapping Class Group: \break
Once Punctured Surfaces \endtitle

\author Lee Mosher \endauthor

\address {Mathematical Sciences Research Institute, Berkeley CA 94720}
\endaddress

\email mosher\@msri.org \endemail

\date 
\vbox{
\centerline{September, 1994}
}
\enddate

\thanks The author was partially supported by NSF grant \# DMS-9204331
\hfill \break
\indent Research at MSRI partially supported by NSF grant \# DMS-9022140
\endthanks

\endtopmatter

\document

\def\centeredepsfbox#1{\centerline{\epsfbox{#1}}}

\def\prmit{{\sl prmit}}

\def\bem{{\sl bem}}
\def\bems{{\sl bems}}

\def\mcg{\func{\cal M \cal C \cal G}}
\def\mcgd{\func{\cal M \cal C \cal G \cal D}}
\def\Homeo{\func{Homeo}}
\def\Ends{\func{\cal E}}
\def\OppEnd{\func{Opp}}

\def\Prongs{\func{\cal P}}

\def\TailEnd{\func{Tail}}

\def\pred{\func{Pred}}
\def\Pred{\pred}
\def\succ{\func{Succ}}
\def\Succ{\succ}
\def\Next{\func{Next}}

\def\Opp{{\hbox{\eightpoint Opp}}}

\def\skel(#1,#2){#1^{(#2)}}

\def\Relabelling#1{@> \func{Rotate}(#1) >>}

\def\omicron{o}

\def\congruent{\equiv}
\def\modulo#1{\pmod{#1}}

An automatic structure for the mapping class group of a surface of finite type
was described in \cite{M}. This document is intended as a practical guide to
computations using a variant of this automatic structure, in the special case of
a once-punctured, oriented surface $S$. As such, we shall try to be more
descriptive and less theoretical than in \cite{M}, leaving the reader to consult
\cite{M} for detailed proofs. Our primary goal is that the reader may learn, as
quickly as possible, how to compute in the mapping class group of a
once-punctured surface: we describe a quadratic time algorithm for the word
problem, henceforth called \em{the algorithm}, which can be implemented with
pencil and paper. A {\it Mathematica} version of the algorithm is (or will soon
be) available; check the software library at the Mathematical Sciences Research
Institute (e-mail address: msri.org), or the Geometry Center (geom.umn.edu).

As with any computational method, it is necessary to learn some of the theory
in order to learn the algorithm. We spend some time developing various
combinatorial tools, with enough justification supplied to aid understanding and
lessen the steepness of the learning curve. There is a trade-off involved here:
time invested understanding theory may be time wasted gaining proficiency; I do
not know if I have found the right balance. Also, despite my stated purpose, in
a few places I have put in perhaps too much detail about items of combinatorica
that interest me, but which are not really to the point, so the reader is
forewarned.

The algorithm described herein can be adapted to arbitrary punctured surfaces,
with or without boundary and orientation. The data structures needed do not
lend themselves quite so nicely to pencil and paper calculation, so we do not
pursue the issue here; details can be found in \cite{M}. And while an automatic
structure for the mapping class group of a closed surface is described in
\cite{M}, in this case the results are \em{not} suited for practical
calculation, because of the non-constructive nature of the proof; hopefully a
practical automatic structure will emerge from a deeper understanding of closed
surfaces.

For the rest of the paper, let $S$ be an oriented, once-punctured surface which
is not the 2-sphere. We regard $S$ as a closed surface with a distinguished
point $p$, the \em{puncture}. The mapping class group is $\mcg(S) = \Homeo(S) /
\Homeo_0(S)$, where $\Homeo(S)$ is the group of all orientation preserving
homeomorphisms of~$S$ fixing~$p$, and $\Homeo_0(S)$ is the normal subgroup of all
homeomorphisms isotopic to the identity leaving $p$ stationary throughout the
isotopy. 

We shall describe an explicit 2-complex $X = X(S)$ whose fundamental group is the
mapping class group $\mcg(S)$. The set of homotopy classes of edge paths in any
complex form a groupoid under the operation of concatenation, called the
\em{edge path groupoid} of that complex; a more descriptive but longer name would
be ``edge path homotopy groupoid'', but we stick with the shorter name. In
the particular case of $X(S)$, the edge path groupoid will be called the
\em{mapping class groupoid}, denoted $\mcgd(S)$. 

Recall that a \em{combing} (with uniqueness) for the edge path groupoid of $X$
consists of a base vertex in $X$, and a choice of a unique representative called
the \em{normal form} for each homotopy class of edge paths in $X$ starting at
the base vertex. When these normal forms satisfy certain computational
properties we say that the set of normal forms is an \em{automatic structure}.
First, there is a finite automaton which checks whether a given path is a normal
form, so the set of normal forms is a regular language. Second, for each edge in
$X$ there is a finite automaton called a \em{multiplier automaton}, which checks
whether two normal forms differ by that edge, up to homotopy. The second
condition can be replaced by the equivalent \tit{fellow traveller property}: if
two normal forms $v,w$ differ up to homotopy by a generator, then (letting
$v(t)$ be the prefix of length $t$ of $v$, or $v(t)=v$ is $t$ is greater than
the length of $v$) we have that $v(t)^\inverse w(t)$ is
homotopic to a path whose length is bounded, by a constant independent of
$v,w,t$. Another important notion is that of an \em{asynchronous automatic
structure}, where the fellow traveller property is replaced by the weaker
\em{asynchronous fellow traveller property}: if two normal forms $v,w$ differ
up to homotopy by a generator, then there are sequences $0=s_0<s_1<\cdots<s_M$
and $0=t_0<t_1<\cdots<t_M$ with bounded differences $s_{i+1}-s_i$,
$t_{i+1}-t_i$, such that $v(s_i)^\inverse w(t_i)$ is homotopic to a path of
bounded length. By contrast, an ordinary automatic structure is sometimes
called a \em{synchronous} automatic structure. The reader
is referred to \cite{ECHLPT} for formal definitions.

In \cite{M}, asynchronous and synchronous automatic structures on the groupoid
$\mcgd(S)$ are described. Associated to the edge path groupoid on a complex is
the group of homotopy classes of closed edge paths, the fundamental group. By a
general result of \cite{ECHLPT}, given an automatic structure on a groupoid one
may obtain an automatic structure on the associated group. By another general
result, given an automatic structure on a group (or groupoid) one may obtain a
quadratic time algorithm for computing normal forms. Combining these results, we
obtain a quadratic time algorithm for the word problem in $\mcgd(S)$ or
$\mcg(S)$. 

Instead of appealing to these general results, we directly construct a quadratic
time algorithm for computing normal forms in the groupoid $\mcgd$, and by
restricting the input to closed edge paths one obtains an algorithm for the
group $\mcg$. Our normal forms will come from the \em{asynchronous} automatic
structured described in \cite{M}. The algorithm we describe for computing these
normal forms will run in quadratic time; this will be proved by comparing the
asynchronous automatic structure to another, synchronous automatic structure. Our
proof will use a special property of the normal forms for $\mcgd$, the ``suffix
uniqueness property'', discussed in section IV. Because of the suffix uniqueness
property, our quadratic time algorithm is more efficient than the one described
in \cite{ECHLPT}.

The automatic structure we describe is very large. As a function of the genus
$g$, the number of states in the word acceptor grows at least as fast as $g^g$;
see figure 19. As such, our algorithm does \em{not} require explicitly computing
and storing the word acceptor and multiplier automata. In a sense, our algorithm
constructs local portions of the automata as they are needed for calculation.

Thus, despite the size of the automatic structures the computations are quite
efficient. With practice, the algorithm can be implemented quite efficiently
with pencil and paper by drawing lots of funny pictures called ``chord
diagrams''; we give profuse illustrations of such calculations. The author is
able to compute the normal forms for a once-punctured surface of genus 2,
starting from an edge path of length $n$, in at most $25 n^2$ minutes, given a
sufficient quantity of paper, pencils, and coffee. In actual practice, the
computations are much faster (although errors, and the correction thereof using
an eraser, may slow down computation time).

Another issue arises from the fact that the complex $X$ is so large, so one
would not want to write down a presentation for $\mcg$ using $X$. This raises
the question of what form the algorithm uses for input. The direct form of input
is an edge path in $X$. However, there are well-known ``small'' presentations
for $\mcg(S)$ described in the literature whose size grows linearly with the
genus. The \tit{Mathematica} implementation of the algorithm allows the user to
input a word in standard generators. This word is converted into an edge path
in $X$ in linear time, and then the algorithm works on the edge path. We shall
not describe this conversion process here.

\head I. The complex $X$ \endhead

In this section we construct a finite complex $X$ whose fundamental group is
$\mcg$. First we construct a contractible complex $Y$ on which $\mcg$ acts with
finite cell stabilizers and finitely many orbits; the complex $Y$ was first
described by Harer \cite{Har}. We are only interested in the 2-skeleton
$\skel(Y,2)$. Then we resolve finite cell stabilizers of $\skel(Y,2)$ to obtain
a 2-complex $\tilde X$ on which $\mcg$ acts freely with finitely many cell
orbits. The quotient of $\tilde X$ by $\mcg$ yields $X$.	

\subhead Ideal arc systems and the complex $Y$ \endsubhead

Let $I$ be the closed unit interval. An \em{ideal arc} is the image $h$ of a
map 
$$(I,\bd I,\interior I) \mapsto (S,p,S-p)
$$ 
which is injective in $\interior I$, such that $h$ does not bound a disc; this
map is called a \em{characteristic map} of $h$. The image of $\interior I$ is
called the \em{interior} of $h$, denoted $\interior(h)$. Two ideal arcs $h,h'$
are \em{isotopic} if there exists $\phi \in \Homeo_0(S)$ such that $\phi(h) =
h'$. In general, for any set of objects on which $\Homeo(S)$ acts, two objects
are \em{isotopic} if they differ by an element of $\Homeo_0(S)$. 

Given two ideal arcs $h_1,h_2$, it is easy to decide whether they are
isotopic. In the first case where the interior are disjoint, then $h_1,h_2$ are
isotopic if and only if they bound a disc. In the second case where the interiors
intersect, perturb $h_1,h_2$ so that the interiors have a finite number of
transverse intersection points, and then successively isotop them to remove any
complementary components which are discs, descreasing the number of intersection
points. If this process stops with a positive number of intersection points then
$h_1,h_2$ are not isotopic; otherwise we have reduced to the first case.

An \em{ideal arc system} $\gamma$ is a collection of non-isotopic ideal arcs with
disjoint interiors, such that each component of $S - \gamma$ is a disc.
Figure 1 gives several examples of ideal arc systems. Figure 1(a) is the
standard method for cutting a surface of genus 2 into an octagon. By adding
arcs to this octagon we obtain figure 1(b), which cuts the surface into six
triangles. Figure 1(c) shows another such ``ideal triangulation''. Note that
(b) and (c) are not isotopic, because (b) has a separating ideal arc while (c)
has none.

\midinsert
\centeredepsfbox{01.SomeIdealArcSystems.eps}
\botcaption{Figure 1}
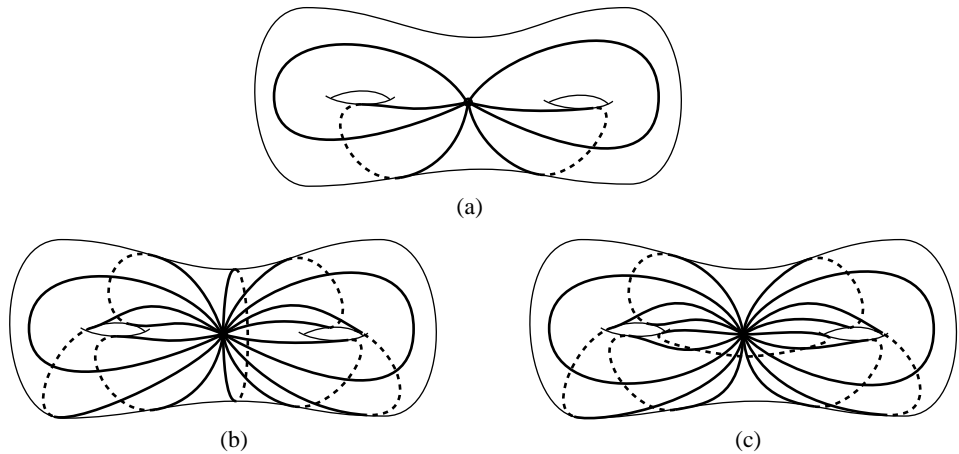
Some ideal arc systems
\endcaption
\endinsert

The group $\Homeo(S)$ acts on the set of ideal arc systems. Given an ideal arc
system $\gamma$ and $\phi \in Homeo_0(S)$, if $\phi(\gamma) = \gamma$ then
$\phi$ fixes each arc of $\gamma$ setwise (because distinct arcs are pairwise
non-isotopic), and preserves the orientation (because $\pi_1(S)$ has no
torsion). Thus, $\phi$ preserves the ends of arcs in $\gamma$. We refer to this
fact as \tit{Rigidity of Ends}.

A complementary component $C$ of an ideal arc system $\gamma$ is called a
\em{polygon} of $\gamma$. There is a characteristic map $D \mapsto S$ for $C$,
where $D$ is a convex Euclidean polygon, so that each vertex of $D$ goes to
$p$, each side of $D$ gives a characteristic map of some ideal arc of $\gamma$,
and $\interior(D)$ goes to $C$. The number of sides of $C$ is defined to be the
number of sides of $D$. Polygons are referred to as \em{triangles,
quadrilaterals, pentagons} etc.\ depending on the number of sides; in general an
$n$-sided polygon is called an $n$-gon. Note that there are no 1-gons or 2-gons.

Examples: The polygons of $\gamma$ are all triangles if and only if $|\gamma| =
6g-3$. At the other extreme, $\gamma$ has a single polygon if and only if
$|\gamma| = 2g$, in which case the polygon is a $4g$-gon. These facts are easily
verified using the Euler characteristic.

Given a polygon, there is a certain well-defined number of ideal arcs that can
be added to triangulate the polygon: for a quadrilateral, add 1 arc; for a
pentagon add 2 arcs; for an $n$-gon add $n-3$ arcs (while there are many
different ways to add these arcs, the number of arcs added is always $n-3$; the
number of distinct ways to add the arcs, up to isotopy, is given by the
Catalan number ${1 \over m+1} {2m \choose m}$ where $m = n-2$ \cite{STT}). It
follows that for an ideal arc system $\gamma$, by adding over all polygons one
obtains a certain well-defined number of ideal arcs that can be added to
triangulate $\gamma$; this number is called the \em{defect} of $\gamma$. Since a
triangulation always has $6g-3$ arcs, the defect is equal to
$6g-3-|\gamma|$. 

If one ideal arc system $\gamma'$ is obtained by adding arcs to another ideal
arc system $\gamma$, then we say that $\gamma'$ is a \em{refinement} of
$\gamma$. For example, every $\gamma$ can be refined to become an ideal
triangulation.

Now we construct $Y$. In general $Y$ has one $k$-cell for each isotopy class of
ideal arc systems of defect $k$. For example, the 0-cells of $Y$ are in 1-1
correspondence with isotopy classes of ideal triangulations. To describe the
attaching maps of cells, suppose the $(k-1)$-skeleton $Y^{(k-1)}$ has been
constructed, and let $\gamma$ be an ideal arc system of defect $k$. The isotopy
classes of all possible refinements of $\gamma$ form a subcomplex of
$\skel(Y,k-1)$. Now check that this subcomplex is a topological $(k-1)$-sphere,
and attach a $k$-cell to this sphere corresponding to $\gamma$. We will
explicitly study the attaching maps for $1$ and $2$-cells below.

The complex $Y$ was first described in \cite{Har}, where $Y$ is proved to be
contractible using Strebel differentials. See \cite{Hat} for an elementary proof
of contracitbility.

Since the number of arcs in an ideal arc system is always at least $2g$, it
follows that the defect is always at most $6g-3-2g=4g-3$, so the dimension of
$Y$ is $4g-3$. Harer proves \cite{Har} that this is the minimum possible, by
showing that the virtual cohomological dimension of $\mcg$ is $4g-3$.

\subhead $0$-cells of $Y$: ideal triangulations \endsubhead

The 0-cells of $Y$ are in 1-1 correspondence with isotopy classes of ideal
triangulations of $S$. Figures 1(b,c) show two non-isotopic ideal
triangulations. In general there are infinitely many isotopy classes of ideal
triangulations, because $\mcg$ is infinite, and for each ideal triangulation
$\delta$ there are only finitely many $\Phi \in \mcg$ such that $\Phi[\delta] =
[\delta]$ (this follows from \tit{Rigidity of Ends}). In fact, we will see that
the stabilizer of $[\delta]$ is a finite cyclic subgroup of $\mcg$ whose order
can take only finitely many values depending on the genus. 

This begs the question: how many orbits of ideal triangulations are there under
$\mcg$? This question is particularly easy to answer on a once-punctured torus.
We shall take this up later.

\subhead $1$-cells of $Y$: elementary moves \endsubhead

The $1$-cells of $Y$ are in 1-1 correspondence with isotopy classes of ideal arc
systems $\gamma$ of defect 1. The polygons of $\gamma$ consist of one
quadrilateral $Q$ and the rest triangles. An example is given in figure 2, with
$Q$ shaded. A quadrilateral can be triangulated by inserting an
ideal arc, and there are two ways to do this insertion, as shown schematically in
figure 3 (a schematic picture like this can be thought of as the domain of a
characteristic map for the quadrilateral, so all vertices will be identified to
the puncture, and there may be certain pairwise side identifications as well).
This gives two ideal triangulations, forming two 0-cells of $Y$, to which a
1-cell is attached corresponding to $\gamma$ (see figure 3). Using the $\gamma$
given in figure 2, the two ideal triangulations refining $\gamma$ are
shown in figures 1(b,c), producing a particular 1-cell in $Y$.

\midinsert
\centeredepsfbox{02.Defect1ArcSystem.eps}
\botcaption{Figure 2} 
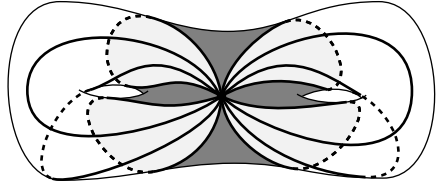
An ideal arc system of defect 1. The quadrilateral, part of which wraps around
back, is shaded. 
\endcaption\endinsert

\midinsert
\centeredepsfbox{03.TriangulatingAQuadrilateral.eps}
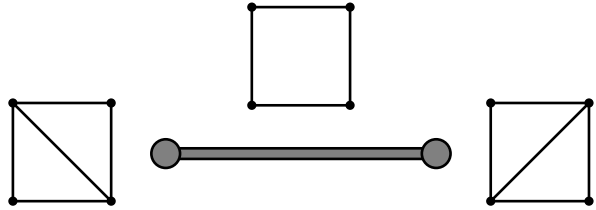
\botcaption{Figure 3} The two triangulations of a quadrilateral,
giving two 0-cells at the ends of a 1-cell in $Y$.
\endcaption\endinsert

Here is a way to think about an \em{oriented} 1-cell of $Y$. Start with an
ideal triangulation $\delta$ and an arc $h$ of $\delta$. Notice that $h$ cannot
lie on two sides of a single triangle of $\delta$, for when two sides of a
triangle are identified then there must be two or more punctures; see figure 4.
Therefore, when $h$ is removed from $\delta$ the two adjacent triangles form a
quadrilateral $Q$ in a defect 1 ideal arc system $\gamma$. The arc $h$ forms one
diagonal of $Q$; let $h'$ be the opposite diagonal. Inserting $h$ into $\gamma$
yields a new ideal triangulation $\delta'$. We shall indicate this operation by
saying that $\delta \to \delta'$ is an \em{elementary move}, and that the
elementary move is \em{performed on $h$}, with \em{opposite diagonal
$h'$}. To emphasize the role of $h$ we also write $\delta \mapright{h} \delta'$,
and in pictures such as figure 5 we thicken $h$. Also, we say that the
quadrilateral $Q$ is the \em{support} of the elementary move $\delta \to
\delta'$. To summarize, there is a natural 1-1 correspondence between: oriented
1-cells of $Y$, isotopy classes of pairs $(\delta,h)$, and pairs of the form
$([\delta],[\delta'])$ where $\delta \to \delta'$ is an elementary move; the
corresponding 1-cell in $Y$ is denoted $[\delta \to \delta']$.

\midinsert
\centeredepsfbox{04.FoldedTriangle.eps}
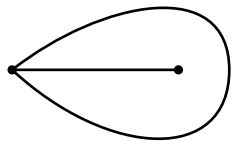
\botcaption{Figure 4} In an ideal triangulation $\delta$ on an oriented
surface, if a triangle has two sides identified then that surface must have at
least two punctures.
\endcaption\endinsert

Note that elementary moves are symmetric: if $\delta \to \delta'$ is an
elementary move then so is $\delta' \to \delta$. Note also that for each
ideal triangulation $\delta$, there are $6g-3$ outgoing elementary moves, one
performed on each ideal arc of $\delta$; by reversing the directions we also
see that there are $6g-3$ incoming elementary moves.

Returning to our only example so far, from figures 1(b,c) we obtain the
elementary move $\delta \mapright{h_1} \delta_1$ shown in figure 5(a), with
intervening defect 1 ideal arc system $\gamma_1$. Figure 5(b)
shows another elementary move $\delta \mapright{h_2} \delta_2$ starting from the
same $\delta$, with intervening arc system $\gamma_2$. Note that $\gamma_2$ has
a quadrilateral with one pair of opposite sides identified. Also, note that
$\delta$ and $\delta_2$ differ by a mapping class, namely the Dehn twist around
the core curve of the handle on the right side of the surface. Therefore, the
0-cells of $Y$ corresponding to $[\delta]$ and $[\delta_2]$ are in the same orbit
under the action of $\mcg$.

\midinsert
\centeredepsfbox{05.SomeElementaryMoves.eps}
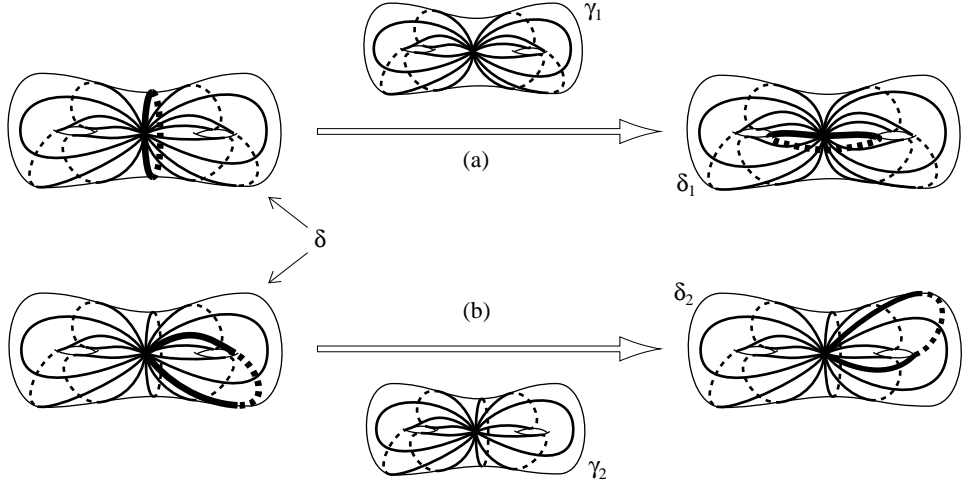
\botcaption{Figure 5} Some elementary moves
\endcaption\endinsert

\subhead $2$-cells of $Y$: commutator and pentagon relators \endsubhead

The $2$-cells of $Y$ are in 1-1 correspondence with ideal arc systems $\gamma$
of defect 2. This can happen in two ways: the polygons of $\gamma$ can consist
of two quadrilaterals and the rest triangles; or one pentagon and the rest
triangles. 

If $\gamma$ has two quadrilaterals, each quadrilateral can be independently
triangulated in one of two ways, yielding four distinct triangulations which
refine $\gamma$. The four corresponding 0-cells in $Y$ are connected up by four
1-cells as shown in figure 6, forming a closed edge path in $Y$ of length four.
Attached to this edge path is a 2-cell of $Y$ corresponding to $\gamma$. If
$\delta$ is one of the four triangulations, and if $h_1,h_2$ are the two
diagonals inserted into $\gamma$ to form $\delta$, then the two adjacent sides
of the 2-cell are given by elementary moves $\delta \mapright{h_n} \delta_n$. There
are also elementary moves $\delta_1 \mapright{h_2} \delta'$ and $\delta_2
\mapright{h_1} \delta'$, as shown in figure 6. For this reason, we can say that
the elementary moves performed on $h_1$ and on $h_2$ commute with each other, and
thus we say that this attached 2-cell is a \em{commutator relator}.

\midinsert
\centeredepsfbox{06.CommutatorRelator.eps}
\botcaption{Figure 6} A commutator relator
\endcaption\endinsert

Now suppose $\gamma$ has one pentagon. This pentagon can be triangulated in one
of five ways as shown in figure 7, forming a closed edge path in $Y$ of
length five, to which a $2$-cell is attached. This 2-cell is called a
\em{pentagon relator}. Note that if $\delta$ is one of the five ideal
triangulations, if $h_1$, $h_2$ are the arcs inserted into $\gamma$ to form
$\delta$, and if $\delta \mapright{h_1} \delta_1$ and $\delta \mapright{h_2}
\delta_2$ are the two sides of the relator incident to $\delta$, then the next
two sides are $\delta_1 \mapright{h_2} \delta'_1$ and $\delta_2 \mapright{h_1}
\delta'_2$, and there is a fifth side $\delta'_1 \to \delta'_2$. 

\midinsert
\centeredepsfbox{07.PentagonRelator.eps}
\botcaption{Figure 7} A pentagon relator
\endcaption\endinsert

One computation which arises over and over is the following. Given an ideal
triangulation $\delta$ and ideal arcs $h_1 \ne h_2 \in \delta$, consider the two
elementary moves $\delta \mapright{h_1} \delta_1$ and $\delta \mapright{h_2}
\delta_2$. Do these two elementary moves lie on a unique relator? If so, is it a
commutator relator or a pentagon relator? 

These questions can be answered by examining the adjacencies of ends of
$h_1,h_2$. Each ideal arc has two ends. Given an end of $h_1$ and an end of
$h_2$, these ends are \em{adjacent} in $\delta$ is they are incident to some
corner of some polygon of $\delta$. It cannot happen that $h_1$ has two ends
adjacent to a single end of $h_2$, for then we obtain a folded triangle as in
figure 4. Thus, $h_1$ and $h_2$ can have zero, one, or two pairs of
adjacent ends, and if two then the pairs are disjoint. 

If $h_1,h_2$ have no pairs of adjacent ends, then removal of $h_1,h_2$ yields an
ideal arc system with two quadrilaterals, and we obtain a commutator relator. If
$h_1,h_2$ have one pair of adjacent ends, their removal produces an ideal arc
system with a pentagon, and we get a pentagon relator. If $h_1,h_2$ have two
pairs of adjacent ends, then removal of $h_1,h_2$ creates a collection of ideal
arcs with an annulus complementary component (see figure 8 for an example).
This violates the definition of ideal arc system, so there is no 2-cell in $Y$
corresponding to this collection of ideal arcs.

\midinsert
\centeredepsfbox{08.AnnulusRegion.eps}
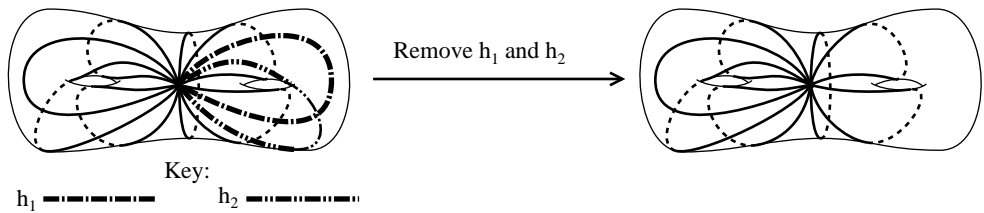
\botcaption{Figure 8} If each end of $h_1$ is adjacent to an end of
$h_2$, removal of $h_1$ and $h_2$ creates an annulus, and no relator is
obtained.
\endcaption\endinsert

\subhead The action of $\mcg$ on $Y$ \endsubhead

In general, whenever there is a set of objects on which $\Homeo(S)$ acts, then
$\mcg$ acts on the isotopy classes. Now $\Homeo(S)$ acts in the obvious way on
the set of ideal arc systems, so $\mcg$ acts on their isotopy classes. Also,
the action of $\mcg$ preserves the relation of refinement, hence $\mcg(S)$ acts
on the complex $Y$ by cellular homeomorphisms.

We need some notation for this action. In general, given an ideal arc system
$\gamma$ the isotopy class of $\gamma$ is denoted $[\gamma]$. Thus, given $\Phi
\in \mcg$ represented by $\phi \in \Homeo(S)$, then $\Phi[\gamma] =
[\phi(\gamma)]$.

\subhead Combinatorial equivalence and chord diagrams \endsubhead

The main goal of this section is to present a calculus which allows us to
understand cell stabilizers and cell orbits of the action of $\mcg$ on
$Y$. This calculus will enable us to understand why cell stabilizers are
finite, why the number of cell orbits is finite, and it will guide us in
``resolving'' finite order cell stabilizers, leading up to the definition of the
complex $X$. 

Two ideal arc systems $\gamma,\gamma'$ are said to be \em{combinatorially
equivalent}, or to have the same \em{combinatorial type}, if they are in the same
orbit under the action of $\Homeo(S)$, i.e.\ $\phi(\gamma) = \gamma'$ for some
$\phi \in \Homeo(S)$. Equivalently, their isotopy classes $[\gamma], [\gamma']$
are in the same orbit under the action of $\mcg$. The combinatorial equivalence
class of $\gamma$ is denoted $\{ \gamma \}$. We now associate to each $\gamma$ a
finitistic object called its \em{combinatorial diagram}, which will encode the
combinatorial type of $\gamma$. Then we show how to represent the combinatorial
diagram pictorially with the \em{chord diagram}.

Let $\gamma$ be an ideal arc system. Let $\Ends(\gamma)$ be the set of ends of
ideal arcs of $\gamma$. If $h$ is an ideal arc, then an end of $h$ is just an
end, in the usual sense, of the topological space $h-p$, so each ideal arc has
two ends. An end of $h$ can be represented by a \em{half-arc} of $h$, which is
the closure of a component of $h - \{p,x\}$ for some $x \in \interior(h)$. Now we
put some extra structure on $\Ends(\gamma)$. 

Recall that a \em{circular ordering} on a finite set is just a permutation with
one cycle. There is a natural way to use the orientation on $S$ to put a
circular ordering on $\Ends(\gamma)$. Choose a disc $D$ containing $p$ so that $D
\intersect \gamma$ is a union of radii of $D$. These radii form half-arcs of
$\gamma$, and they are in 1-1 correspondence with $\Ends(\gamma)$. The
orientation on $S$ determines a boundary orientation on $\boundary D$, which
determines in turn the circular ordering on $\Ends(\gamma)$. Denote this
circular ordering by $\Ends(\gamma) \mapright{\Succ} \Ends(\gamma)$, the
\em{successor} map. The inverse permutation is called the \em{predecessor} map,
denoted $\Ends(\gamma) \mapright{\Pred} \Ends(\gamma)$. Next recall that a
\em{transposition} on a finite set is a permutation where every cycle has length
2. The correspondence between opposite ends of the same arc determines a
transposition on $\Ends(\gamma)$ denoted $\Ends(\gamma) \mapright{\OppEnd}
\Ends(\gamma)$, the \em{opposite end} map. The \em{combinatorial diagram} of
$\gamma$ is defined to be the ordered triple $(\Ends(\gamma),
\OppEnd, \Succ)$. 

Given ordered triples $(\Ends_i, \omicron_i, \sigma_i)$, $i=1,2$, where $\Ends_i$
is a finite set and $\omicron_i, \sigma_i$ are permutations of $\Ends_i$, we say
these triples are \em{isomorphic} if there is a bijection $\Ends_1 \mapright{\phi}
\Ends_2$ such that $\phi \composed \sigma_1 = \sigma_2 \composed \phi$ and
$\phi \composed \omicron_1 = \omicron_2 \composed \phi$.

The fact we need is that two ideal arc systems are combinatorially equivalent
if and only if their combinatorial diagrams are isomorphic. For if $\gamma_1,
\gamma_2$ are combinatorially equivalent then the homeomorphism between them
induces an isomorphism of their combinatorial diagrams. Conversely, if their
combinatorial diagrams are isomorphic then one can construct the desired
homeomorphism up through the skeleta by induction: since the opposite end maps
correspond the homeomorphism can be extended over the 1-skeleta, and since the
successor maps correspond it can be extended over the polygons preserving
orientation.

We can immediately see why there are finitely many combinatorial equivalence
classes of ideal arc systems on $S$: the size of $\Ends(\gamma)$ is bounded by
$12g - 6$, and for a set $\Ends$ of bounded size there are only finitely many
equivalence classes of triples $(\Ends,\omicron,\sigma)$. Also, using
\tit{Rigidity of Ends} we can see why cell stabilizers are finite, because the
subgroup of $\mcg$ stabilizing $[\gamma]$ is isomorphic to the set of
automorphisms of the combinatorial diagram of $\gamma$, which is cyclic of
order bounded by $12g - 6$, because an automorphism must commute with the
successor map. The optimal order bound is somewhat smaller than $12g-6$, because
an automorphism also commutes with the opposite end map; on a surface of genus 2
the optimal order bound is 3.

The combinatorial diagram of $\gamma$ can be represented pictorially by the
\em{chord diagram}. Draw a circle on a piece of paper, oriented
counterclockwise, and draw $12g-6$ points on the circle corresponding to
$\Ends(\gamma)$, so that the counterclockwise ordering corresponds to $\Succ$.
Now draw chords connecting up the points in pairs, using the transposition
$\OppEnd$; in diagrams we use chords which are hyperbolic geodesics in the
\Poincare\ disc model, i.e.\ arcs of circles orthogonal to the boundary. Figure
9 shows three examples, two ideal triangulations and a defect 1 ideal arc
system, taken from the first elementary move pictured in figure 5. In order for
the reader to get used to the chord diagrams, in these figures we have indexed
the arc ends with integers, but we will not use any indexing in future chord
diagrams. 

In the second elementary move $\delta \to \delta_2$ of figure 5, one can check
that $\delta$ and $\delta_2$ have the same chord diagram, verifying the earlier
statement that they are combinatorially equivalent.

\topinsert
\centeredepsfbox{09.SomeChordDiagrams.eps}
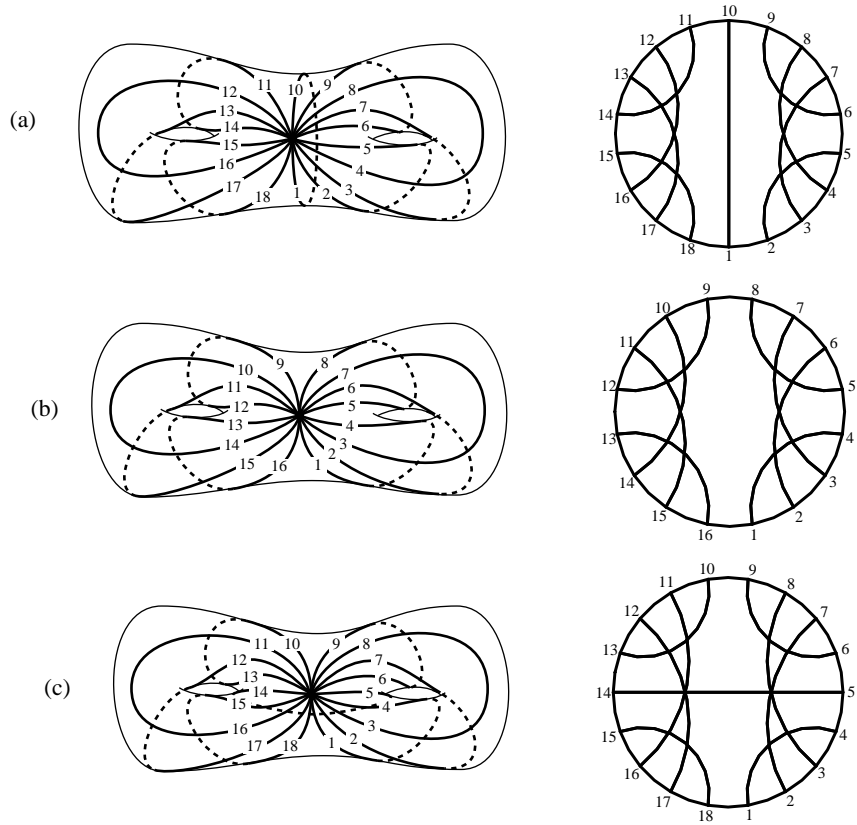
\botcaption{Figure 9} Chord diagrams of some ideal arc systems
\endcaption\endinsert

It is easy to distinguish combinatorial types by viewing the chord diagram. If
the points representing ends are spaced regularly around the circle, and if the
chords are drawn with hyperbolic geodesics, then the chord diagram itself is a
complete invariant of the combinatorial type, regarding two chord diagrams as
being isomorphic if they differ by a Euclidean similarity. It is also easy to
recognize the automorphism group of a chord diagram, by just looking for
circular symmetries of the diagram. These tasks are easily accomplished even
with slightly sloppy hand drawings of chord diagrams.

Since the combinatorial diagram or the chord diagram completely determines
the combinatorial type, one can derive from either of them any combinatorial
properties of ideal arc systems. In general, for any set of objects on which
$\Homeo(S)$ acts, a \em{combinatorial property} defined on those objects is a
property invariant under the action of $\Homeo(S)$. For example, the \em{polygon
type} of an ideal arc system $\gamma$ is a combinatorial property: this is the
sequence $(i_3,i_4,\ldots)$ where $i_n$ is the number of $n$-gons in $\gamma$. To
determine the polygon type from the combinatorial diagram $(\Ends(\gamma),
\OppEnd, \Succ)$, define a permutation $\Ends(\gamma) \mapright{\Next}
\Ends(\gamma)$ as $\Next = \Succ \composed \OppEnd)$. Then the $n$-cycles of
$\Next$ are in $1-1$ correspondence with the $n$-gons of $\gamma$, for each $n
\ge 3$. For example, in the defect 1 chord diagram of figure 10, the permutation
$\Next$ has cycle structure $\{(1,5,9,13), (2,7,4), (3,8,6), (10,15,12),
(11,16,14) \}$, showing one 4-gon and four 3-gons. Tracing out the boundary of
the 4-gon starting with end 1, and then successively applying $\OppEnd$ and
$\Succ$, we obtain: 
$$1 \mapright{\OppEnd} 4 \mapright{\Succ} 5 \mapright{\OppEnd}	8 \mapright{\Succ} 9
\mapright{\OppEnd} 12 \mapright{\Succ} 13 \mapright{\OppEnd} 16 \mapright{\Succ} 1
$$
as shown in figure 10(a).

\midinsert
\centeredepsfbox{10.ChordDiagramPolygon.eps}
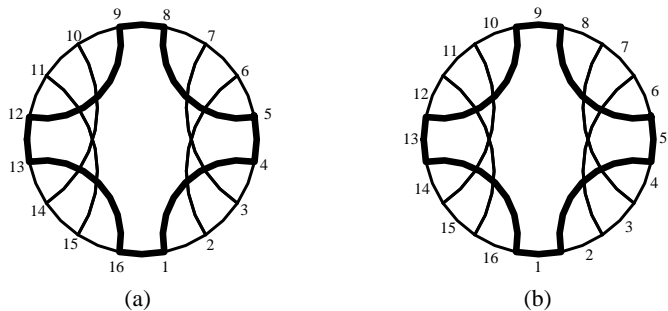
\botcaption{Figure 10} A 4-gon in a defect 1 chord diagram
\endcaption\endinsert

An important concept which plays a central role later on is that of a \em{prong}
of an ideal arc system $\gamma$. Informally, a prong of $\gamma$ is a corner of
a polygon of $\gamma$. Formally, a prong is an ordered pair $(e,e')$ in
$\Ends(\gamma)$ such that $e' = \Succ(e)$. In a chord diagram, a prong is
represented as the circular arc between adjacent ends; we shall call this an
\em{end gap}. Thus, when we represent a polygon in a chord diagram as in figure
10, what is actually drawn are the chords representing the sides of the polygon
and the end gaps representing the prongs of the polygon. If we index each prong
$(e,e')$ using the index of the second end $e'$ in the pair, then each
cycle of $\Next$ lists the prong indices of the corresponding polygon.
Figure 10(b) shows prong indices, making clear the correspondence
between the 4-gon and the cycle $(1,5,9,13)$.

Here is an exercise: prove that given a set $\Ends$ of size $12g-6$, if $\Succ$
is a cyclic permutation and $\OppEnd$ is a transposition, then the triple
$(\Ends,\OppEnd,\Succ)$ is isomorphic to the combinatorial type of an ideal
triangulation on a surface of genus $g$ if and only if the permutation $\Next =
\Succ \composed \OppEnd$ has a cycle structure consisting solely of 3-cycles.
This idea is used in appendix 2 of \cite{P} to obtain an asymptotic formula,
given in the next section, for the number of combinatorial types of ideal
triangulations on a surface of genus $g$.

\subhead Chord diagrams of ideal triangulations \endsubhead

In this section we shall make several observations about chord diagrams of
ideal triangulations. These observations serve two purposes: they help in
learning to recognize features of chord diagrams; and they can be used to
enumerate the chord diagrams on a surface of genus 2. This enumeration was
first obtained by \cite{Jorgensen, Martineen}. We shall also report on
enumerations for higher genus, and an asymptotic formula.

Suppose $\delta$ is an ideal triangulation. Let $T$ be a triangle of $\delta$,
and let $(e_1,e_2,e_3)$ be the corresponding 3-cycle of $\Next$. Every distinct
ordered triple in $\Ends(\delta)$ is either \em{positive} or
\em{negative}. We say that $(e_1,e_2,e_3)$ is positive if
there is a circular enumeration $\Ends(\delta) = \{f_1,\ldots,f_K\}$, i.e.\ an
enumeration with $\Succ(f_k) = f_{k+1}$ for $k \in \integers / K$, so that if
$e_i = f_{k_i}$ then $k_1 < k_2 < k_3$. 

Figure 11 shows how the two types of triangles appear in a chord diagram: a
positive 3-cycle of $\Next$ yields an \em{untwisted} triangle, and a
negative 3-cycle yields a \em{twisted} triangle. Note that $T$ is
untwisted if and only if its regular neighborhood is homeomorphic to a three
holed sphere, and $T$ is twisted if and only if the regular neighborhood is
homeomorphic to a one holed torus. 

\midinsert
\centeredepsfbox{11.TwistedUntwisted.eps}
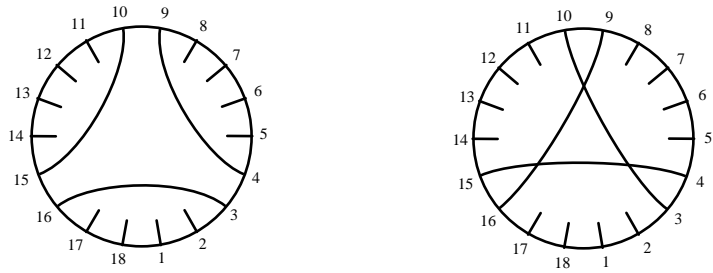
\botcaption{Figure 11} The untwisted triangle corresponds to the
positive or ``increasing'' 3-cycle $(4,10,16)$, while the twisted
triangle corresponds to the negative or ``decreasing'' 3-cycle
$(16,10,4)$
\endcaption\endinsert

\proclaim{Fact} If $S$ is a once-punctured surface of genus $g$, then in any
ideal triangulation $\delta$ the number of twisted triangles is
$2g$ and the number of untwisted triangles is $2g-2$. \endproclaim

\demo{Proof \#1} Construct a single ``base'' example of an ideal triangulation
with $2g$ twisted and $2g-2$ untwisted triangles, observe that the count of
twisted and untwisted triangles is unchanged when doing an elementary move, and
then apply connectivity of the complx $Y$. \qed\enddemo

This proof has the disadvantage that it is not intrinsic to $\delta$: one must
find a path of elementary moves from $\delta$ to the base example. Here is an
intrinsic proof:

\demo{Proof \#2} Let $d$ be a disc neighborhood of $p$ chosen so that $\delta
\intersect d$ is a union of radii of $d$. Let $\hat S = S - \interior(d)$. We
shall use $\delta$ to put a piecewise Euclidean metric on $\hat S$. This metric
will have concentrated negative curvature corresponding to the twisted
triangles. By applying the Gauss-Bonnet theorem we obtain a count of the number
of twisted triangles. To set up the metric requires some alterations on $\delta$.

Let $\delta_1$ be obtained from $\delta$ by replacing each ideal
arc $h$ of $\delta$ with two copies of $h$ bounding a bigon, still intersecting
$d$ in radii. Now consider a triangle $T$ of $\delta_1$. If $T$ is twisted, then
alter $T$ near each prong by taking two half-arcs incident to that prong
extending slightly beyond $d$, and pinching their ends together, as shown in
figure A1; make sure that the half-arcs are left unchanged in $d$, still
intersecting $d$ in radii. The triangle $T$ is divided into four regions: a
pinched triangle, and three pinched prongs. Making this alteration for each
twisted triangle of $\delta_1$, the resulting collection of ideal
arcs is denoted $\delta_2$. Now let $\hat\delta$ be obtained by intersecting
$\delta_2$ with $\hat S$. This has the effect of truncating each untwisted
triangle, each bigon, and each pinched prong; pinched triangles are left intact.

\midinsert
\centeredepsfbox{A1.PinchingProngs.eps}
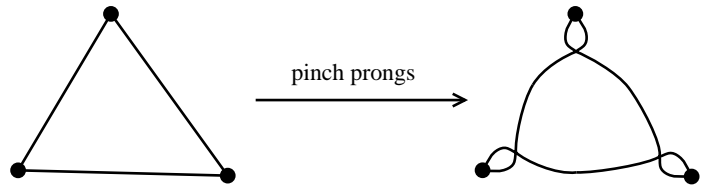
\botcaption{Figure A1} Pinch the prongs of each twisted triangle
\endcaption\endinsert

Consider the chord diagram $D$ of $\delta$. Let $\hat D$ be obtained from $D$
by doubling each chord, replacing it with two parallel chords, then
straightening all chords to become Euclidean segments instead of hyperbolic
lines (see figure A2).

\midinsert
\centeredepsfbox{A2.Doubling.eps}
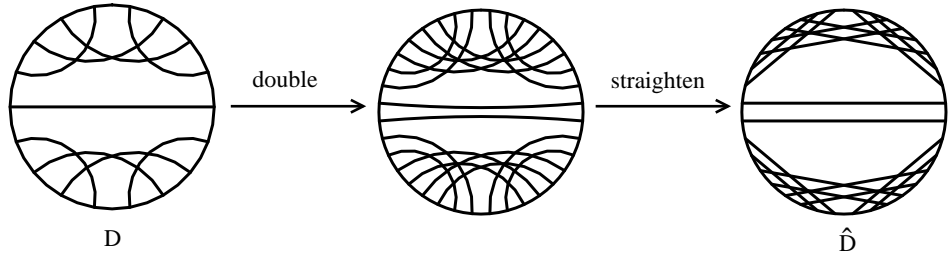
\botcaption{Figure A2} To obtain $\hat D$, double each chord of $D$ then
straighten
\endcaption\endinsert

We regard $\hat D$ as lying in the Euclidean plane $\Euclidean^2$, and we
construct a map $f \from \hat S \to \Euclidean^2$ whose picture is given by
$\hat D$, as follows. The boundary of $\hat S$ goes to the boundary circle of
$\hat D$. Each truncated arc of $\hat \delta$ goes to the corresponding chord of
$\hat D$. Each component of $\hat S - \hat \delta$ is either a truncated
untwisted triangle, truncated bigon, truncated pinched prong, or pinched
triangle; for each of these regions the boundary is already mapped to a simple
closed curve in $\Euclidean^2$, and there is an extension of $f$ to an embedding
of the region, as shown in figure A3. 

Notice that for each pinch point $x$, the map $f$ creates a ``pleat'' at $x$;
see figure A4. Another way to say this is that each twisted triangle is
``twisted'' by the map $f$, at each pinch point of the triangle.

\midinsert
\centeredepsfbox{A3.EmbeddingRegions.eps}
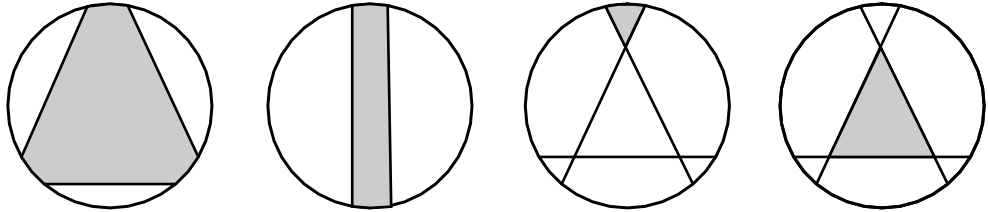
\botcaption{Figure A3} Embedding a truncated untwisted triangle,
truncated bigon, truncated pinched prong, and pinched triangle into $\hat D$
\endcaption\endinsert

By pulling back the Euclidean metric from $\Euclidean^2$ to $\hat D$, we obtain
a piecewise Euclidean metric on $S$. The boundary has total geodesic curvature
$2\pi$. Consider a pinch point $x$. We must compute the Euclidean cone angle
$\theta_x$. The pinched triangle incident to $x$ has a certain interior angle
$\alpha_x$, and as figure A4 shows we have $\theta_x = 2\pi + 2\alpha_x$.
Therefore at $x$ there is an angle defect of $2\pi - \theta_x = -2\alpha_x$. 

\midinsert
\centeredepsfbox{A4.ComputingConeAngle.eps}
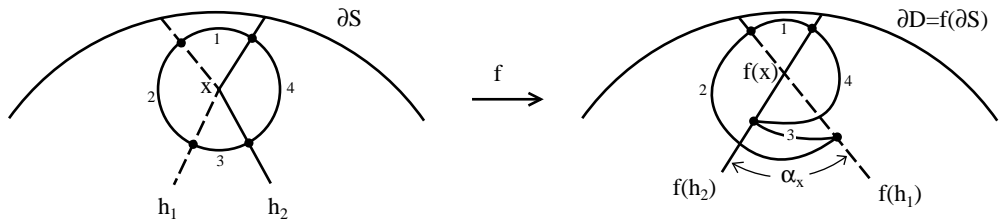
\botcaption{Figure A4} At a pinch point $x$, if the pinched triangle has an
interior angle
$\alpha_x$, then the cone angle at $x$ is $2\pi + 2\alpha_x$.
\endcaption\endinsert

A pinched triangle with vertices $x,y,z$ therefore contributes an angle defect of
$-2(\alpha_x + \alpha_y + \alpha_z)$. But this equals $-4\pi$, since $\alpha_x,
\alpha_y, \alpha_z$ are the interior angles of a Euclidean triangle. Therefore,
if $K$ is the number of twisted triangles, then by the Gauss-Bonnet theorem we
have
$$2\pi - 4\pi K = 2\pi\chi(\hat S) = 2\pi(1 - 2g)
$$
so $K = 2g$.
\qed\enddemo

The twisted and untwisted triangles in an ideal triangulation arrange
themselves into several larger structures, whose visualization helps in
recognizing a chord diagram.

If two untwisted triangles share a side, then they cannot share any other side,
and their union forms an ``untwisted 4-gon''; the chord diagrams of ideal
triangulations in figures 9(a,c) each have an untwisted 4-gon. Continuing
inductively, if an untwisted $n$-gon shares a side with an untwisted triangle,
then they cannot share another side, and their union forms an untwisted
$n+1$-gon. Maximal untwisted polygons are called \em{untwisted islands}. Figure
12(a) shows a genus 2 chord diagram whose two untwisted islands are both
triangles. As an exercise, check that this chord diagram is obtained from figure
9(c) by an elementary move on the arc with ends labelled 1 and 4 (hint: see
figure 24(a) ). Figure 12(b) shows a genus 3 chord diagram with one untwisted
island, a hexagon.

\midinsert
\centeredepsfbox{12.UntwistedPolygons.eps}
\botcaption{Figure 12} Untwisted polygons
\endcaption\endinsert

On the other hand, two twisted triangles can share either one, two, or all
three sides. If they share all three sides, then they close up to form a torus,
with the chord diagram shown in figure 13. Thus, on a higher genus surface a
pair of twisted triangles can share at most two sides. If two twisted triangles
share two sides, then they form a \em{1-handle piece}. The triangulations in
figures 9(a) and (c) each have two 1-handle pieces, figure 12(a) has one, and
figure 12(b) has three.

\midinsert
\centeredepsfbox{13.TorusTriang.eps}
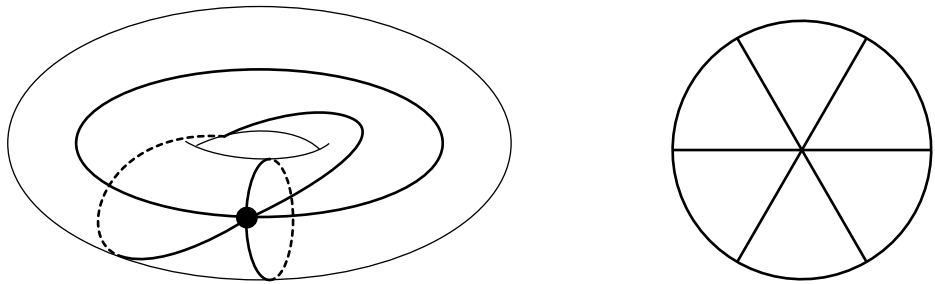
\botcaption{Figure 13} Two twisted triangles sharing three sides
form a torus.
\endcaption\endinsert

If two twisted triangles share only one side, then they form a twisted 4-gon,
see, for example, figure 14(a) which shows the same triangulation as 12(a).
Continuing inductively, if $n \ge 4$ and a twisted $n$-gon shares a side with a
twisted triangle, then they cannot share another side, and their union forms a
twisted $(n+1)$-gon. A maximal twisted polygon is called a \em{twisted island}.
Figure 14(b) shows a genus 2 chord diagram with a twisted hexagon; this
triangulation comes from figure 9(c) by doing two elementary moves one after
another, first on the arc with ends $1,4$ and then on the arc with ends $10,13$.

\midinsert
\centeredepsfbox{14.TwistedPolygons.eps}
\botcaption{Figure 14} Twisted polygons
\endcaption\endinsert

Thus, any ideal triangulation can be decomposed into untwisted islands, twisted
islands, and 1-handle pieces. We may therefore enumerate chord diagrams by the
``island'' method, as follows. First choose a partition of the $2g-2$ untwisted
triangles into islands. Then choose how the prongs of these islands interleave
in the circular ordering. This choice determines the twisted islands and the
number of 1-handle pieces. Now choose triangulations of the twisted and
untwisted islands, using the enumeration by Catalan numbers. 

Now we use the island method to enumerate chord diagrams for low genus
surfaces. A ``connectivity'' proof is given later, using connectivity of
$Y$.

Suppose first that $S$ has genus 1. Then there are no untwisted triangles and 
two twisted triangles $T_1,T_2$. The sides of $T_1$ and $T_2$ must be
glued in 1-1 correspondence, and the prongs must interleave on the chord
diagram. Thus, the chord diagram is forced to be the one shown in figure 13.
This shows that all ideal triangulations of a once-punctured torus are
combinatorially equivalent, and the automorphism group of each one is cyclic of
order 6.

Now suppose $S$ has genus 2. An ideal triangulation has two untwisted and four
twisted triangles. The untwisted triangles can form either a 4-gon island or
two triangle islands.

Suppose first that there is a 4-gon island, which contains a diagonal arc
separating it into two triangles. The four sides of this island must bound two
1-handle pieces, which can arrange themselves in one of three ways as shown in
figure 15. The 1-handle pieces may be parallel to the diagonal arc as in $T_1$;
they may cross the diagonal arc but not cross each other as in $T_2$; or they
may cross the diagonal arc \em{and} each other as in $T_3$. The orders of the
automorphism groups are also shown in figure 15. Note that the unoriented
automorphism groups are dihedral groups of twice the size.

\midinsert
\centeredepsfbox{15.Two1-handleChordDiagrams.eps}
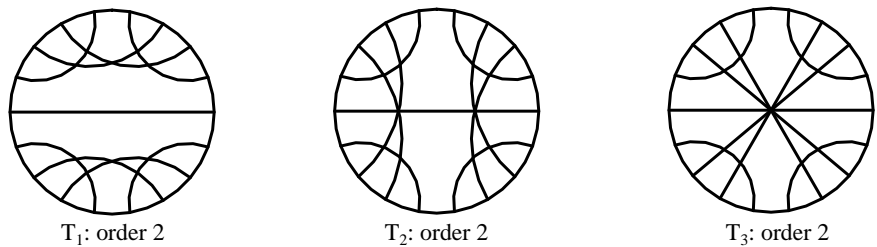
\botcaption{Figure 15} Genus 2 chord diagrams with two 1-handle
pieces
\endcaption\endinsert

Now suppose there are two triangle islands. The prongs of these islands may
interleave in one of two ways, as shown in figure 16.

\midinsert
\centeredepsfbox{16.UntwistedIslands.eps}
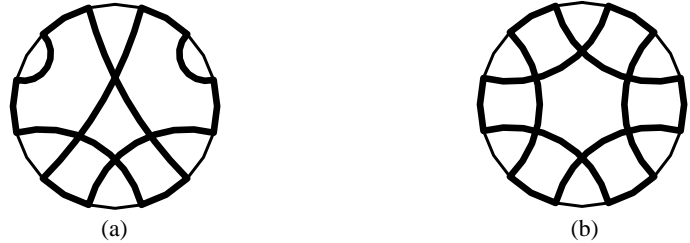
\botcaption{Figure 16} Untwisted triangle islands in genus 2
\endcaption\endinsert

In figure 16(a), there must be one 1-handle piece and one twisted 4-gon island.
There are two ways to insert a diagonal in the twisted 4-gon, yielding the two
chord diagrams $T_4$ and $T_5$ shown in figure 17; these are the only chord
diagrams in genus 2 with one 1-handle piece. Their automorphism groups have
order 1; their unoriented automorphism groups have order 2.

\midinsert
\centeredepsfbox{17.One1-handleChordDiagrams.eps}
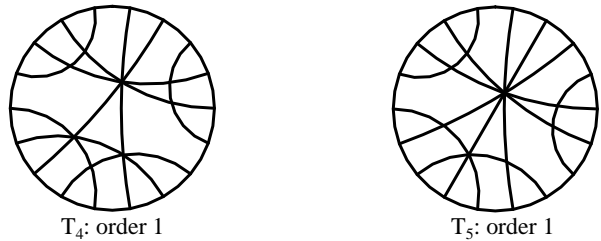
\botcaption{Figure 17} Two genus 2 chord diagrams with one 1-handle piece
\endcaption\endinsert

In figure 16(b), there is a single twisted island, a hexagon. Up to rotation
there are four ways to triangulate this hexagon, yielding the four chord
diagrams in figure 18, the only chord diagrams in genus 2 with no 1-handle
pieces. Both $T_6$ and $T_7$, which are orientation reversals of each other,
have oriented and unoriented automorphism groups cyclic of order 2. The diagram
$T_8$ has trivial automorphism group, and the unoriented automorphism group is
dihedral of order two. The diagram $T_9$ has automorphism group cyclic of order
3, and unoriented automorphism group dihedral of order 6.

\midinsert
\centeredepsfbox{18.No1-handleChordDiagrams.eps}
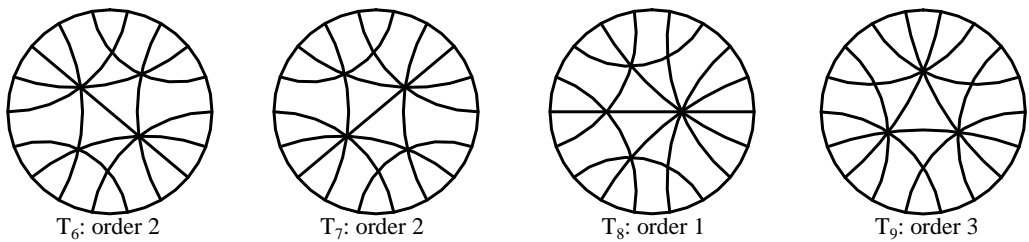
\botcaption{Figure 18} Four genus 2 chord diagrams with no 1-handle
pieces
\endcaption\endinsert

To summarize, figures 15,17 and 18 show the nine combinatorial types of ideal
triangulations on a once-punctured surface of genus 2.

The island method used to obtain this enumeration is rather inefficient,
although it is good for learning to recognize chord diagrams. Despite this
inefficiency, when I was young and energetic I used the island method to
enumerate the chord diagrams on a once-punctured surface of genus 3. I then
wrote a computer program implementing the connectivity method (explained later),
obtaining 1726 combinatorial types. This did not accord exactly with the island
method, so I went through and found some errors, correcting the result of the
island method, still not obtaining the same answer. After iterating this process
a few time, I obtained a count of 1726 chord diagrams, and quit.

An asymptotic formula for $n_g$, the number of distinct chord diagrams in genus
$g$, is given in appendix B of \cite{P}:
$$ n_g \sim {(2g)! \over 6g-3} \left( {e \over g} \right)^{-2g}
$$
where $x(g) \sim y(g)$ means that $x(g)/y(g) - 1 = O(1/g)$.

Figure 19 summarizes what I know about the numbers of combinatorial types of
ideal triangulations, compared to the above asymptotic formula. The purpose of
this table is to drive home the point that one would not want to enumerate chord
diagrams, or any objects derived from them such as the states of the automatic
structure, and store them all in one place, if it were not absolutely necessary
for computation. On the other hand, methods for generating the objects as
needed are very useful; this is how our algorithm for the word problem works.

\midinsert
\topcaption{Figure 19} A table of chord diagrams \endcaption
$$\vbox{
\settabs\+Genus\quad&1726\quad&1,081,820&\cr 
\+Genus&$n_g$&${(2g)! \over 6g-3} \left( {g \over e} \right)^{2g}$\cr
\+1&1&3\cr
\+2&9&63\cr
\+3&1726&5551\cr
\+4&?&1,081,820\cr
}$$
\endinsert

\subhead Labelling ideal triangulations: the zero skeleton of $X$ \endsubhead

In this section we introduce the machinery needed to define the zero-skeleton
of $\tilde X$ and of $X$ itself. The point is this: we already have an action
of $\mcg$ on the zero-skeleton of $Y$, but that action has some non-trivial point
stabilizers. In order to obtain the zero-skeleton of $\tilde X$ we need an
action of $\mcg$ with trivial point stabilizers. Thus, we must somehow break
the symmetries of an ideal triangulation, by labelling it with extra data.

A \em{labelled ideal triangulation} consists formally of an ordered pair
$(\delta,e)$ where $\delta$ is an ideal triangulation and $e$ is an arc end of
$\delta$. The mapping class group acts on isotopy classes of labelled ideal
triangulations, with trivial stabilizers. In later sections, we will often
suppress the labelling $e$, and speak of ``a labelled ideal triangulation
$\delta$''. For now, we stick with the formal notation $(\delta,e)$.

We now define the zero skeleton of $\tilde X$ to be the set of isotopy classes of
labelled ideal triangulations on $S$. The zero skeleton of $X$ is therefore the
set of combinatorial types of labelled ideal triangulations. The combinatorial
type of a labelled ideal triangulation $(\delta,e)$ is described by a
\em{labelled chord diagram}, obtained from the chord diagram for $\delta$ by
drawing a solid dot at the chord end corresponding to $e$. An example is shown
in figure 20, where the labelled end is represented by a shaded half-arc. 

\midinsert
\centeredepsfbox{20.LabelledIdealTriang.eps}
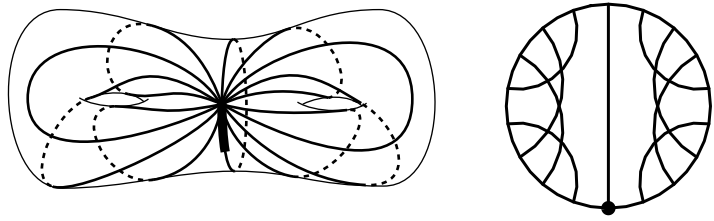
\botcaption{Figure 20} A labelled ideal triangulation and the corresponding
labelled chord diagram
\endcaption\endinsert

Now we enumerate the 0-cells of $X$ in low dimensions, using the fact that a
chord diagram, with $12g-6$ ends and automorphism group cyclic of order $k$,
yields $(12g-6)/k$ labelled chord diagrams.

If $S$ has genus 1, the unique chord diagram has six chord ends and the
automorphism group permutes them transitively, so there is a unique chord diagram
of a labelled ideal triangulation. Therefore, $X$ has one 0-cell.

If $S$ has genus 2, each chord diagram has 18 chord ends, so using the orders
of the automorphism groups given in figures 15,17,18 the number of labelled
chord diagrams is
$$3 \cdot 18 + 5 \cdot {18 \over 2} + 1 \cdot {18 \over 3} = 105
$$
so there are 105 0-cells in $X$.

If $S$ has genus 3, since $12g-6=30$ then an upper bound on the number of
0-cells is $30 \cdot 1726 = 51780$, but this number will be strictly smaller
after taking automorphism groups into account. A computer calculation has yielded
$50050$ 0-cells in $X$.

\subhead The construction of $\tilde X$, $X$, and the mapping class groupoid
\endsubhead

It is now possible to give a purely abstract definition of the mapping class
groupoid $\mcgd$, as was done in \cite{M}. Recall that an abstract groupoid is a
category with invertible morphisms. Let $\tilde D$ be the set of isotopy classes
of labelled ideal triangulations on $S$. Then $\mcg$ acts freely on $\tilde D$,
and the diagonal action on $\tilde D \cross \tilde D$ is also free. The objects
of $\mcgd$ are the orbits of the action of $\mcg$ on $\tilde D$, i.e.\ the
combinatorial types of labelled ideal triangulations. The morphisms of $\mcgd$
are the orbits of the diagonal action of $\mcg$ on $\tilde D \cross \tilde D$.
If $\delta$, $\delta'$ are labelled ideal triangulations (suppressing the
labellings), then the orbit of the pair $([\delta], [\delta'])$,
denoted $\{ \delta, \delta' \}$, has as its initial object $\{ \delta
\}$ and as its terminal object $\{ \delta' \}$. The composition rule is as
follows. Given morphisms $m_1 = \{ \delta_1,\delta'_1 \}$ and $m_2 =
\{\delta_2,\delta'_2\}$ such that the terminal object $\{ \delta'_1 \}$ of $m_1$
equals the initial object $\{ \delta_2 \}$ of $m_2$, then $\delta'_1$ and
$\delta_2$ are combinatorially equivalent. Therefore there exists $\Phi \in \mcg$
such that $\Phi[\delta'_1] = [\delta_2]$. Note that $m_2 = \{ \Phi(\delta_2),
\Phi(\delta'_2)\}$, abusing notation. Then $m_1 \composed m_2$ is defined to be
$\{ \delta_1, \Phi(\delta'_2) \}$.

One easily checks that if $\tilde D$ is the 0-skeleton of a simply connected
complex $\tilde X$ on which $\mcg$ acts freely, then the abstract groupoid
constructed above is naturally isomorphic to the edge path groupoid of
the quotient complex $X = \tilde X / \mcg$. We now proceed to the construction
of $\tilde X$.

We have already constructed $\skel(\tilde X,0) = \tilde D$, the set of isotopy
classes of labelled ideal triangulations. Now we construct $\tilde X$, together
with a cellular map $q \from \tilde X \to Y$ which is useful in proving simple
connectivity of $X$. The construction is by ``abstract nonsense''. Each $k$-cell
will come equipped with a ``boundary certificate'', which is a total ordering of
the cells of all dimensions on its boundary, and the boundary certificate
determines the $k$-cell. This convention will allow us to define the action of
$\mcg$ up through the skeleta of $\tilde X$ by induction; we will similarly prove
that the action has trivial cell stabilizers. 

Any unordered pair of 0-cells in $\skel(\tilde X,0)$ to which an edge is
attached will be called a \em{boundary 0-cycle}. Given any boundary 0-cycle,
each of the two possible ordering will be the boundary certificate of some
1-cell. With this convention, to specify the 1-cells of $X$ we need only
specify which vertex pairs are boundary 0-cycles.

Similarly, any edge cycle $\gamma$ in $\skel(\tilde X,1)$ to which a 2-cell is
attached will be called a \em{boundary 1-cycle}. Given any boundary 1-cycle
$\gamma$, if you choose a vertex and an orientation, then you obtain a boundary
certificate by starting with that vertex and reading off the cells encountered
by going around in that order. If $\gamma$ has $k$ vertices then there are $2k$
possible choices, each of which is the boundary certificate of some 2-cell,
and these are the only 2-cells attached to $\gamma$. Again, with this
convention we need only specify what the boundary 1-cycles are, in order to
specify the 2-cells of $X$.

Consider a vertex $[\delta]$ of $Y$; we begin by constructing the part of
$\tilde X$ lying over $[\delta]$, also known as $q^\inverse[\delta]$. We
already know that the 0-cells lying over $[\delta]$ are the isotopy classes
$[\delta,e]$ of labellings of $\delta$. 

Every unordered pair of vertices in $q^\inverse[\delta]$ will be a boundary
0-cycle. Thus, for every ordered pair of ends $e_1,e_2$ of $\delta$ there is a
1-cell with boundary certificate $([\delta,e_1],[\delta,e_2])$. We denote this
1-cell as $[\delta,e_1,e_2]$. The image of this cell downstairs in $X$ is called
a \em{relabelling generator} of $X$. 

Every edge cycle of length 2 or 3 in the 1-skeleton of $q^\inverse[\delta]$ will
be a boundary 1-cycle. To determine a boundary 1-cycle $\gamma$ of length 2,
choose a set of two labels $\{e_1,e_2\}$, so the edges of $\gamma$ are
$[\delta,e_1,e_2]$ and $[\delta,e_2,e_1]$. To determine a boundary 1-cycle
$\gamma$ of length 3, choose a set of three labels $\{e_1,e_2,e_3\}$, and choose
an ordering for each of the sets $\{e_1,e_2\}$, $\{e_2,e_3\}$, $\{e_3,e_1\}$ to
obtain the edges of $\gamma$. The images downstairs in $X$ of these 2-cells
will be called \em{relabelling relators}.

Now consider an elementary move $\delta \to \delta'$, yielding a 1-cell
$[\delta \to \delta']$ of $Y$. We already know that the 0-skeleton of
$q^\inverse[\delta \to \delta']$ is $q^\inverse[\delta] \union
q^\inverse[\delta']$. Each pair of a vertex in $q^\inverse[\delta]$ and a vertex
in $q^\inverse[\delta']$ is a boundary 0-cell. The images downstairs in $X$ of
the attached 1-cells will be called \em{elementary move generators}. 

To determine the 2-cells in $q^\inverse[\delta \to \delta']$, note
that in the 1-skeleton of $q^\inverse[\delta \to \delta']$, each edge cycle
contains an even number of edges mapping to $[\delta\to\delta']$ under $q$;
between two such edges the cycle may wander around for a while in
$q^\inverse[\delta]$ or $q^\inverse[\delta']$. The boundary 1-cycles in
$q^\inverse[\delta \to \delta']$ are the ones which have exactly two edges
mapping to $[\delta \to \delta']$, and which contain at most one edge each in
$q^\inverse[\delta]$ and $q^\inverse[\delta']$. The images downstairs in $X$ of
the attached 2-cells will be called \em{elementary move relabellings}. 

This finishes the description of $q^\inverse(\skel(Y,1))$.

Finally, given a 2-cell $c$ of $Y$, each edge cycle in $\skel(\tilde X,1)$
projecting homeomorphically to $\bd c$ is a boundary 1-cycle, and the attached
2-cell maps homeomorphically to $c$. Depending on the nature of the 2-cell $c$,
the images of these 2-cells in $X$ will be called \em{labelled commutator
relators} or \em{labelled pentagon relators}.

This completes the description of $X$.

The action of $\mcg$ on $\skel(\tilde X,0)$ is already defined, and each 0-cell
has trivial stabilizer. A 1-cell is determined by its boundary certificate,
which is an ordered pair of 0-cells, and clearly the set of such ordered pairs is
invariant under the action of $\mcg$, so the action of $\mcg$ extends over
$\skel(\tilde X,1)$ with trivial 1-cell stabilizers. Again, a 2-cell is
determined by its boundary certificate, which is a sequence of 0-cells and
1-cells, and the set of such sequences is invariant under the action of $\mcg$,
so the action of $\mcg$ extends over 2-cells with trivial cell stabilizers.

The quotient complex $X = \tilde X / \mcg$ is now defined. It is evident that
the map $q \from \tilde X \to Y$ has the path lifting property as well as the
homotopy lifting property for paths, so $\tilde X$ is simply connected. It
follows that $\pi_1(X) \isomorphic \mcg$. Also, for any cell of $Y$ the inverse
image in $\tilde X$ is a finite cell complex, and since $Y$ has finitely many
cell orbits it follows that $X$ is a finite complex. We can
now define the mapping class groupoid $\mcgd$ as the edge path groupoid
of $X$. Given an edge path $w$ in $X$, the corresponding homotopy class is
denoted $\overline w \in \mcgd$.

We now have a finite presentation for $\mcgd$, with the edges of $X$ as
generators and the 2-cells of $X$ as relators. There are two types of edges:
relabelling generators and elementary move generators. There are several types of
relators: relabelling relators, elementary move relabellings, and labelled
commutator and pentagon relators.

In the next two sections we whittle down the elementary move generators to a
smaller subset called the ``labelled elementary move generators'', which
together with the relabelling generators will still generate $\mcgd$ (this is
the generating set used in \cite{M}). In order to understand labelled elementary
move generators, we first initiate a study of chord diagrams of elementary moves.

\subhead Chord diagrams of elementary moves
\endsubhead

Consider an elementary move $\delta \to \delta'$ performed on the ideal arc $h$
of $\delta$, with opposite diagonal $h'$ in $\delta'$, and with support $Q$. Let
$D$ be the chord diagram of $\delta$; by abuse of notation we use $h$ to stand
for the chord representing the ideal arc $h$, and this chord will be shaded in
the diagrams. Now we show how, using $D$ and $h$ as input, we may compute the
chord diagram $D'$ of $\delta'$.

Look at the two triangles adjacent to $h$. There are several possibilities for
these two triangles: both untwisted; one untwisted and one twisted, also known
as \em{mixed}; or both twisted, with either one, two, or three side pair
idenfications. For each of these five cases, figure 21 shows how the elementary
move appears ``locally'' in a chord diagram, i.e.\ the figure shows the
support of the elementary move, with missing chord ends indicated by a
twiddle "$\sim$". The five cases can also be enumerated according to the chord
diagram of the quadrilateral $Q$. We now go through the cases one by one.

\midinsert
\centeredepsfbox{21.LocalElemMoves.eps}
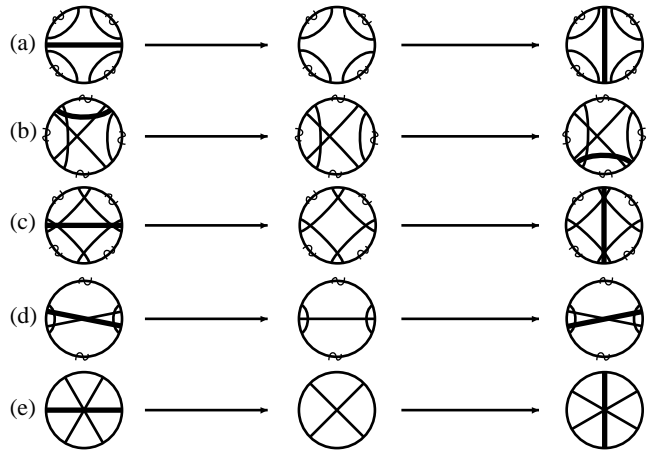
\botcaption{Figure 21} Chord diagrams of the support of an elementary
move. Missing chord ends are indicated with a $\sim$.
\endcaption\endinsert

First we dispose of the case where both triangles are twisted and there are
three side pair identifications. This occurs only on the punctured torus, and
is given in figure 21(e). The quadrilateral $Q$ has two side pair
identifications. Any elementary move on the punctured torus has this
chord diagram. Thus, in genus 1 there is only one orbit of edges of $Y$ under the
action of $\mcg$. 

With both triangles untwisted, then $Q$ is an untwisted 4-gon, and we have an
\em{untwisted-untwisted} elementary move (figure 21(a)). With both triangles
twisted and one side pair identification, then $Q$ is a twisted 4-gon, and we
have a \em{twisted-twisted} elementary move (figure 21(c)). These two cases are
the easiest to visualize. An untwisted-untwisted example is shown in figure
22(a); this is the chord diagram of the elementary move from figure 5(a). A
twisted-twisted example is given in figure 22(b). In these chord diagrams the
intervening defect 1 chord diagram is shown; after this section we will not
usually show this. One general feature to note is that when one chord is removed
and another inserted in a chord diagram, the chord endpoints should be
repositioned so that they are evenly spaced; space will have to be contracted
near the removed chord ends, and it will be expanded near the inserted chord
ends.

\midinsert
\centeredepsfbox{22.Untw-Untw_Tw-Tw.eps}
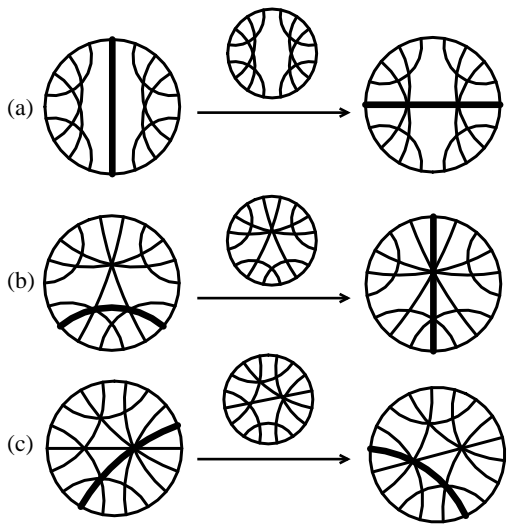
\botcaption{Figure 22} An untwisted-untwisted and two twisted-twisted
elementary moves
\endcaption\endinsert

With both triangles twisted and two side pair identifications, the triangles form
a 1-handle piece and $h$ is one of the two interior arcs of the 1-handle piece.
The support quadrilateral $Q$ has one side identification (figure 21(d)). An
example is given in figure 23, which is the chord diagram of the elementary move
from figure 5(b). Notice that after the elementary move, a 1-handle piece forms
again, and outside the 1-handle piece the chord diagram is unchanged, so
$\delta$ and $\delta'$ have the same chord diagram. This means that there is a
mapping class $\Phi \in \mcg$ such that $\Phi[\delta] = [\delta']$. This mapping
class may always be taken to be a Dehn twist about the core of the 1-handle, as
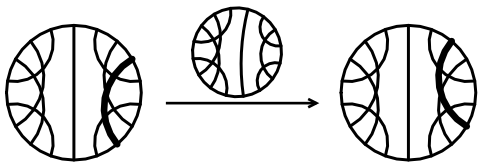
noted earlier for figure 5(b), and we call this a \em{Dehn twist} elementary move
(if the chord diagram has symmetries, as in figure 23, we can also post-multiply
the Dehn twist $\Phi$ by any mapping class which stabilizes $[\delta']$). The
phenomenon of Dehn twist elementary moves is the tip of a big iceberg; in
section V we make a general study of sequences of elementary moves
representing Dehn twists, and this is used to obtain an automatic structure and
prove quadratic computation time of our algorithm for the word problem in
$\mcg$. 

\midinsert
\centeredepsfbox{23.DehnTwistElemMove.eps}
\botcaption{Figure 23} A Dehn twist elementary move
\endcaption\endinsert

In the mixed case, $Q$ is neither a twisted nor untwisted 4-gon and we say that
$Q$ is a \em{mixed} 4-gon (figure 21(b)). A \em{mixed} elementary move is usually
the most difficult to visualize. Several examples of mixed elementary moves are
shown in figure 24 and 25. Figure 24 shows examples where $h$ lies on the
boundary of a 1-handle piece, and figure 25 shows examples where $h$ lies on the
boundary of a twisted island; the opposite diagonal may be of one or the other
type. 

Exercise: Figures 24(c) and 25(b) are interesting because they return to the
same chord diagram, hence $\delta$ and $\delta'$ differ by a mapping class. What
are these mapping classes? To be more precise, how do they fit into Thurston's
classification scheme of finite order, reducible, or pseudo-Anosov?

\midinsert
\centeredepsfbox{24.MixedElemMoves1.eps}
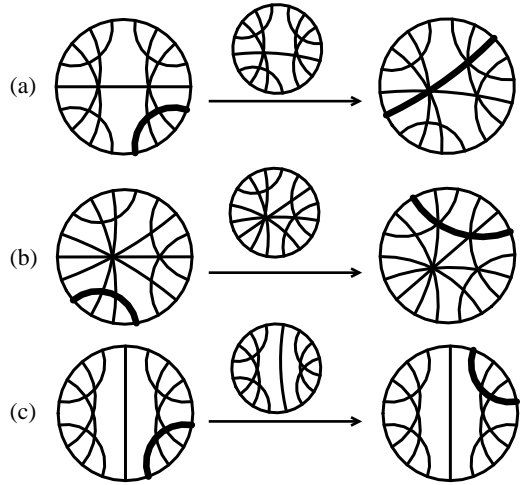
\botcaption{Figure 24} Mixed elementary moves where $h$ bounds an
untwisted island and a 1-handle piece
\endcaption\endinsert

\midinsert
\centeredepsfbox{25.MixedElemMoves2.eps}
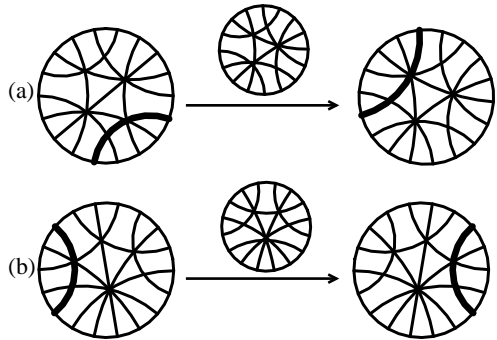
\botcaption{Figure 25} Mixed elementary moves where $h$ bounds an untwisted
island and a twisted island
\endcaption\endinsert

Exercise: Notice that among figures 22--25, we have managed to produce paths of
elementary moves from $T_1$ to $T_2$, $T_3$, $T_4$, $T_7$, $T_8$ and $T_9$. For
example, follow the path
$$ T_1 \mapright{22(a)} T_2 \mapright{24(a)} T_3 \mapright{22(b)} T_4 
\mapright{24(b)} T_8 \mapright{22(c)} T_7 \mapright{25(a)} T_9
$$
Construct enough chord diagrams of elementary moves to obtain the remaining
chord diagrams $T_5$, $T_6$. Getting $T_5$ is slightly tricky, because there is
only one other chord diagram that $T_5$ may be accessed from by a single
elementary move. 

Exercise: How long is the shortest path from $T_1$ to $T_7$?

Having described in some detail how elementary moves are represented with
chord diagrams, we remark that the computer representation of elementary moves
using combinatorial diagrams is easily implemented. There is a simple
algorithm which takes as input the combinatorial diagram $(\Ends, \OppEnd,
\Succ)$ of $\delta$, together with the cycle of $\OppEnd$ representing $h$, and
outputs the combinatorial diagram of $\delta'$ where $\delta \mapright{h}
\delta'$.

Once we know how to generate elementary moves on combinatorial diagrams, there
is a simple algorithm for enumerating combinatorial types of ideal
triangulations. The quotient of $Y$ under the action of $\mcg$ can be regarded
as a finite, connected 1-complex, with a 0-cell for every combinatorial type
of ideal triangulation, and with a 1-cell for every combinatorial type of
defect 1 arc system (more properly, the quotient complex should be thought of
as an ``orbi-complex'' in the sense of Haefliger; for instance, if an edge in
$Y$ has an orientation reversing stabilizing element then its image in the
quotient should be regarded as a half-edge with a ``mirrored endpoint''). It is
then easy to write an algorithm for constructing this 1-complex, say using a
breadth first search: construct an initial chord diagram; initialize a queue
with one entry for each chord of the initial chord diagram; now process the
queue inductively, taking the first entry of the queue and performing the
indicated elementary move; check if the new chord diagram has already been
found, and if not add an entry to the end of the queue for each chord of the new
chord diagram. 

The result of this algorithm for a surface of genus 2 is shown in figure 26.
This is what we call the ``Connectivity Proof'' for the enumeration of the nine
chord diagrams in genus 2. Each chord is labelled, with chords in a diagram
having the same label when there is an automorphism carrying one chord to the
other. Elementary moves are also given labels; for instance, the elementary move
between $T_2$ and $T_3$ labelled $1$--$5$ means that chord 1 is removed from
$T_2$ and chord 5 is inserted in $T_3$ (or vice versa). There is also a single
mirrored 1-cell, which is drawn unmirrored as the edge labelled $7$--$7$ going
from $T_5$ to itself; note that this edge corresponds to a defect 1 chord
diagram with an order 4 symmetry group that rotates the quadrilateral by 1/4.

\midinsert
\centeredepsfbox{26.ChordDiagramGraph.eps}
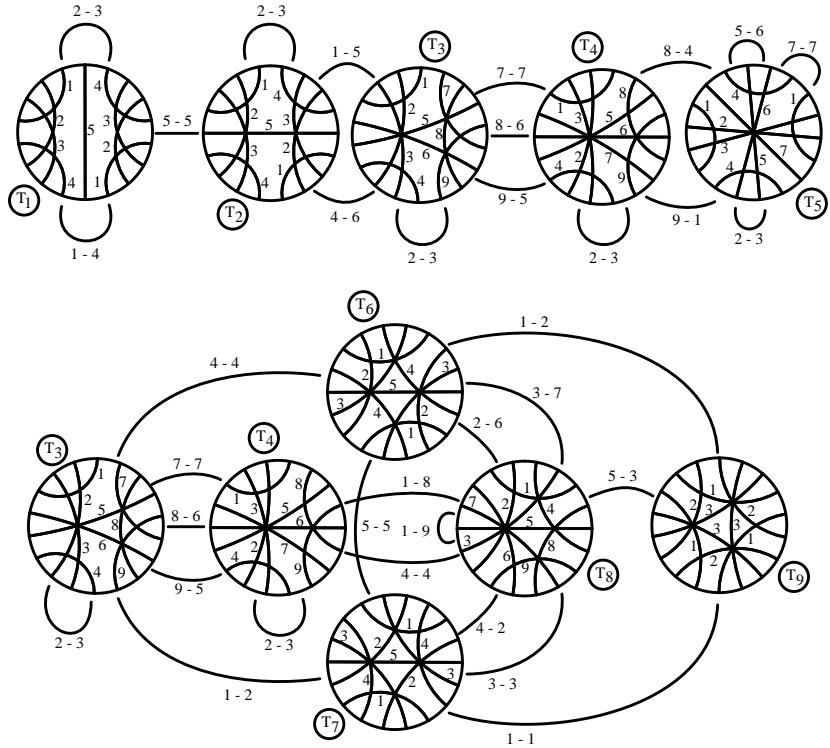
\botcaption{Figure 26} The chord diagrams of ideal triangulations and
elementary moves on a once-punctured surface of genus 2. The top part of the
diagram overlaps with the bottom part in triangulations $T_3$ and $T_4$.
\endcaption\endinsert

\subhead The end map of an elementary move \endsubhead

For any elementary move $\delta \mapright{h} \delta'$ where $h'$ is the
opposite diagonal of $h$, the ideal arc systems $\delta - h$ and $\delta' -
h'$ are isotopic, hence by the lemma \tit{Rigidity of ends} we have a
well-defined bijection $\Ends(\delta - h) \to \Ends(\delta'-h')$ called the
\em{end map}. This map may be implemented in chord diagrams as follows. Suppose
that $D$ is the chord diagram for $\delta$, and let the chord corresponding to
$h$ also be denoted $h$. Index the chord ends of $D$ except for $h$, starting at
an arbitrary end with 1 and increasing in counter-clockwise order, skipping over
the ends of $h$. Now when the chord $h$ is erased and the opposite chord $h'$ is
inserted, resulting in the chord diagram $D'$, we have an indexing of the chord
ends of $D$ except for $h$. This indexing gives the end map, a bijection
between chord ends of $D$ except for $h$ and chord ends of $D'$ except for
$h'$. Examples are shown in figures 27 and 28. The end map plays an important
role in what follows.

\midinsert
\centeredepsfbox{27.EndMap.eps}
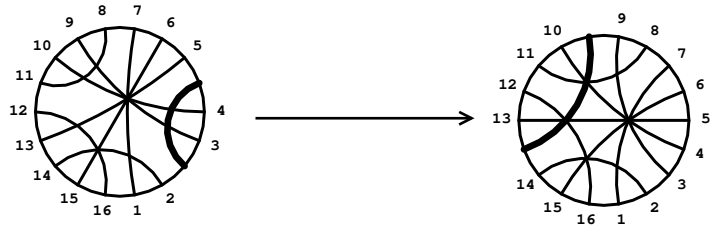
\botcaption{Figure 27} Chord ends with the same index correspond under the end
map
\endcaption\endinsert

\subhead Labelling elementary moves \endsubhead
Consider a labelled ideal triangulation $(\delta,e)$, so that $e$ is an end of
an arc $g$ of $\delta$. Consider also an elementary move $\delta \to \delta'$
performed on an arc $h$ of $\delta$. We adopt the following convention for
determining a labelling $e'$ of $\delta'$. If $g \ne h$ then set $e' = e$;
whereas if $g = h$, let $e'$ be the predecessor of $e$ in $\Ends(\delta)$. In
either case, $e'$ lies on an arc of $\delta$ which is also an arc of $\delta'$,
hence the labelled ideal triangulation $(\delta',e')$ is defined. The complex
$\tilde X$ has a 1-cell with boundary certificate $([\delta,e], [\delta',e'])$.
The image of this 1-cell downstairs in $X$ is denoted $\{ \delta,e \}
\mapright{h} \{ \delta',e' \}$, and is called a \em{labelled elementary move};
in this notation, $h$ should be regarded as a chord in the chord diagram for $\{
\delta,e \}$. 

In order to understand chord diagrams of labelled elementary moves, suppose that
$D \to D'$ is the chord diagram of the labelled elementary move $\delta \to
\delta'$. Let $h$ be the removed chord of $D$, and let $h'$ be the inserted
chord of $D'$. Suppose that $e$ is the labelled chord end in $D$. If $e$ is
not an end of $h$, then the labelled chord end $e'$ of $D'$ is just the image of
$e$ under the end map. On the other hand, if $e$ is an end of $h$, then $e'$ is
the image under the end map of the predecessor of $e$, obtained by rotating $e$
one notch clockwise. Examples are given in figure 27. 

\midinsert
\centeredepsfbox{28.LabelledElemMoves.eps}
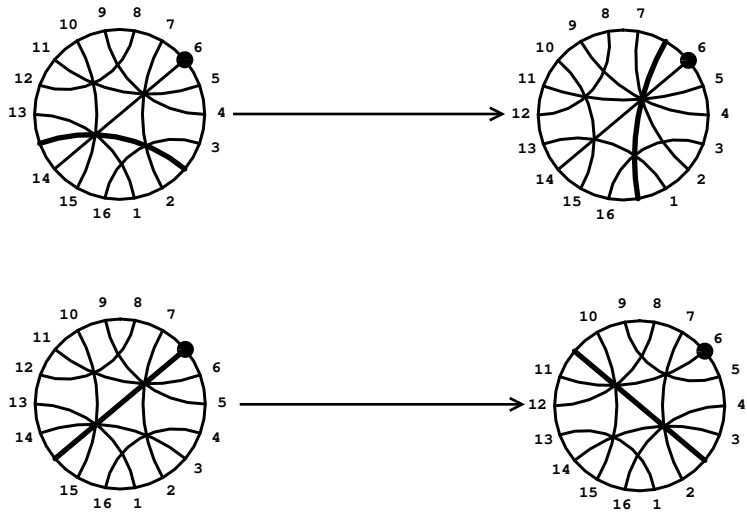
\botcaption{Figure 28} Some labelled elementary moves
\endcaption\endinsert

\subhead Relabelling moves \endsubhead

We have already defined relabelling generators: given an ideal triangulation
$\delta$ and two distinct arc ends $e_1,e_2$ of $\delta$, there is an edge
$[\delta,e_1,e_2]$ in $X$ with boundary certificate $[\delta,e_1],
[\delta,e_2]$. The image of this edge downstairs in $X$ is denoted $\{
\delta, e_1, e_2 \}$ and is called a \em{relabelling generator}; it points from
the vertex $\{\delta,e_1 \}$ to $\{ \delta,e_2 \}$. Since the chord diagram $D$
for $\delta$ has $12g-6$ chord ends arrayed in circular order, then we can write
$e_2 =\Succ^r(e_1)$ for a unique $r \in \integers / 12g-6$. Then we say that
$e_2$ is obtained from $e_1$ by \em{rotating $r$ notches}, and we denote the
relabelling generator as
$$ \{ \delta,e_1 \} \Relabelling{r} \{ \delta,e_2 \}
$$
Figure 29 gives an example of a relabelling generator.

\midinsert
\centeredepsfbox{29.RelabellingMove.eps}
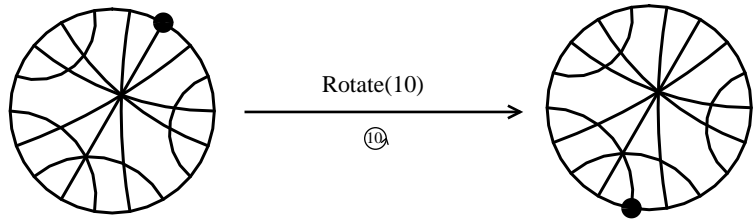
\botcaption{Figure 29} A relabelling generator (both notations will be used)
\endcaption\endinsert

\proclaim{Proposition} The groupoid $\mcgd$ is generated by labelled elementary
moves and relabelling generators. In fact, every groupoid element may be
represented by a string of labelled elementary moves followed by a single
relabelling generator.
\endproclaim

\demo{Proof} A few observations about relators make this obvious. First, by
using relabelling relators, any consecutive sequence of relabelling generators
may be replaced by a single relabelling generator. Second, by using elementary
move relabellings, any relabelling generator followed by an elementary move
generator may be replaced by an elementary move generator, followed by at most
one relabelling generator. Third, also by using elementary move relabellings,
any elementary move generator $g_1$, may be replaced by a labelled elementary
move $g_2$ followed by a relabelling generator. Given an arbitrary word,
conglomerate all initial relabelling generators into one and push it past the
first elementary move generator, then replace that by a labelled elementary move
if necessary; now repeat the procedure starting with the next block of
relabelling generators. \qed\enddemo

\head II. Asynchronous normal forms \endhead

In this section we describe normal forms for elements of $\mcgd$. The normal
forms will be defined over the alphabet $\cal A_0$ consisting of all labelled
elementary moves and relabelling generators. We shall construct a finite
automaton $\cal M_0$ defined over $\cal A_0$, and the normal forms will be the
language $\cal L_0$ accepted by this automaton. First a quick review of finite
automata over groupoid generators. 

The automaton $\cal M_0$ will be a directed graph, whose vertices are called
\em{states} and whose directed edges are called \em{arrows}. Each arrow will be
named with an element of $\cal A_0$. There will be a cellular map $p \from \cal
M_0\to X$; each state $s$ goes to a vertex $ps=D$, and each arrow going out of
$s$ is named by a generator going out of $D$ which is identified with the image
of that arrow under $p$. One state of $\cal M_0$ is specified as the \em{start
state}, and it will map to the base vertex of $X$. Some subset of states are
specified as the \em{accept states}. An \em{accept path} is a directed path from
the start state to an accept state, and by reading off the names of edges along
that path we obtain a word in $\cal A_0$, called an \em{accepted word}. The
language $\cal L_0$ will be the set of all accepted words. Note that each
accepted word represents an edge path in $X$ starting at the base vertex.

The proof that $\cal L_0$ represents each element uniquely is given in
section II.5 of \cite{M}, culminating in the proposition \tit{Normal forms are
regular} (the language $\cal L_0$ is defined differently in \cite{M}, but from
the proof of \tit{Normal forms are regular} the two definitions clearly give
the same language). Section III describes an algorithm for computing the
normal form representing a given word in the generators $\cal A_0$; from this
description it is straightforward to show that $\cal L_0$ satisfies the
asynchronous fellow traveller property, hence is an asynchronous automatic
structure for $\mcgd$. 

Remark 1: The words in $\cal L_0$ will each have at most one relabelling
generator, and it is always the last letter. This is different from the
convention adopted in the original version of \cite{M}, where the relabelling
generator comes first. This change makes no difference in proving the
asynchronous fellow traveller property, but it does simplify the description of
an algorithm for computing normal forms.

Remark 2: From now on, we usually suppress the label in our notation for a
labelled ideal triangulation, writing $\delta$ instead of the more formal
$(\delta,e)$.

\subhead The states of $\cal M_0$ \endsubhead

Recall that normal forms for an asynchronous automatic structure on a groupoid
must all start at some chosen base vertex. Once and for all, pick some labelled
ideal triangulation $\delta_B$ as a base vertex of $\tilde X$, and the
combinatorial type $\{\delta_B\}$ will be the base vertex of $X$. This choice
is quite arbitrary, but for the figures to come we choose the base vertex given
in figure 30. This pattern may be generalized to any genus $g$: take a chord
diagram with $g$ 1-handle pieces, none crossing any other, and then take the
``fan triangulation'' of the resulting untwisted $2g$-gon, putting the labelled
end just clockwise of the base of the fan.

\midinsert
\centeredepsfbox{30.BaseVertex.eps}
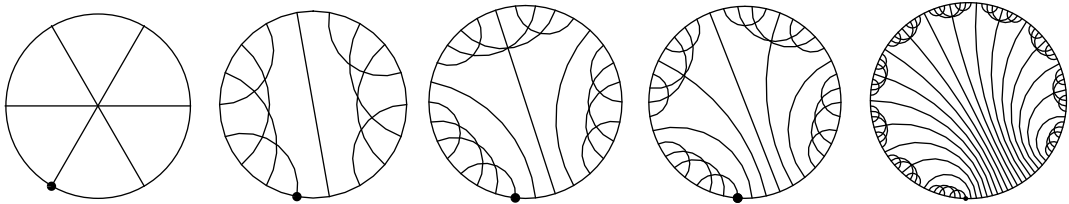
\botcaption{Figure 30} Our choice for a base vertex, in genus 1,2,3,4 and 10
\endcaption\endinsert

A finite deterministic automaton can be thought of as a machine with a fixed,
finite number of states and a memory of a fixed finite size. Usually the memory
is incorporated into the state set, but thinking of the machine in this
way allows us to focus on the question: what information does the machine need
to remember? Answering this question tells us how to define the states.

The machine $\cal M_0$ that we construct will read edge paths in $X$ starting
at the base vertex. To motivate our description of $\cal M_0$, we start with an
informal description of the information that $\cal M_0$ has to remember as it
reads along an edge path. Consider an edge path in $X$ starting at $\{ \delta_B
\}$:
$$ \{ \delta_B \} \mapright{e_1} \{ \delta_1 \} \mapright{e_2} \cdots 
\mapright{e_I} \{ \delta_I \}
$$
This may be lifted uniquely to an edge path in $\tilde X$ starting at
$\delta_B$:
$$ \delta_B \to \delta_1 \to \cdots \to \delta_I
$$
As $\cal M_0$ reads the edge path in $X$, it will keep track of several pieces
of information. It keeps track of the combinatorial type of $\delta_i$, a
finite amount of information. It also remembers some information about how
$\delta_B$ is related to $\delta_i$. In some cases certain arcs of $\delta_k$
will be isotopic to arcs of $\delta_B$, and the automaton will remember these
arcs, again a finite amount of information. After a while, one would expect that
there are none of these arcs left. But the automaton will still keep track of a
tiny bit of information: where the end of an arc of $\delta_B$ is situated with
respect to $\delta_i$, again only a finite amount of information. This is
formalized as follow.

In order to define the states of the automaton, pick once and for all an
enumeration of the arcs of $\delta_B$, $\{g_1,\ldots,g_\kappa\}$ where $\kappa =
\kappa(g) = 12g-6$, and pick an orientation of each $g_k$, so we may speak of
the \em{tail end} and \em{head end} of each $g_k$. (Note: we shall \em{not} need
to list and orient the arcs of any other ideal triangulation; this choice is made
only for $\delta_B$.) 

Consider an arbitrary ideal triangulation $\delta$. As described in the
\tit{Tightness proposition} of \cite{M}, we may pull $\delta_B$ and $\delta$
tight with respect to one another, so that the following conditions are
satisfied:
\roster
\item"(1)" If some arc $g_k$ of $\delta_B$ is isotopic to an arc of $\delta$,
then $g_k$ is the same as some arc of $\delta$. In this case we say that
$\delta$ is \em{combed} along $g_k$. 
\item"(2)" If $\delta$ is uncombed along $g_k$ then $g_k$ is transverse to
$\delta$, and for each arc $h$ of $\delta$, there are no bigons of $g_k$ and
$h$. A \em{bigon} is a segment $\alpha \subset g_k$ and $\beta \subset h$,
neither $\alpha$ nor $\beta$ having $p$ in its interior, such that $\alpha
\union \beta$ is a simple closed curve bounding a disc.
\endroster
Furthermore, once $\delta$ is pulled tight with respect to $\delta_B$, then
$\delta$ is uniquely determined up to an isotopy preserving each arc of
$\delta_B$. 

Having pulled $\delta$ tight with respect to $\delta_B$, we may ask: Along which
arcs of $\delta_B$ is $\delta$ combed? Furthermore, if $\delta$ is not combed
along $g_k$, then the tail end of $g_k$ must emerge from some prong of $\delta$
as shown in figure 31. From which prong of $\delta$ does $\TailEnd(g_k)$ emerge?

\midinsert
\centeredepsfbox{31.ArcEmergingFromProng.eps}
\botcaption{Figure 31} $\TailEnd(g_k)$ emerges from $\pi = (e,e')$
\endcaption\endinsert

The automaton $\cal M_0$ will keep track of the answers to these questions, but
only up to the first uncombed arc in the list $\{g_1,\ldots,g_\kappa\}$. That
is, if $g_k$ is the first uncombed arc, then the automaton marks the initial ends
of $g_1,\ldots,g_{k-1}$ in $\delta$, and it marks the prong of $\delta$ from
which $\TailEnd(g_k)$ emerges. Now we formalize the concept of a ``marking'' of
$\delta$.

Let $\delta$ be a labelled ideal triangulation. Let $\Ends(\delta)$ be the set
of ends of $\delta$, and let $\Prongs(\delta)$ be the set of prongs. A
\em{marking} of $\delta$ is an injective map $\mu$, whose domain is an initial
segment $\{1,\ldots,k\}$ of $\{1,\ldots,\kappa = 12g-6\}$, and whose range is
$\Ends(\delta) \union \Prongs(\delta)$. Additional conditions will be imposed on
a marking below, but first we motivate these conditions by considering the
marking that $\delta_B$ induces on $\delta$.

Every ideal triangulation $\delta$ will have a \em{base marking}, which is
induced by $\delta_B$ as follows. First of all, if $\delta = \delta_B$, then the
initial end of each $g_k$ is marked with a $k$, and no prong is marked. Next,
assuming that $\delta \ne \delta_B$, then some arc of $\delta_B$ will be uncombed
with respect to $\delta$; let $g_k$ be the first uncombed arc. Then $\delta$ is
combed along $g_i$ for $1 \le i < k$ and the initial end of $g_i$ is marked
with an $i$ in $\delta$. Furthermore, the prong from which $g_k$ emerges will be
marked with a $k$. This defines the base marking of an ideal triangulation. 

Note that the base marking of $\delta$ is \em{not} a combinatorial invariant of
$\delta$. Nonetheless, base markings enjoy certain combinatorial properties,
which we impose as defining conditions on a general marking. For a general
marking $\mu$ of $\delta$, we require:
\roster
\item"(1)" For each arc $h$ of $\delta$, at most one end of $h$ is marked by
$\mu$.
\item"(2)" At most one prong is marked by $\mu$.
\item"(3)" If there is no marked prong, then every arc of $\delta$ has a marked
end.
\item"(4)" If there is a marked prong, it must be the last item marked.
\item"(5)" If there is a marked prong, the arc opposite the marked prong must not
have a marked end. 
\endroster
It may be that there is no marked end at all; in some sense this will be the
generic case. 

To see why property (5) is satisfied by the base marking, note in figure 31
that the arc opposite $\pi = (e,e')$ cannot have a marked end, because
otherwise that arc would then be some $g_i$, and its interior would intersect the
interior of $g_k$, contradicting the fact that $g_i, g_k \in \delta_B$ have
disjoint interiors.

A \em{marked ideal triangulation} is a pair $(\delta,\mu)$ where $\delta$ is a
labelled ideal triangulation, and $\mu$ is a marking of $\delta$. Note that
properties (1--4) are combinatorial properties, hence $\Homeo(S)$ acts on marked
ideal triangulations, and we may speak about the combinatorial type of a marked
ideal triangulation. The combinatorial type of $(\delta,\mu)$ is denoted
$\{\delta,\mu\}$.

The combinatorial type of a marked ideal triangulation $(\delta,\mu)$ may be
represented by a marked chord diagram. Let $D$ be the chord diagram of
$\delta$. If $\mu(i)$ is a marked end of $\delta$, we write the numeral $i$
adjacent to the corresponding chord end of $D$. If $\mu(k)$ is a marked prong, we
write a star $*$ next to the corresponding end gap of $D$; it is unnecessary to
actually write the $k$, because the value of $k$ can be recovered as the least
natural number which is not an end marking. Figure 32 gives an example of a
marked ideal triangulation and its chord diagram; the end marked 1 is also
the labelled end in this example.

\midinsert
\centeredepsfbox{32.MarkedChordDiagram.eps}
\botcaption{Figure 32} A marked ideal triangulation and its chord diagram
\endcaption\endinsert

The states of $\cal M_0$ may now be defined. There will be one state for each
combinatorial type of marked ideal triangulation, all of which are accept
states. Furthermore there is one failure state $F_v$ for each vertex $v$ of $X$.
The start state of $\cal M$ is the combinatorial type of the base marking
on $\delta_B$ itself.

We use the following convention for choosing the start state. We have already
chosen a labelled end for $\delta_B$. Mark that end with $k=1$, and go around the
ends in the clockwise direction; as an end of a new arc is encountered, mark
that end with the next value of $k$. The end marked $k$ is $\TailEnd(g_k)$. Figure
33 shows start states for genus 1, 2, and 3, with markings chosen by this
convention.

\midinsert
\centeredepsfbox{33.StartStates.eps}
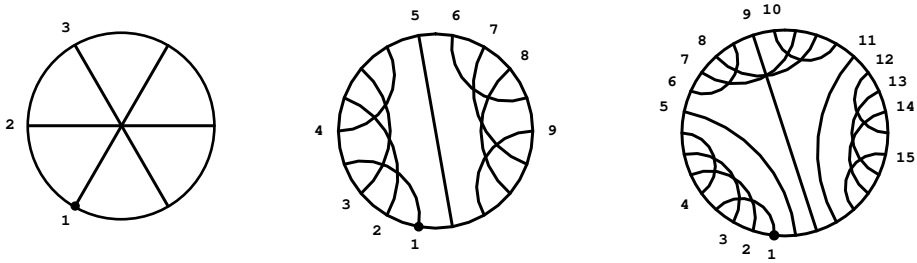
\botcaption{Figure 33} The start states for genus 1, 2, and 3
\endcaption\endinsert

\subhead Arrows of $\cal M_0$ \endsubhead

Before discussing arrows of $\cal M_0$, we need some observations and notation
concerning states. First observe how the mapping $p \from \cal M_0 \to X$ is
defined on a state $s$: identifying $s$ with a marked chord diagram, erasing
the marking leaves a labelled chord diagram which is identified with the image
vertex $D = ps$. As mentioned earlier, for any arrow $s \mapright{w} s'$ where
$w \in \cal A_0$, the image of this arrow under $p$ will be the edge in $X$
identified with $w$. The arrows can be denoted in shorthand, using the fact
that each element of $\cal A_0$ is either a relabelling generator or a labelled
elementary move generator. If $w$ is a $\Relabelling{k}$ relabelling
generator then the arrow is denoted $s \Relabelling{k} s'$, and if $w$ is a
labelled elementary move generator performed on a chord $h$ then the arrow is
denoted $s \mapright{h} s'$. In the next several paragraphs we will describe
further shorthand for determining the chord $h$.

Consider a marked ideal triangulation $(\delta,\mu)$. If there is a marked
prong, the triangle having that prong as a corner is called the \em{marked
triangle}, and the marked prong is indicated with a $*$. We usually orient the
marked triangle so that the marked prong forms a downward pointing angle,
bisected by the $-y$ direction, as in figure 34. The ends to the Left and Right
of the marked prong are denoted $e^L,e^R$ as in figure 34; formally the marked
prong is equal to the ordered pair $(e^L,e^R)$. The arcs with these ends are
denoted $h^L,h^R$. The third side of the marked triangle, the ``arc opposite the
$*$'', is denoted $h^\Opp$. In a marked chord diagram, we usually place the
``$*$'' at the bottom of the diagram, so that the chord ends $e^L,e^R$ and the
chords $h^L,h^R,h^\Opp$ are as shown in figure 34, depending on whether the
marked triangle is twisted or untwisted. 

\midinsert
\centeredepsfbox{34.LeftAndRightEnds.eps}
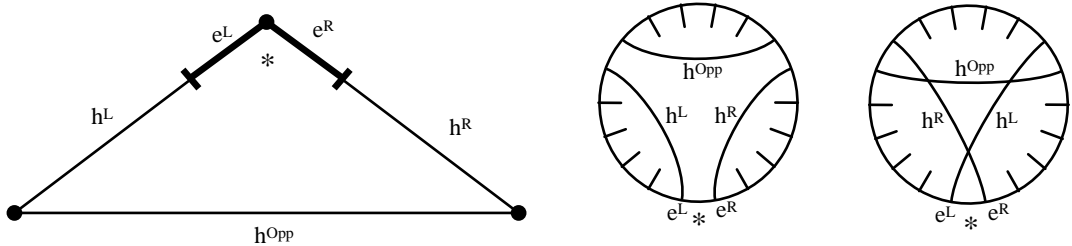
\botcaption{Figure 34} Items associated with the marked triangle
\endcaption\endinsert

Given a marked ideal triangulation $(\delta,\mu)$, we define the combinatorial
property of \em{consistency}:
\roster
\item If there is a marked end, then the end $\mu(1)$ is the labelled end.
\item If there is no marked end, then the end $e^L$ is the labelled end.
\endroster
In a marked chord diagram, consistency means that if there is a marked chord
end then the labelling dot is located at the same chord end as the numeral 1;
and if there is no marked chord end then the labelling dot is just to the left of
the $*$. Thus, in a consistent marked chord diagram the marking
determines the labelling; for this reason, we shall often leave the labelling
dot out of our chord diagrams, so any unlabelled but marked chord diagram is
assumed by default to be consistent. As we shall see, as long as an accept word
consists of labelled elementary moves, it stays among the consistent states; as
soon as the word has a relabelling generator it moves to an inconsistent state;
and if any further letters occur it moves to the failure states.

Note that the start state is consistent, by our convention for marking
$\delta_B$.

Now we describe arrows of $\cal M_0$, first the arrows leading out of failure
states. For each generator $\alpha \in \cal A_0$ leading from a vertex $v$ to a
vertex $w$ of $X$, there is an arrow $F_v \mapright{\alpha} F_w$. Thus, the set
of failure states forms a dead end set: all arrows leading out of this set lead
back into it (from which it follows that the language $\cal L_0$ of accepted
words is prefix closed). 

Next we describe relabelling arrows. Consider a relabelling generator $v
\Relabelling{r} w$ and an accept state $s$ lying over $v$. If $s$ is
inconsistent, there is a failure arrow $s \Relabelling{r} F_w$. If
$s$ is consistent, then rotate the labelling dot $r$ notches counterclockwise,
leaving the marking unsullied, to obtain an inconsistent state $s'$ lying over
$w$, and define an arrow $s \Relabelling{r} s'$. See figure 35 for an
example.

\midinsert
\centeredepsfbox{35.RelabellingArrow.eps}
\botcaption{Figure 35} A relabelling arrow
\endcaption\endinsert

To describe labelled elementary move arrows, consider a consistent marked ideal
triangulation $(\delta,\mu)$, and a labelled elementary move $\delta
\mapright{h} \delta'$. There will be an arrow $\{\delta,\mu \} \mapright{h} s$,
whose tip $s$ is some state lying over the vertex $\{ \delta' \}$, either the
failure state $F_{ \{ \delta' \} }$ or the accept state $\{ \delta',\mu' \}$ for
some marking $\mu'$ of $\delta'$. The rule for specifying $s$ is given in three
cases: 
\roster 
\item $h$ has a marked end.
\item $h$ has no marked end but one of $e^L$ or $e^R$ is an end of $h$.
\item None of the above.
\endroster
Note that these cases are combinatorial properties of the triple
$(\delta,\mu,h)$. 

Case (3) can be dispensed with immediately: the tip of the arrow is the failure
state $s=F_{ \{ \delta' \} }$. 

In each of cases (1) and (2), we must specify a marking $\mu'$ of $\delta'$, and
the tip of the arrow will be the state $\{ \delta',\mu' \}$. Despite the fact
that case (3) has been dispensed with, in certain subcases of case (3) we will
also specify a marking $\mu'$; this will be useful in describing the algorithm
for computing normal forms.

The rule for specifying $\mu'$ must be combinatorially invariant, that is, the
chord diagram for $\{ \delta,\mu \}$, together with the chord for $h$, must
determine the chord diagram for $\{ \delta',\mu' \}$. Also, assuming that $\mu$
is the base marking of $\delta$, then $\mu'$ must be the base marking of
$\delta'$. We shall keep these considerations in mind in defining the rule for
$\mu'$.

\subsubhead Case 1: $h$ has a marked end \endsubsubhead 
Suppose that an end of $h$ is marked with $j$. Assuming for the moment that
$\mu$ is the base marking of $\delta$, we may derive the base marking $\mu'$ of
$\delta'$ as follows. Since $h$ has an end marked $j$, then $h=g_j$. Also,
$\delta$ is combed along the arcs $g_1,\ldots,g_{j-1}$; since $h$ is distinct
from these arcs, then $\delta'$ is also combed along these arcs, so $\mu'$
places the marks $1,\ldots,j-1$ on the same arc ends that $\mu$ placed them.
Finally, $g_j$ is the first arc along which $\delta'$ is uncombed, and figure 36
shows the prong of $\delta'$ from which $g_j$ emerges, which is therefore the
marked prong of $\delta'$.	

\midinsert
\centeredepsfbox{36.MoveOnMarkedArc.eps}
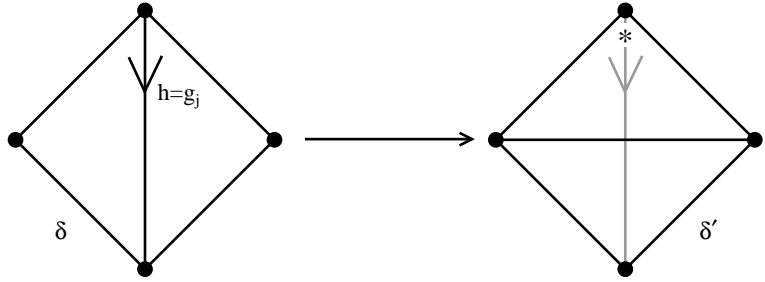
\botcaption{Figure 36} An elementary move $\delta \mapright{h} \delta'$ where
$h = g_j$ has an end marked $j$
\endcaption\endinsert

These properties of $\mu'$ may be stated as a combinatorial property in the
following manner.
\roster
\item The domain of $\mu'$ is $1 \le i \le j$.
\item If $1 \le i < j$, then $\mu(i)$ is an end of $\delta'$ as well as of
$\delta$, and we define $\mu'(i) = \mu(i)$.
\item In $\Ends(\delta)$, let $e' = \Pred(\mu(j))$ and let $e'' =
\Succ(\mu(j))$. Then the prong $(e',e'')$ is an element of $\Prongs(\delta')$,
and we define $\mu'(j) = (e',e'')$.
\endroster
We may now define an arrow $\{ \delta,\mu \} \mapright{h} \{ \delta',\mu' \}$ in
$\cal M_0$. This arrow is called a \em{$j$-marked elementary move} or a
\em{$j$-marked arrow}, and is denoted in shorthand as $\{ \delta,\mu \}
\mapright{j} \{ \delta',\mu' \}$. 

Observe that the marked ideal triangulation $(\delta',\mu')$ is still consistent.
To see why, if $j>1$ then by consistency of $(\delta,\mu)$, the labelled end of
$\delta$ is $\mu(1)$; since this is also an arc end of $\delta'$ then it is the
labelled end of $\delta'$ by definition of a labelled elementary move; but this
end is also $\mu'(1)$ by definition of $\mu'$, proving consistency. And if $j=1$
then the rule for a labelled elementary move says that the labelling is moved to
$e'$, the predecessor of $\mu(1)$, which is then the labelled end of $\delta'$;
but this equals $e^L$ in $\delta'$, and since $\delta'$ has no marked ends then
consistency is proved.

To implement a $j$-marked arrow on chord diagrams, suppose $s$ is the marked
chord diagram for $\{ \delta,\mu \}$, and the elementary move is performed on
the chord $h$ with an end marked $j$. Now erase all end markings greater than
$j$, erase the chord $h$, place a $*$ in the gap vacated by the end marked $j$,
then draw in the new chord to obtain the marked chord diagram $s'$ for
$\{\delta',\mu'\}$. Some examples are given in figure 37.

\midinsert
\centeredepsfbox{37.NonparityGoodEM.eps}
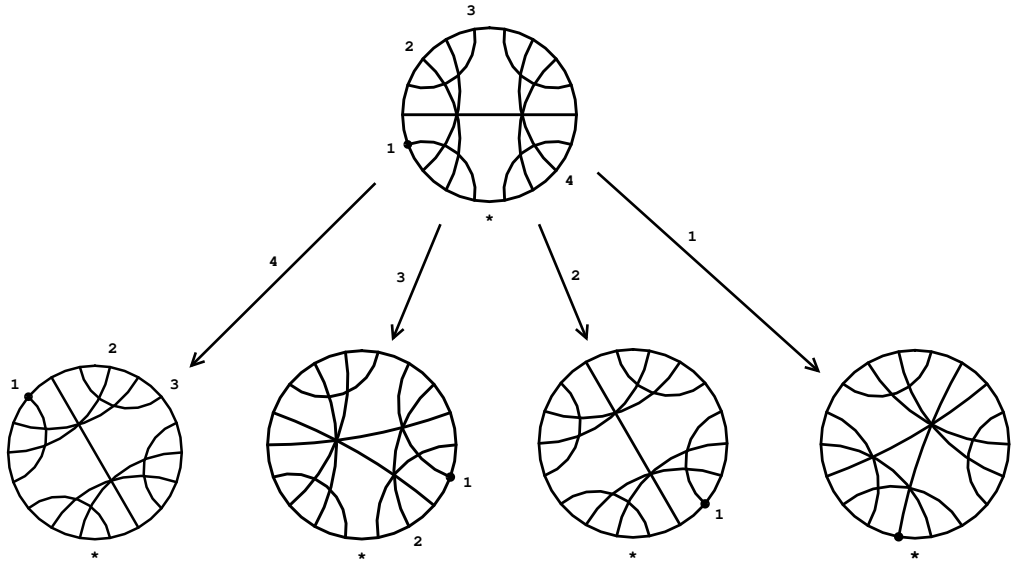
\botcaption{Figure 37} Leaving from the same state as the arrow in figure 35,
there are $j$-marked arrows for each $j=1,2,3,4$. After doing the elementary
move, putting the $*$ in the gap vacated by $j$, and erasing end marks greater
than $j$, then the chord diagram is rotated so that the $*$ appears at the
bottom.
\endcaption\endinsert

\subsubhead Case 2: $h$ has no marked end, but one of $e^L$ or $e^R$ is an end of
$h$ \endsubsubhead
We define a \em{parity} to be an element of the set $\{L,R\}$, and we often use
the variable $d$ to represent a parity. Let $* = \mu(j)$ be the
marked prong, so $* = (e^L,e^R)$. Fix $d \in \{L,R\}$ so that $e^d$ is an end
of $h$.

Assuming for the moment that $\mu$ is the base marking of $\delta$, then the
base marking $\mu'$ of $\delta'$ may be derived as follows. For $1 \le i < j$
then $\delta$ is combed along $g_i$ and $h \ne g_i$, therefore $\delta'$ is
combed along $g_i$, so $\mu'$ places the marks $1,\ldots,j-1$ on the same arc
ends that $\mu$ placed them. Figure 38 shows the two situations where $e^L$ and
$e^R$ are ends of $h$: the prong $*$ is the one from which $g_j$ emerges in
$\delta$, and when the elementary move is performed then $*$ coalesces with
another prong of $\delta$ to form a prong of $\delta'$, from which $g_j$ emerges
in $\delta'$, hence this prong is $\mu'(j)$, marked $*$ in $\delta'$. This
determines the marking $\mu'$, as shown in figure 38.

\midinsert
\centeredepsfbox{38.MoveOnLeftOrRightArc.eps}
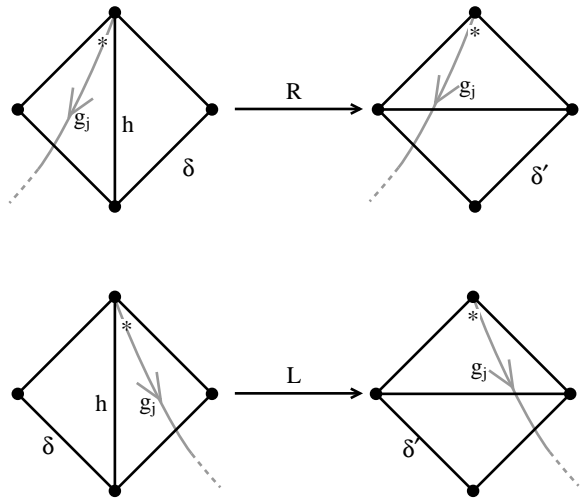
\botcaption{Figure 38} An elementary move $\delta \mapright{h} \delta'$ where
$h$ has no marked end, but $h$ has either $e^L$ or $e^R$ as an end
\endcaption\endinsert

The rule for $\mu'$ may be stated in a combinatorially invariant way as follows. 
\roster
\item The domain of $\mu'$ is equal to the domain of $\mu$, namely
$\{1,\ldots,j\}$.
\item If $1 \le i < j$ then $\mu(i)$ is an end of $\delta'$ as well as of
$\delta$, and we define $\mu'(i) = \mu(i)$.
\item If $e^d$ is an end of $h$ for some parity $d \in \{L,R\}$, in
$\Ends(\delta)$ let $e' = \Pred(e^d)$ and let $e'' = \Succ(e^d)$. Then the prong
$(e',e'')$ is an element of $\Prongs(\delta')$ and we set $\mu'(j) = (e',e'')$.
\endroster
We now define an arrow $\{ \delta,\mu \} \mapright{h} \{ \delta',\mu' \}$.
This arrow is called a \em{$d$-marked elementary move} or just a \em{$d$-arrow},
and is denoted $\{ \delta,\mu \} \mapright{d} \{ \delta',\mu' \}$. We also say
that this is a \em{parity} arrow, and a $j$-marked arrow is a \em{non-parity}
arrow. 

Observe again that the marked ideal triangulation $(\delta',\mu')$ is consistent.
If $j>1$ then by consistency of $(\delta,\mu)$ the labelled end of $\delta$ is
$\mu(1)$; and since this is not an end of $h$ then this is also the labelled end of
$\delta'$, and it is equal to $\mu'(1)$, proving consistency. Whereas if $j=1$,
then $e^L$ is the labelled end in $\delta$ by consistency; then if $d=R$, then
$e^L$ is still an end of $\delta'$ so it is the labelled end of $\delta'$, but
it is also $e^L$ in $\delta'$ proving consistency; whereas if $d=L$ then the
label is first moved to the predecessor $e'$ of $e^L$, which is then the
labelled end of $\delta'$, but this is also $e^L$ in $\delta'$, proving
consistency.	

To implement a parity $d$ elementary move on the chord diagram $s$ for
$\{\delta,\mu\}$, first locate the chord end $e^d$ adjacent to the $*$, then
erase that chord, coalescing $*$ and another gap of $s$ into a larger gap, leave
the marking $*$ in the larger gap, leave all end markings where they are,
then insert the opposite chord to form the chord diagram $s'$ for
$\{\delta',\mu'\}$. Some examples are shown in figure 39.

\midinsert
\centeredepsfbox{39.ParityGoodEM.eps}
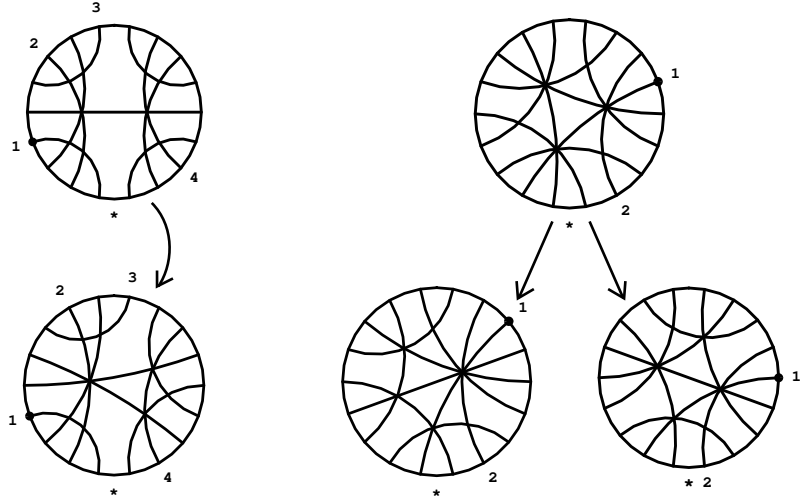
\botcaption{Figure 39} Some parity arrows. Parity is indicated by which side
of the $*$ the tail of the arrow is closest to. The first example uses, once
again, the same initial state as figure 35; there is only a Right parity arrow
with this initial state, because only $h^R$ has no marked ends. The second
example has both Left and Right parity arrows.
\endcaption\endinsert

In each of cases 1 and 2, we have constructed an accept arrow of $\cal M_0$.
Collectively, these arrows, the parity and non-parity arrows, will be
called \em{good elementary moves}. Here are a few random comments about good
elementary moves.

Comment 1: Observe that the start state is consistent, and a good elementary
move always leads from a consistent state to a consistent state. Therefore:

\proclaim{Consistency lemma} Every path of good elementary moves in $\cal M_0$,
beginning at the start state, stays among consistent states. \qed\endproclaim

\noindent This observation is what prompted us to change the normal forms from
the early version of \cite{M}, where the relabelling arrow was located at the
beginning. When the relabelling arrow is located at the end, then the
consistency lemma makes it easy to keep track of the labelled end, and there
are some simplifications in the algorithm for computing normal forms.

Comment 2: In figures, we depict good elementary moves vertically, drawing $D$
above $D'$. Bad elementary moves, defined below in case 3, are depicted as
nonvertical arrows. 

Comment 3: If $D \to D'$ is a good elementary move, then the chord diagram $D'$
always has a marked prong $*$, and the inserted chord in $D'$ is always the
chord opposite the $*$.

Comment 4: As shown by example in figure 39, given a consistent state $D$ of
$\cal M_0$, there may be both Right and Left arrows leading from $D$, or just
one, or neither, depending on which of the chords $h^L,h^R$ have marked ends.
For example, the start state never has parity arrows, since every chord has a
marked end.

\subsubhead Case 3: $h$ has no marked end, and neither $e^L$ nor $e^R$ is an
end of $h$ \endsubsubhead
We have already constructed a failure arrow $\{ \delta,\mu \} \mapright{h} F_{
\{\delta'\} }$. Nonetheless, in one subcase of case 3 we shall specify a
marking $\mu'$ of $\delta'$, and we shall say that $\{ \delta,\mu \}
\mapright{h} \{ \delta',\mu' \}$ is a \em{bad elementary move}. This will
\em{not} be an arrow in the automaton $\cal M_0$, but it will be a useful
relation among states of $\cal M_0$. In another subcase, we shall also see how
\em{inverse good elementary moves} arise. Bad and inverse good elementary moves
will be useful in describing the algorithm for computing normal forms.

Recall the notation $h^\Opp$ for the arc of $\delta$ opposite the marked prong
$*$. We consider two subcases of case 3, distinguished by whether or not $h =
h^\Opp$. Note that ``$h = h^\Opp$'' is a combinatorial property of the triple
$(\delta,\mu,h)$.

\subsubhead Case 3.1: $h \ne h^\Opp$ \endsubsubhead
Assuming $\mu$ is the base marking of $\delta$, the base marking $\mu'$ for
$\delta'$ is determined as follows. Let $* = \mu(j)$ be the marked prong of
$\delta$. Then $\delta$ is combed along the arcs $g_1,\ldots,g_{j-1}$. Not
being in case 1, then $h$ has no marked end, so $\delta'$ is also combed along
the arcs $g_1,\ldots,g_{j-1}$, and $\mu'$ places the marks $1,\ldots,j-1$ on
the same arcs that $\mu$ placed them. Not being in case 2, then $h \ne
h^L,h^R$, and being in case 3.1 then $h \ne h^\Opp$, hence the prong $*$, from
which $g_j$ emerges in $\delta$, is still a prong of $\delta'$ and $g_j$ still
emerges from it, so $\mu'$ places the mark $j$ on $*$.

This rule may be stated in a combinatorially invariant manner as follows:
\roster
\item The domain of $\mu'$ is the same as $\mu$, namely $\{1,\ldots,j\}$.
\item $\mu(i) = \mu'(i)$ for $1 \le i \le j$.
\endroster
It is evident that $(\delta',\mu')$ is consistent. We shall say that $\{
\delta,\mu\} \to \{ \delta',\mu' \}$ is a \em{bad elementary move}. We
emphasize: this does \em{not} define an arrow in $\cal M_0$, merely a relation
among states in $\cal M_0$.

The chord diagram for a bad elementary move is easily implemented: starting with
the chord diagram $D$ for $\{\delta,\mu\}$, the removed chord $h$ has no marked
end, no end adjacent to the $*$, and is not opposite the $*$; hence the end
markings and the $*$ may all be left in place as the chord is removed and the
opposite chord is inserted, yielding the chord diagram $D'$ for
$\{\delta',\mu'\}$. In figures, bad elementary moves are depicted as
nonvertical arrows, with $D'$ usually to the right of $D$, and with the removed
chord of $D$ darkened. An example is shown in figure 40. In depicting bad
elementary moves, there are no special conventions for specifying the chord on
which the move is performed, so we adopt the convention of darkening that chord.

\midinsert
\centeredepsfbox{40.BadEM.eps}
\botcaption{Figure 40} A bad elementary move
\endcaption\endinsert

\subsubhead Case 3.2: $h = h^\Opp$ \endsubsubhead
In this case, assuming that $\mu$ is the base marking of $\delta$, it is
impossible to give a combinatorially invariant description of the base marking
$\mu'$ of $\delta'$. The reason is that the elementary move $\delta' \to
\delta$, in which $h^\Opp$ is the inserted arc, gives rise to a good elementary
move $\{ \delta',\mu' \} \to \{ \delta,\mu \}$,	 and this could be either a
parity marked elementary move of either parity $d \in \{L,R\}$, or a nonparity
elementary move marked by some value $j = 1,\ldots,\kappa$; the combinatorial
type $\{ \delta',\mu' \}$ depends on the value of $d$ or $j$, and on the
placing of additional end markings in $\{ \delta',\mu' \}$. Whatever case
applies, we refer to $\{\delta,\mu\}\to\{\delta',\mu' \}$ as an \em{inverse good
elementary move}.

\subhead Normal forms \endsubhead
We have finished the construction of $\cal M_0$, and the language $\cal L_0$
accepted by $\cal M_0$ is our language of normal forms. In \cite{M} it is
proved that $\cal L_0$ contains a unique representative for each element of
$\mcgd$ whose initial vertex is the base vertex of $X$. Here is a quick summary
of the proof. 

Suppose that $\delta$ is an arbitrary unlabelled ideal triangulation on $S$.
Our task is to define a path of elementary moves $\delta_B = \delta_N \to \cdots
\to \delta_0 = \delta$, so that when the ideal triangulation $\delta_i$ is
equipped with its base marking $\mu_i$, then we obtain an accept edge in the
automaton $\{ \delta_N,\mu_N \} \to \cdots \to \{ \delta_0,\mu_0 \}$. The path
is defined in reverse order, as follows. Assume by induction that
$\delta_0,\ldots,\delta_n$ have been defined. Recall the \tit{Tightness
proposition} of \cite{M}, which says that $\delta_n$ and $\delta_B$ may be
pulled tight with respect to one another. Let $g_i$ be the first arc of
$\delta_B$ along which $\delta_n$ is uncombed. Then $g_i$ must emerge from a
certain prong of $\delta_n$, as shown in figure 31 (this is the marked prong
``$*$'', using the base marking). Now let $h$ be the first arc of $\delta$
crossed by $g_i$ (this is the arc $h^\Opp$, using the base marking). Then
$\delta_{n+1}$ is obtained by performing the elementary move $\delta_n
\mapright{h} \delta_{n+1}$. The main work of the proof is to show (1) this
sequence eventually stops at $\delta_B$; and (2) $\{ \delta_{n+1},\mu_{n+1}
\} \to \{ \delta_n,\mu_n \}$ is an arrow in $\cal M_0$.

\subhead The structure of the automaton $\cal M_0$ \endsubhead
Before proceeding with the description of the algorithm for computing normal
forms, we discuss the automaton $\cal M_0$ and its language $\cal L_0$. As
mentioned earlier, the language $\cal L_0$ forms an asynchronous automatic
structure for $\mcgd$, and we must understand certain of its properties in
order to construct a synchronous automatic structure.

We have seen that the set of failure states is a dead end set. The accept states
are partitioned into consistent and inconsistent states; each relabelling arrow
from a consistent state goes to an inconsistent state, and each arrow from an
inconsistent state goes to a failure state. All arrows leading between
consistent states are labelled elementary move arrows. This forces each word
in $\cal L_0$ to consist of a sequence of zero or more labelled
elementary moves, followed by zero or one relabelling move.

The consistent accept states can be partitioned into subsets
called \em{levels}, forming a sequence $\cal M^B_0, \cal M^\kappa_0, \ldots,
\cal M^1_0$ where $\kappa = 12g-6$. First there is the \em{base level} $\cal
M^B_0$, consisting of those consistent marked chord diagrams where every chord
has a marked end and no prong is marked. The start state lies in $\cal M^B_0$.
Going to deeper levels, each consistent accept state not in $\cal M^B_0$ has a
prong marked by some $k=1,\ldots,\kappa$, and this state lies in $\cal M^k_0$.
Recall that our convention in diagrams is to mark the prong with a $*$, and $k$
is characterized as the least integer $\ge 1$ which is not an end marking. Thus,
the deepest level $\cal M^1_0$ consists of states with no end markings. 

If $1 \le k < \kappa$, and if $s$ is a state in some level above $\cal M^k_0$,
then there is a $k$-arrow leading from $s$ into the level $\cal
M^k_0$; just do an elementary move on the chord of $s$ with end marked $k$.
Notice that each elementary move arrow leading out of the start state (or any
state in $\cal M^B_0$) leads to an accept state, because in the chord diagram
for the start state each chord has a marked end. Given $1 \le k \le \kappa$ and
a state $D$ in $\cal M^k_0$, and given $d \in \{L,R\}$, if the chord $h^d$ has no
marked end then there is a $d$-elementary move arrow leading from $D$ to another
state in $\cal M^k_0$. Thus, each arrow either stays in the same level or leads
to a lower level, the lowest level being $\cal M^1_0$. 

Note that there are many inaccessible states in $\cal M_0$. In particular, every
state in $\cal M^B_0$ except the start state is inaccessible. However, note also
that our choice of a base vertex in $X$ and a start state in $\cal M^B_0$ is
somewhat arbitrary: we could choose any state in $\cal M^B_0$ as the start state,
thereby choosing the image vertex in $X$ as the base vertex. These will all
lead to different asynchronous automatic structures on $\mcgd$.

We may count states in each level of $\cal M_0$ as follows. Let $m_g$ be the
number of vertices in $X$, i.e.\ the number of labelled chord diagrams. We have
seen earlier that $m_2 = 105$ and $m_3 = 50050$. The number of states in $\cal
M^1_0$ is $m_g$, because there is a unique way to insert a marked prong in a
labelled chord diagram to make a consistent marked chord diagram: insert the
prong just counterclockwise of the labelled end. The number of states in $\cal
M^2_0$ is $m_g (\kappa - 2)$, because the end marking 1 must be on the labelled
end, and the marked prong may be chosen freely among all $\kappa$ prongs except
that it may not be one of the two prongs opposite the marked chord. The number
of states in $\cal M^3_0$ is $m_g (\kappa-2)(\kappa-4)$, because the end
marking 1 is determined, the end marking 2 may be chosen freely among the
ends of the remaining unmarked chords, and the prong marking may be chosen
freely among the prongs not opposite the two marked chords. In general, the
number of states in $\cal M^k_0$ is $m_g (\kappa - 2)(\kappa - 4)\ldots(\kappa
- 2k + 2)$. By far most of these states are inaccessible, especially in the
base level and the highest levels. On the other hand, it is possible to show
that $\cal M^1_0$ is a strongly connected diagraph, hence all of its states are
accepssible. If one wanted to efficiently construct all the accessible states, it
would be best to use a breadth or depth first search algorithm beginning with
the start state.

From the structure of $\cal M_0$ just described, every word in $\cal L_0$ may be
factored as follows. If $w$ is an accepted word, then we may factor
$w$ into subwords as $w = w^\kappa \composed \cdots \composed w^1 \composed r$,
where the subword $w^j$, if it is not empty, begins with a $j$-arrow and is
followed by parity arrows in $\cal M^j$, and the subword $r$ is either empty or
is a single relabelling move. We call this the factorization of $w$ into
\em{uncombing blocks}; the idea is that as a new uncombing block is entered, a
new ideal arc is being uncombed. Any non-empty uncombing block $w^j$ may be
written uniquely as a single $j$-elementary move, followed by maximal subwords
of constant parity; the non-parity move is absorbed into the following subword
of contant parity, and we obtain the factorization of $w^j$ into \em{parity
blocks}. This factorization is crucial to understanding synchronization.

\head III. An algorithm for computing normal forms \endhead

To start the algorithm:

\proclaim{Input} A path of labelled elementary moves and relabelling moves,
starting at the base vertex $D_0$ of $X$:
$$w := D_0 \mapright{w_1} D_1 \mapright{w_2} \cdots \mapright{w_N} D_N$$
\endproclaim

This path can be described with pencil and paper as a sequence of labelled
chord diagrams, with the initial diagram $D_0$ chosen, say, by the convention
given in figure 30. An example is given in figure 41.

Throughout the description of the algorithm, we use the following notational
conventions. Capital letters like $D$ or $V,W,U$ will be used to denote
labelled but unmarked chord diagrams, and the lower case letter $s$ will be used
to denote marked chord diagrams, and on a consistent state we will often omit
the labelling. Given a marked chord diagram denoted with subscripts or primes,
such as $s'_0$, the corresponding unmarked chord diagram $ps'_0$ will be denoted
$V'_0$. We will also use lower case letters like $w,v$ for paths of marked or
unmarked chord diagrams.

\midinsert
\centeredepsfbox{41.InputWord.eps}
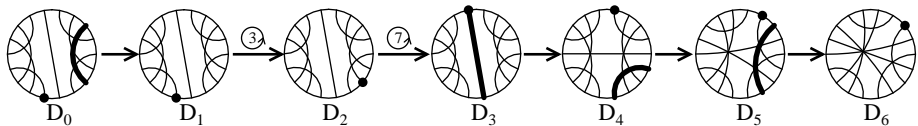
\botcaption{Figure 41} Example input for the algorithm
\endcaption\endinsert

The algorithm will work by successively computing, for $t=0,\ldots,N$, the
normal form $v^t$ representing the groupoid element $\overline w(t)$ where
$w(t)$ is the length $t$ prefix of $w$:
$$w(t) := D_0 \mapright{w_1} D_1 \mapright{w_2} \cdots \mapright{w_t} D_t
$$
Once $v^{t-1}$ is computed, then $v^{t}$ will be computed by homotoping the path
$v^{t-1} w_t$ through a sequence of relators in $X$. In section V we
shall estimate the number of relators used, and we will prove that the total
number of relators needed to calculate $v^N$ is bounded by $(12g-6)N^2$; from the
results of this section it will be clear that the task of deciding which relator
to apply at any moment takes constant time.

The algorithm is initialized by computing $v^0$ and $v^1$:

\proclaim{Initialization, step 1} Let $s_0$ be the start state of $\cal M_0$. Set
$v^0$ to be the empty path in $\cal M_0$ based at $s_0$. 
\endproclaim

This step can be implemented by marking the chord ends of the base vertex
$D_0$, say by the convention chosen in figure 33, which in genus 2 is
reproduced in figure 42.

\midinsert
\centeredepsfbox{42.Genus2StartState.eps}
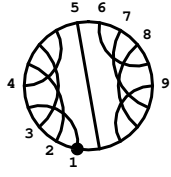
\botcaption{Figure 42} The normal form $v^0$ is the empty path based at the
start state $s_0$.
\endcaption\endinsert

\proclaim{Initialization, step 2} Given the start state $s_0$, the generator
$w_1 := (ps_0 = D_0) \to D_1$ lifts to an arrow $s_0 \mapright{w_1}
s_0'$, which is the normal form $v^1$. \endproclaim

To justify this step, note that \em{any} arrow $s_0 \mapright{w_1} s'_0$ is an
accept arrow when $s_0$ is the start state. If $w_1 := D_0 \Relabelling{r} D_1$
is a relabelling arrow this follows because $s_0$ is consistent, and we get a
relabelling generator $s_0 \Relabelling{r} s'_0$. If $w_1 := D_0
\mapright{h} D_1$ is a labelled elementary move, this follows because every
chord in $s_0$ has a labelled end, so we get a $j$-marked arrow $s_0
\mapright{j} s'_0$ where $h$ has an end marked $j$.

Using $D_0 \to D_1$ as given in figure 41, we see that this yields a $j$-marked
arrow with $j=8$, as shown in figure 43. In figure 43 and later examples, we
represent normal forms as vertical paths going downward; unmarked paths in $X$
are represented as nonvertical paths going rightward.

\midinsert
\centeredepsfbox{43.v1.eps}
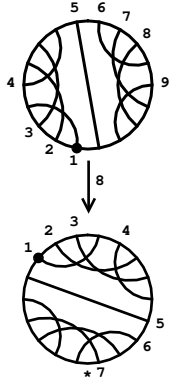
\botcaption{Figure 43} $v^1$ consists of a single arrow, the $j$-marked arrow
with $j=8$, coming out of the start state
\endcaption\endinsert

\proclaim{Main loop} For each $n=1,\ldots,N$, let $v^{n-1} := s_K \to \cdots \to
s_0$ be the normal form representing $\overline w(n-1)$. Using the subroutine
\tit{Do one move}, compute the normal form $v^n$ representing $\overline w(n)$.
\endproclaim

The reason for backwards indexing of states in a normal form, such as $v^{n-1}
:= s_K \to \cdots \to s_0$, is that our algorithm will process normal forms
from back to front.

\proclaim{Subroutine: Do one move} Let $s_K \to \cdots \to s_0$ be a normal
form, and let $V = ps_0$. Let $V\mapright{w} V'$ be a generator in $\cal A_0$,
and break into cases: for a relabelling generator, use the subroutine \tit{Do a
relabelling generator}; otherwise use \tit{Do an elementary move generator}. The
result is to compute the normal form representing the same groupoid element as
$ps_K \to
\cdots \to (ps_0 = V) \mapright{w} V'$. \endproclaim

\proclaim{Subroutine: Do a relabelling generator} Starting with a normal
form $s_K\to\cdots\to s_0$ and a relabelling generator $ps_0 = V
\Relabelling{r} V'$, break into cases depending on whether or not $s_1
\to s_0$ is a relabelling arrow. \endproclaim

\subsubhead Case 1 \endsubsubhead If $s_1 \to s_0$ is not a relabelling arrow,
then the edge $V \Relabelling{r} V'$ lifts to a relabelling arrow $s_0
\Relabelling{r} s'$, and the required normal form is $s_K \to \cdots \to
s_0 \Relabelling{r} s'$. 

An example of case 1 is shown in figure 44. We have already computed the normal
form $v^1$ for $w(1)$ in figure 43. The next generator from figure 42 is a
relabelling move $D_1 \Relabelling{3} D_2$, and the last
arrow of $v^1$ is not a relabelling arrow, so case 1 applies. The resulting
normal form $v^2$ is shown in figure 44.

\midinsert
\centeredepsfbox{44.v2.eps}
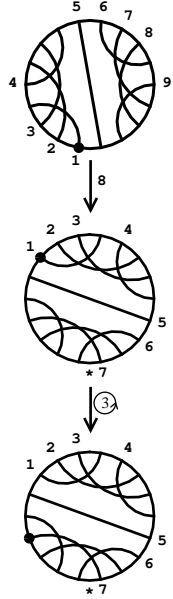
\botcaption{Figure 44} $v^2$ is obtained from $v^1$ by concatenating with a
$\Relabelling{3}$ relabelling arrow.
\endcaption\endinsert

\subsubhead Case 2 \endsubsubhead If $s_1 \Relabelling{a} s_0$ is a
relabelling arrow, then we can apply a pure relabelling relator to replace
$p(s_1) \Relabelling{a} p(s_0)=V \Relabelling{r} V'$ with $p(s_1)
\Relabelling{n} V'$, where $n \congruent r+a \modulo\kappa$.
If $n\not\congruent 0 \modulo\kappa$ then this generator lifts to a relabelling
arrow $s_1 \Relabelling{n} s'$, and the required normal form is
$s_K\to\cdots \to s_1 \Relabelling{n} s'$. If $n \congruent 0
\modulo\kappa$ then the required normal form is $s_K \to \cdots \to s_1$.

An example of case 2 is shown in figure 45. Starting from $v^2$ as in figure 44,
and using the next generator $D_2
\Relabelling{7} D_3$ from figure 42, then
the last letter of $v^2$ is a relabelling arrow, so case 2 applies and we get
$v^3$ as in figure 45.

\midinsert
\centeredepsfbox{45.v3.eps}
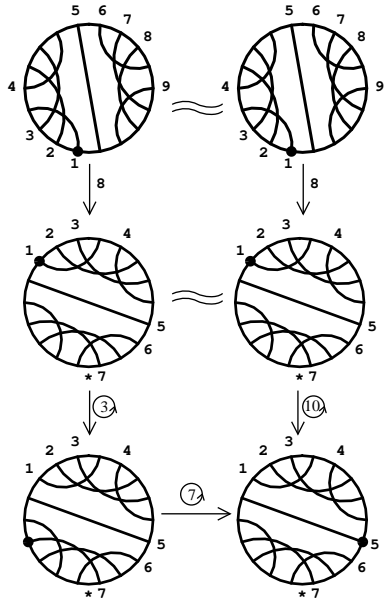
\botcaption{Figure 45} $v^3$ is obtained from $v^2$ by applying a relabelling
relator, which replaces the final $\Relabelling{3}$ arrow, followed by the
generator $D_2\Relabelling{7} D_3$, with a $\Relabelling{10}$ arrow.
\endcaption\endinsert

A word of explanation about figure 45. In paper and pencil computations, we
shall sometimes redraw a certain state, connecting the two copies of
that state via a $\approx$ sign, as in figure 45. These computations will also
explicitly show the relators that are applied, such as the relator in figure 45
which, despite the fourth side labelled $\approx$, is a three sided relabelling
relator. In general, the computations carried out by the subroutine
\tit{Do one move} are represented by such diagrams, whose input is the left hand
vertical side of the diagram followed by the bottom horizontal edge, and whose
output is the right hand side. When I do computations by hand, I usually do as
little extra copying as necessary by not copying the merged portion of the
two normal forms, but for clarity's sake the figures here will always copy the
entire merged portion. 

This finishes the subroutine \tit{Do one relabelling generator}. That was easy!

\proclaim{Subroutine: Do an elementary move generator} Starting with a normal
form $v := s_K\to\cdots\to s_0$ and a labelled elementary move generator $ps_0
= V \mapright{h} V'$, do the following steps to compute the normal form $v' :=
s'_{K'} \to \cdots \to s'_0$ representing the same groupoid element as $ps_K \to
\cdots\to ps_0 = V\mapright{h} V'$.	\endproclaim

Before giving a detailed explanation we give a brief overview, summarized
schematically in figure 46. The normal form will be computed in the backwards
direction. If the final arrow $s_1 \to s_0$ is a relabelling generator we process
that by applying an elementary move--relabelling relator, replacing $ps_1
\Relabelling{n} ps_0 = V \mapright{h} V'$ with an elementary
move $ps_1 \to V''$ followed by a relabelling move $V'' \to V'$. Then we analyze
the move $ps_1 \to V''$ into three cases: good, bad, or inverse good elementary
move. In the good and inverse good cases we quickly complete the computation of
$v'$. To handle the bad case, we compute markings on $V''$ and $V'$ to produce a
bad elementary move $s_1 \to s'_1$ followed by a relabelling arrow $s'_1 \to
s'_0$. Then we enter a loop. Typically the loop will take a sequence of good
elementary moves $s_i \to \cdots \to s_j$ followed by a bad one $s_j \to
s'_{j'}$ and, by applying either a commutator or pentagon relator, replace it
with a bad elementary move $s_i \to s'_{i'}$ followed by a sequence of good ones
$s'_{i'} \to\cdots \to s'_{j'}$; the differences $i-j$ and $i'-j'$ will always be
either 1 or 2 (untill the final relator, when 3 can also occur). This has the
effect of ``raising'' the bad elementary move $s_j \to s'_{j'}$ to a higher one
$s_i \to s'_{i'}$, closer to the end $s_K = s'_{K'}$. Eventually, one final
relator will be used to produce not a bad elementary move but instead an
equation $s_i = s'_{i'}$. A schematic diagram of the computation is given in
figure 46.

\midinsert
\vbox{
\hbox{
	\hfill
	\vbox{
  \hsize=2.0truein
		\vss
		\centeredepsfbox{46.SchematicComputation.eps}
	}
	\hfill
	\vbox{
		\vfill
		\hsize=4.0truein
		\noindent \em{Figure 46\/} 
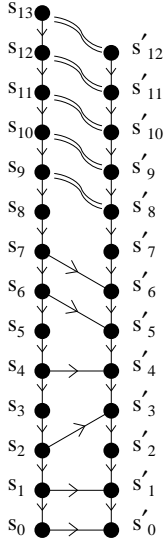
Starting from an input normal form $v := s_K\to\cdots\to s_0$ and an elementary
move generator $\alpha := ps_0 \to V'$, we typically produce an output normal
form $v' := s'_{K'} \to\cdots\to s'_0$ connected to $v$ by a sequence of
relators and bad elementary moves, so the words $v$ and $v'\alpha$ represent the
same groupoid elements. In these computations, good elementary move arrows and
relabelling arrows will be drawn vertically, whereas bad elementary moves, like
rungs of a deformed ladder, will always be nonvertical, possibly with a non-zero
vertical component. The deformation is caused by the different rate at which the
endpoints of the rungs are raised on the two sides of the ladder.
		\vfill
	}
	\hss
}
}
\endinsert

\subsubhead Step 1: Process a final relabelling arrow \endsubsubhead 
In this step, suppose that $s_1 \Relabelling{r} s_0$ is a relabelling arrow.
Apply an elementary move--relabelling relator, to replace the sequence
$ps_1 \Relabelling{r} ps_0 = V \mapright{h} V'$ with a sequence $ps_1
\mapright{h'} V'' \Relabelling{n} V'$. To explain how this is done, under the
relabelling move $s_1 \Relabelling{r} s_0$ there is a 1-1 correspondence between
chords of $s_0$ and of $s_1$, and also between chord ends. In particular, the
chord $h$ of $s_0$ corresponds to a chord $h'$ of $s_1$, yielding the elementary
move $ps_1 \mapright{h'} V''$. To see how $n$ is computed, let $e_i$ be the
labelled chord end in $s_i$. Enumerate the chord ends of $ps_1$ as
$\eta_i = \Succ^i(e_1)$, so the labelled end of $s_1$ is $\eta_0$ and the
labelled end of $s_0$ is $\eta_r$. Consider the ends
$\eta_1,\ldots,\eta_r$, and let~$a$ be the number of them which are ends of the
chord $h$. Enumerate the chord gaps as $\pi_i = (\eta_{i-1},\eta_i)$, and among
the gaps $\pi_1,\ldots,\pi_r$ let $b$ the the number into which an end of the
opposite diagonal to $h$ will be inserted. Then $n=r-a+b$ modulo
$\kappa$. It may happen that $n=0$ modulo $\kappa$, in which case the sequence 
$ps_1 \Relabelling{r} ps_0 = V \mapright{h} V'$ is replaced just with 
$ps_1 \mapright{h'} V'$.

Figure 47 gives several examples, showing how different configurations of the
ends of $h$ and its opposite diagonal can give rise to different values of $a$
and $b$. In figure 48, we continue our main example, where $v^3 := s_2
\mapright{8} s_1 \Relabelling{10} s_0$ is followed by the elementary move
$(ps_0 = D_3) \to D_4$. Applying step 1, we obtain an elementary move $ps_1 \to
V'_1$ followed by a relabelling arrow $V'_1 \to (V'_0 = D_4)$.

\midinsert
\centeredepsfbox{47.ElemMove--RelabRelators.eps}
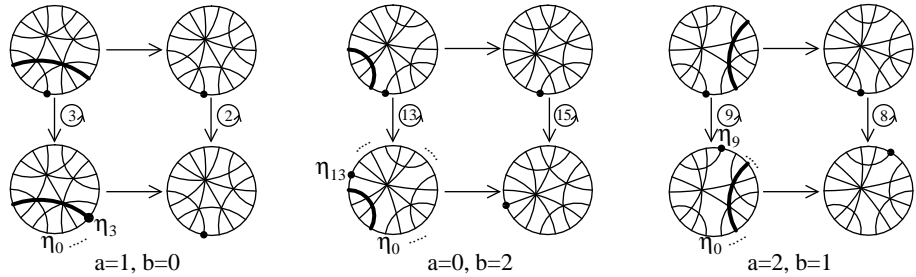
\botcaption{Figure 47} Processing a final relabelling move with an elementary
move--relabelling relator
\endcaption\endinsert

\midinsert
\centeredepsfbox{48.EndOfv4.eps}
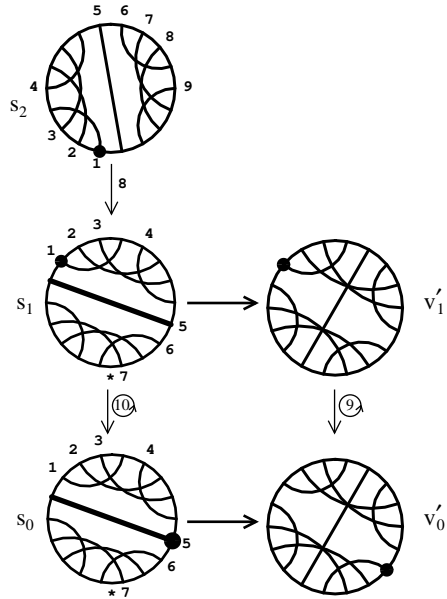
\botcaption{Figure 48} Applying step 1 to start the computation of $v^4$. The
normal form $v^3$, read off from the right side of figure 45, is reproduced
here as the left side.
\endcaption\endinsert

Having completed step 1, we now rename everything to obtain the following data:
an accepted path $s_K \to \ldots s_1$, followed by a labelled elementary move
$(ps_1 = V_1) \mapright{h} V_1'$, followed possibly by a relabelling move $V'_1
\Relabelling{r} V'_0$.

\subsubhead Step 2: Classify the elementary move $(ps_1 = V_1) \mapright{h} V'_1$
\endsubsubhead
This move will be classified as either a good elementary move, an inverse good
elementary move, or a bad elementary move. Locate $h$ in the chord diagram for
$s_1$. If $h$ has a end marked $j$, or one of the ends $e^L,e^R$, then the move
is good. If $h$ is the chord opposite the marked prong, then the move is
inverse good. Otherwise, the move is bad. 

\subsubhead Step 3: Process the elementary move $(ps_1 = V_1) \mapright{h}
V'_1$, using whichever of the three subroutines applies: \tit{Do a good
elementary move}, \tit{Do an inverse good elementary move}, or \tit{Do a bad
elementary move} \endsubsubhead

\proclaim{Subroutine: Do a good elementary move} Consider the state $s_1$ and the
chord $h$. If $h$ has an end marked $j$, then compute the $j$-marked
arrow $s_1 \to s'_1$; and if $h$ has no marked end but $e^d$ is an end of $h$ for
$d \in \{L,R\}$, compute the parity $d$ arrow $s_1 \to s'_1$.
If a relabelling move $V'_1 \Relabelling{m} V'_0$ is appended, compute the
relabelling arrow $s'_1 \Relabelling{m} s'_0$. Then $s_K \to \cdots \to s_1
\to s'_1$, with $s'_1 \Relabelling{m} s'_0$ appended if necessary, is the
normal form required to finish the subroutine \tit{Do one move}.
\endproclaim

For example, in figure 48 the elementary move $ps_1 \to V'_1$ is performed on
the chord with end marked $j=5$, so we can compute the marking on $V'_1$ by doing
a $j$-marked elementary move on $s_1$, as shown in figure 49. The marking on
$V'_0$ is then obtained by computing the $\Relabelling{9}$ arrow on
$V'_1$, also shown in figure 49. This completes the computation of $v^4 := s'_3
\mapright{8} s'_2 \mapright{5} s'_1 \Relabelling{9} s'_0$.

\midinsert
\centeredepsfbox{49.v4.eps}
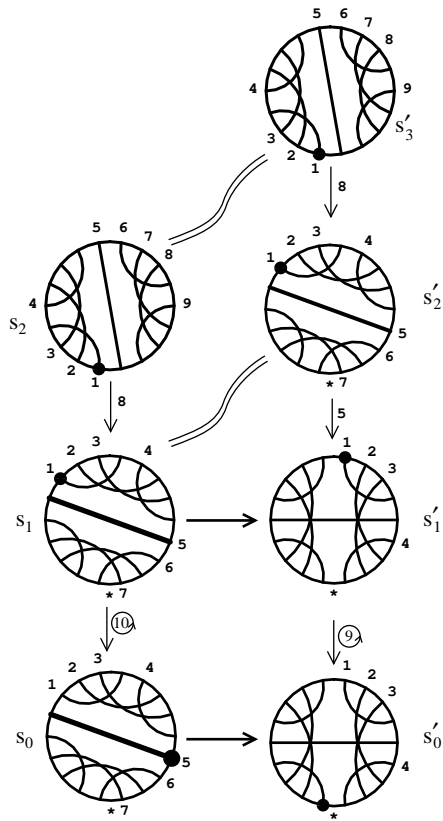
\botcaption{Figure 49} Finishing the computation of $v^4$
\endcaption\endinsert

The subroutine \tit{Do an inverse good elementary move} is next. Roughly
speaking, all we do is cancel the inverse good elementary move with the last
good elementary move. However, there is one problem: the label might not return
to its original position. More precisely, consider a labelled elementary move
$\delta \mapright{h} \delta'$ with inserted chord $h'$, and consider the
labelled elementary move on $\delta'$ performed on $h'$. If the label of $\delta$
is not on an end of $h$, then the inverse move results in $\delta$ with the label
in the same position, so there is a relator $\delta \mapright{h} \delta'
\mapright{h}\delta$; in our earlier terminology this is called an elementary
move--relabelling relator, although there are no relabelling moves in this
particular relator. On the other hand, if the label of $\delta$ is on an end of
$h$ then the inverse move $\delta' \mapright{h'} \delta''$ results in a labelled
ideal triangulation $\delta''$ obtained from $\delta$ by rotating the label one
notch clockwise, as shown in figure 50. To cancel the effect of this rotation,
we must apply a $\Relabelling{1}$ relabelling move, producing an elementary
move--relabelling relator $\delta \mapright{h} \delta'\mapright{h'}\delta''
\Relabelling{1} \delta$. This relator can be applied to replace the sequence
$\delta \mapright{h} \delta' \mapright{h'} \delta''$ by the sequence $\delta
\Relabelling{-1} \delta''$.

\midinsert
\centeredepsfbox{50.CancellingInverseMoves.eps}
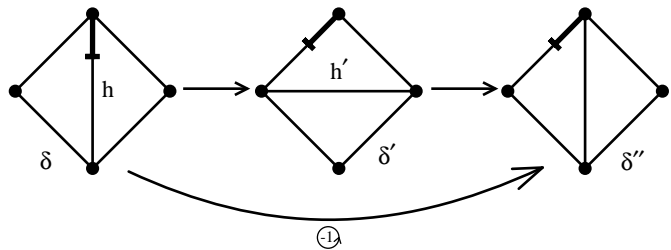
\botcaption{Figure 50} Cancelling inverse labelled elementary moves, when the
first move is performed on the labelled arc
\endcaption\endinsert

Applying these ideas to an inverse good elementary move, we obtain:

\proclaim{Subroutine: Do an inverse good elementary move, case 1} If the arrow
$s_2 \to s_1$ is a Right arrow or a $j$-marked arrow with $j \ne 1$,
so the move is not performed on the labelled chord, then apply an elementary
move--relabelling relator to replace the sequence $p(s_2) \to p(s_1) = V_1 \to
V'_1$ with the constant sequence at $p(s_2) = V'_1$. Then if there is an
appended relabelling move $V'_1 \Relabelling{r} V'_0$, append the relabelling
arrow $s_2 \Relabelling{r} s''$ to obtain the required normal form $s_K \to
\cdots \to s_2 \Relabelling{m} s''$. Otherwise, if there is no appended
relabelling move, then $s_K \to \cdots \to s_2$ is the required normal form.
\endproclaim

\proclaim{Subroutine: Do an inverse good elementary move, case 2} If the arrow
$s_2 \to s_1$ is a Left arrow or a 1-marked arrow, so the
move is performed on the labelled chord, apply an elementary move--relabelling
relator to replace the sequence $p(s_2) \to p(s_1) = V_1 \to V'_1$ with the
relabelling move $p(s_2) \Relabelling{-1} V'_1$. Then if there is an appended
relabelling move $V'_1 \Relabelling{r} V'_0$ with $r \not\congruent 1
\modulo\kappa$, apply another relabelling relator,
replacing $p(s_2) \Relabelling{-1} V'_1 \Relabelling{r} V'_0$ with $p(s_2)
\Relabelling{m} V'_0$ where $m \congruent r-1 \modulo\kappa$, and we obtain the
required normal form $s_K \to\cdots\to s_2\Relabelling{m} s''$; whereas if $r
\congruent 1 \modulo\kappa$ then the effect of the relabelling relator will be
to cancel the two relabelling moves, resulting in the normal form $s_K \to
\cdots \to s_2$. If there is no appended relabelling move, then
$s_K\to\cdots\to s_2 \Relabelling{-1} s''$ is the required normal form.
\endproclaim

An example of case 2 is shown in figure 51. The example shows the last arrow of
a normal form, a Left arrow, followed by a labelled elementary move. Applying
step 2, the elementary move is classified as an inverse good, and then \tit{Do
an inverse good elementary move} is applied.

\midinsert
\centeredepsfbox{51.DoingAnInverseGoodMove.eps}
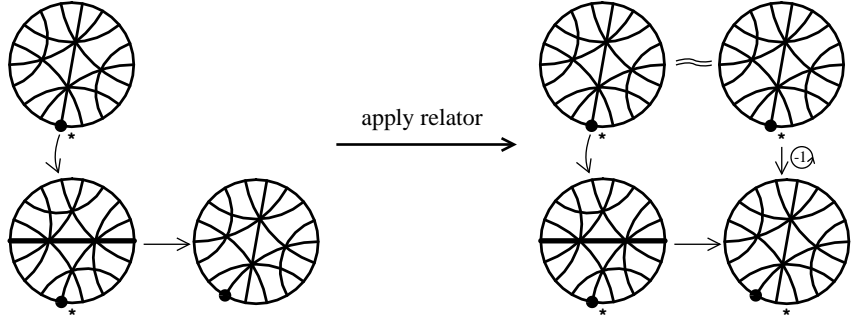
\botcaption{Figure 51} Doing an inverse good elementary move
\endcaption\endinsert

Now we come to the most laborious portion of the algorithm:

\proclaim{Subroutine: Do a bad elementary move} Initialize by computing the
state $s'_1$: in the chord diagram $s_1$, since the chord $h$ has no marked end,
and is distinct from $h^L$, $h^R$, and $h^\Opp$, then we may compute the bad
elementary move $s_1 \mapright{h} s'_1$ as described earlier. Now loop through
the subroutine \tit{Raising a bad elementary move}.
\endproclaim

\proclaim{Subroutine: Raising a bad elementary move} Given a normal form $s_K
\to \cdots \to s_1$, and a bad elementary move $s_i \to s'_j$ for some $1 \le i
< K$, apply the case analysis below, with the following effect. There is a unique
relator in $X$ having as two of its sides the labelled elementary moves $V_{i+1}
\to V_i \to V'_j$. Moreover, the vertices on this relator may be marked in a
unique way so that one of the following is true: 
\roster
\item""\tit{(Another bad elementary move)} For some $(a,b) \in
\{(1,1), (1,2), (2,1)\}$, there is a bad elementary move $s_{i+a} \to s'_{j+b}$
and a path in the automaton $s'_{j+b}\to\cdots \to s'_j$.
\smallskip
\item""\tit{(Normal forms merging)} For some $(a,b) \in \{(2,2), (1,2), (2,1),
(1,3), (3,1)\}$ there is a path in the automaton $s_{i+a} = s'_{j+b} \to \cdots
\to s'_j$.
\endroster\endproclaim

The possible outcomes are illustrated schematically in figure 52. 

\midinsert
\centeredepsfbox{52.RaisingBadElemMove.eps}
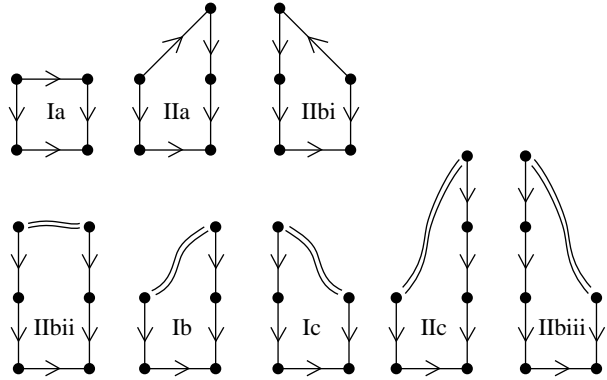
\botcaption{Figure 52} Raising a bad elementary move. The figures are labelled
according to the case analysis below.
\endcaption\endinsert

Before describing the case analysis, it should be clear how the subroutine
\tit{Do a bad elementary move} will proceed: as long as the bad elementary move
is raised to another bad elementary move using relators of types Ia, IIa, or
IIbi in figure 52, we are left with a shorter and shorter
initial segment of the original normal form $s_K \to \ldots \to s_1$. Once the
normal forms merge using a relator of type IIbii, Ib, Ic, IIc, or IIbiii, we are
done. The normal forms must eventually merge, because there is no bad elementary
move leading out of the start state $s_K$.

Now we describe the case analysis for \tit{Do a bad elementary move}. To
simplify the notation, we assume that $i=j=1$. Suppose the move $s_1 \to s'_1$ is
performed on the chord $h_0$ of $s_1$, and suppose that under the move $s_2 \to
s_1$ the inserted chord is $h_1$. If $h_0$ and $h_1$ have no adjacent ends in
$s_1$ then a commutator relator applies, otherwise a pentagon relator applies. 

\subhead Case I: $h_0$ and $h_1$ have no adjacent ends in $s_1$ \endsubhead
The commutator relator in $X$ may be described as 
$$V_1 \mapright{h_0} V'_1 \mapright{h_1} U \mapright{h_2} V_2 \mapright{h_3} V_1
$$
in one direction, and in the opposite direction as
$$V_1 \mapright{h_1} V_2 \mapright{h_0} U \mapright{h_3} V'_1 \mapright{h_2} V_1
$$
as shown in figure 53. We do not yet say in which direction these elementary
moves are labelled. 

\midinsert
\centeredepsfbox{53.CommutatorRelator.eps}
\botcaption{Figure 53} A commutator relator
\endcaption\endinsert

Now compute as follows. Locate $h_0$ and $h_3$ in the marked chord diagram
$s_2$. Decide whether: (a) $h_0$ has neither a marked end nor an end adjacent to
the marked prong of $s_2$, nor is $h_0$ opposite the marked prong of $s_2$; (b)
$h_0$ has a marked end or an end adjacent to the marked prong of $s_2$; or (c)
$h_0$ is opposite the marked prong of $s_2$. 

Note: if $s_2 \to s_1$ is a parity arrow then case (a) applies. For the end map
induces up a 1-1 correspondence between marked ends of $s_2$ and of $s_1$, and
since $h_0$ is unmarked in $s_1$ then it is unmarked in $s_2$. Also, all prongs
outside the support of a parity elementary move are unmarked in both the source
and target of the elementary move, and since $h_0$ is not adjacent to a marked
chord then the two triangles adjacent to $h_0$ are outside the support of $s_2
\to s_1$, so the prongs opposite $h_0$ are unmarked.

\subsubhead Case Ia: Neither---nor \endsubsubhead Compute the bad elementary move
$s_2 \mapright{h_0} s'_2$, so $ps'_2=U$. Now locate $h_3$ in $s'_2$, and compute
the arrow $s'_2 \mapright{h_3} s'_1$. Two examples are given in figure
54, one with parity arrows and another with non-parity arrows.

In this case as in all later cases, the algorithm computes certain arrows, but
to be formally correct we must justify that these arrows exist; the reader may
want to skip these justifications at first.

The arrow $s'_2 \mapright{h_3} s'_1$ exists because $h_3$ has the
same relationship with the marking in $s'_2$ as in $s_2$, i.e.\ it either has a
marked end or has an end adjacent to the marked prong. Moreover the arrows
$s_2 \to s_1$ and $s'_2 \to s'_1$ are of the same type: both have the same
parity, or both are $j$-marked for the same $j$. 

\midinsert
\centeredepsfbox{54.CaseIa.eps}
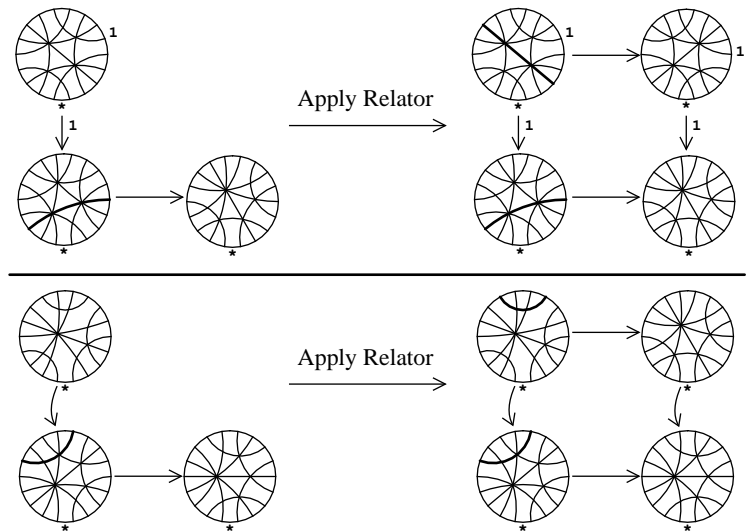
\botcaption{Figure 54} Examples of case Ia
\endcaption\endinsert

\subsubhead Case Ib: $h_0$ has either a marked end or end adjacent to the marked
prong of $s_2$ \endsubsubhead 
As noted above, this happens only if the arrow $s_2 \to s_1$ is $j$-marked for
some $j$, hence $h_3$ has an end marked $j$. Compute the arrow $s_2 = s'_3
\mapright{h_0} s'_2$ in the automaton, so $ps'_2=U$. Now locate $h_3$ in $s'_2$,
and compute the arrow $s'_2 \mapright{h_3} s'_1$. Two examples are given in
figure 55, one where $s'_3 \to s'_2$ is a parity arrow and one where it is
non-parity.

To justify why the arrows $s'_2 \mapright{h_3} s'_1$ exists, the end map sets up
a 1-1 correspondence between ends in $s_2$ and in $s_1$ which are marked by some
$i<j$. Therefore, since $h_0$ is unmarked in $s_1$, then if $h_0$ has a marked
end in $s_2$ that marking must be $\ge j$, and it must be $>j$ since $h_3$ is
marked with $j$; it follows that after the arrow $s'_3 \mapright{h_0} s'_2$ then
$h_3$ is still marked $j$ in $s'_2$. Also, if $h_0$ is adjacent to the marked
prong of $s_2$ that prong must be marked $>j$, so $h_3$ is still marked $j$ in
$s'_2$. Thus, $s'_2 \mapright{h_3} s'_1$ is a $j$-marked arrow. 

\midinsert
\centeredepsfbox{55.CaseIb.eps}
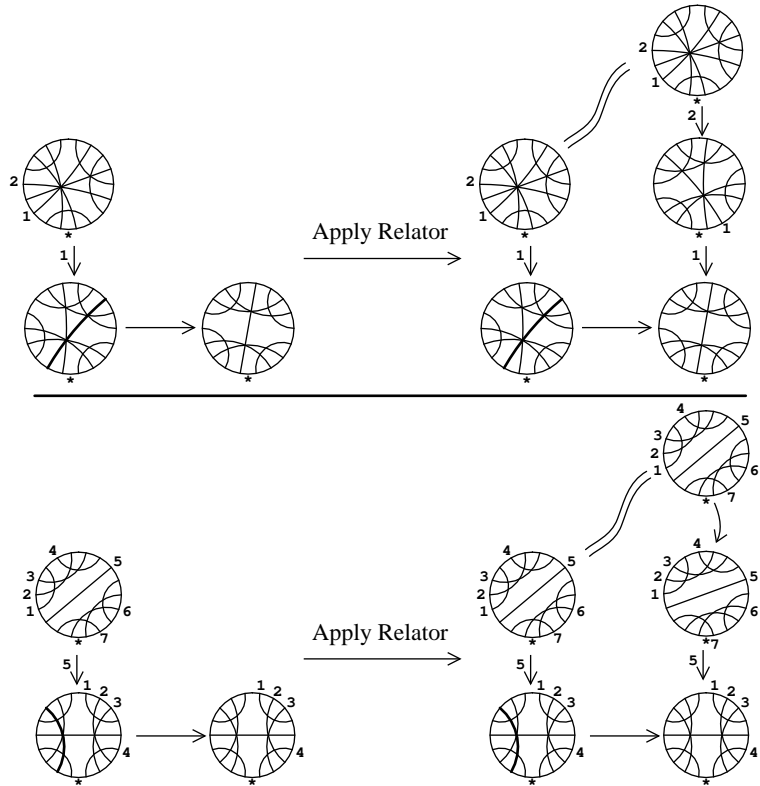
\botcaption{Figure 55} Examples of case Ib
\endcaption\endinsert

\subsubhead Case Ic: $h_0$ is opposite the marked prong of $s_2$ \endsubsubhead
Again, this happens only if the arrow $s_2 \to s_1$ is $j$-marked for
some $j$, hence $h_3$ has an end marked $j$. Since $h_0$ is opposite the marked
prong of $s_2$, then $h_0$ is the inserted chord under the arrow $s_3 \to s_2$,
and $ps_3=U$. Now locate $h_3$ in $s_3$, and compute the arrow $(s_3 = s'_2)
\mapright{h_3} s'_1$. Examples are given in figure 56. Note that occurences of
case Ib and Ic are ``orientation reversals'' of each other; c.f.\ figures 55,56.

The arrow $(s_3 = s'_2) \mapright{h_3} s'_1$ exists because $h_3$, being marked
by $j$ in $s_2$, is also marked by $j$ in $s_3$. Thus, the arrow $s'_2
\mapright{h_3} s'_1$ is $j$-marked. 

\midinsert
\centeredepsfbox{56.CaseIc.eps}
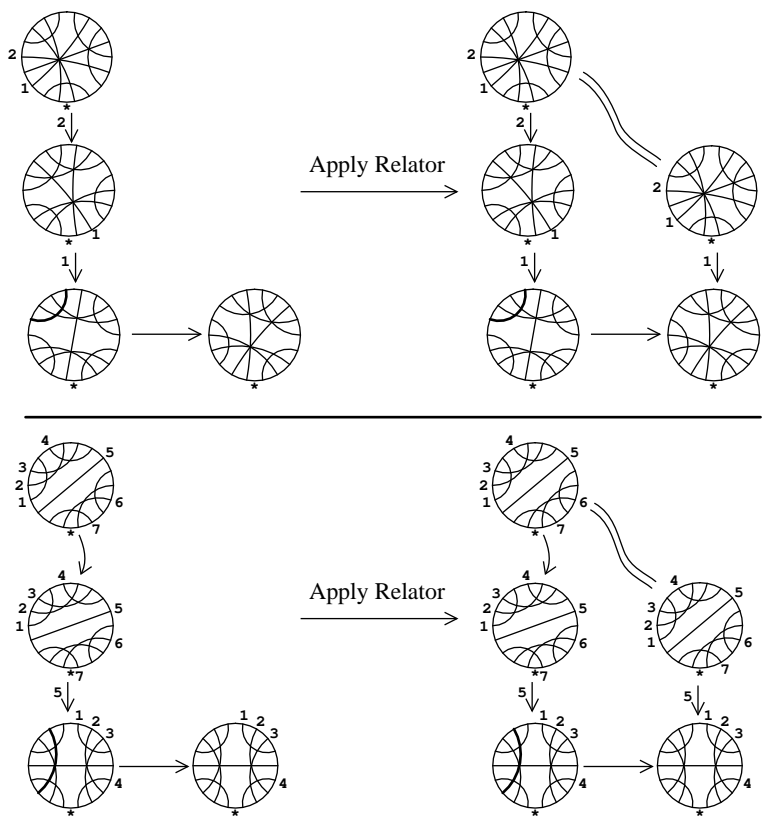
\botcaption{Figure 56} Examples of case Ic
\endcaption\endinsert

\subhead Case II: $h_0$ and $h_1$ have adjacent ends in $s_1$ \endsubhead
The pentagon relator in $X$ may be described as
$$V_1 \mapright{h_0} V'_1 \mapright{h_1} W \mapright{h_2} U \mapright{h_3} V_2
\mapright{h_4} V_1
$$
and in the other direction
$$V_1 \mapright{h_1} V_2 \mapright{h_0} U \mapright{h_4} W \mapright{h_3} V'_1
\mapright{h_2} V_1
$$
as shown in figure 57.

\midinsert
\centeredepsfbox{57.PentagonRelator.eps}
\botcaption{Figure 57} A pentagon relator
\endcaption\endinsert

Now compute. Locate $h_0$ and $h_4$ in the marked chord diagram $s_2$. Decide
whether: (a) $h_0$ has no marked end nor end adjacent to the marked prong of
$s_2$, nor is $h_0$ opposite the marked prong of $s_2$; (b) $h_0$ is opposite
the marked prong of $s_2$; or (c) $h_0$ has a marked end or an end adjacent to
the marked prong of $s_2$.

\subsubhead Case IIa: Neither---nor \endsubsubhead Compute the bad elementary
move $s_2 \mapright{h_0} s'_3$, whose inserted chord is $h_3$, so $ps'_3=U$. Now
locate the chord $h_4$ in $s'_3$ and compute the arrow $s'_3 \mapright{h_4}
s'_2$, so $ps'_2=W$. Finally, locate the chord $h_3$ in $s'_2$, and compute the
arrow $s'_2 \mapright{h_3} s'_1$. Examples are shown in figure 58. 

To see why the arrows $s'_3 \mapright{h_4} s'_2$ and $s'_2 \mapright{h_3} s'_1$
exist, first note that in figure 57, the marking on $s_2$ must include either
$a$ or $b$ in $V_2$: the marking must include one of $a$--$f$, since $s_2
\mapright{h_4} s_1$ is an arrow; but $c$ and $f$ are excluded because no mark
can be on a prong adjacent to or opposite $h_0$; also $d$ and $e$ are excluded
because then the marking on $s_1$ would include $b$ in $V_1$, violating the fact
that $s_1 \mapright{h_0} s'_1$ is a bad elementary move. It follows that the
marking on $s'_3$ includes $a$ or $b$ in $U$, showing that $s'_3 \mapright{h_4}
s'_2$ is an arrow, of the same type as $s_2 \to s_1$. Also, it follows that the
marking on $s'_2$ includes $a$ in $W$, showing that there is an arrow $s'_2
\mapright{h_3} s'_1$; this can be a parity or nonparity arrow, depending on
whether $h_3$ has a marked end. 

\midinsert
\centeredepsfbox{58.CaseIIa.eps}
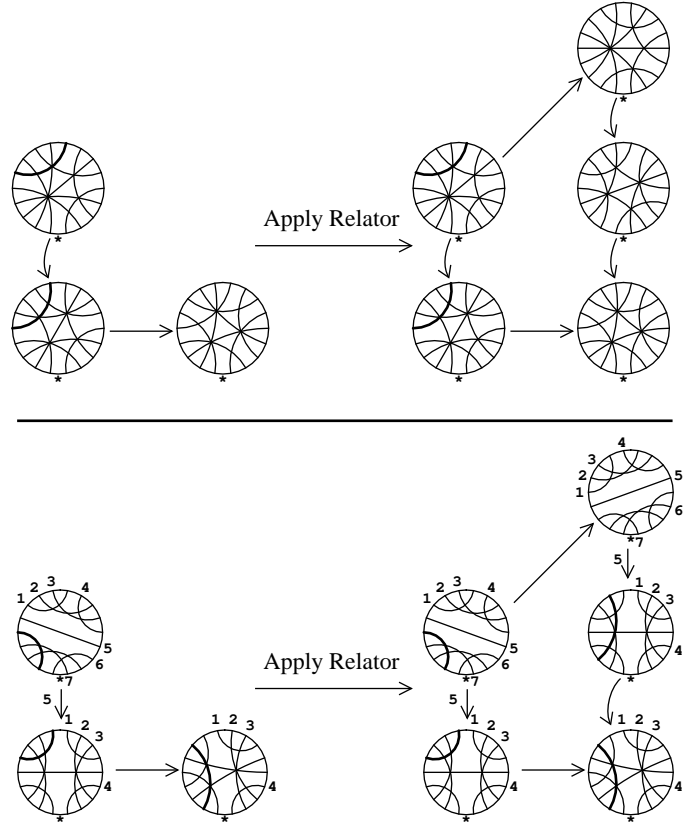
\botcaption{Figure 58} Examples of case IIa
\endcaption\endinsert

\subsubhead Case IIb: $h_0$ is opposite the marked prong of $s_2$
\endsubsubhead	Then the arrow $s_3 \to s_2$ has inserted chord $h_0$, so
$ps_3=U$, and $s_3 \to s_2$ is performed on $h_3$. Locate the chord
$h_4$ in $s_3$. Decide whether: (i) $h_4$ has no marked end, no end
adjacent to the marked prong of $s_3$, and $h_4$ is not opposite the marked prong
of $s_3$; (ii) $h_4$ has a marked end or end adjacent to the marked prong of
$s_3$; (iii) $h_4$ is opposite the marked prong of $s_3$; one of these must
happen.

\subsubhead Case IIbi \endsubsubhead Compute the bad elementary move $s_3
\mapright{h_4} s'_2$, so $ps'_2 = W$. Now locate the chord $h_3$ in $s'_2$, and
compute the arrow $s'_2 \mapright{h_3} s'_1$. Examples are given in figure 59.
Note that cases of IIa and IIbi are orientation reversals of each other; c.f.\
figures 58,59.

To see why the arrow $s'_2 \mapright{h_3} s'_1$ exists, note that the marking on
$s_3$ must include one of $d,e,j,k$ in $U$: since $s_3 \mapright{h_3} s_2$ is an
arrow then one of the marks $c$--$d,i$--$k$ must be included, but the marks
$c,i$ are forbidden because no prong adjacent to or opposite $h_4$ is marked; it
follows that the marking on $s'_2$ must include one of $b,c,h,i$ in $W$, so
$s'_2 \mapright{h_3} s'_1$ is an arrow. 

\midinsert
\centeredepsfbox{59.CaseIIbi.eps}
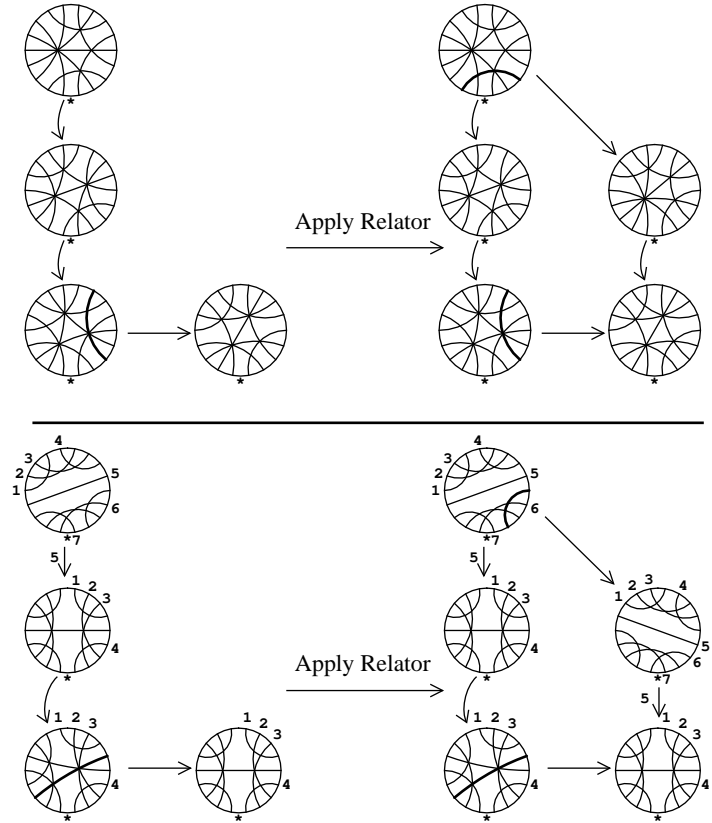
\botcaption{Figure 59} Examples of case IIbi
\endcaption\endinsert

\subsubhead Case IIbii: $h_4$ has a marked end or end adjacent to the marked
prong of $s_3$ \endsubsubhead Compute the arrow $s_3 = s'_3 \mapright{h_4}
s'_2$, so $ps'_2 = W$. Now locate the chord $h_3$ in $s'_2$, and compute the
arrow $s'_2 \mapright{h_3} s'_1$. Examples are given in figure 60. Note that the
orientation reversal of an example of case IIbii is another example of case
IIbii; it is instructive to study the orientation reversals of the examples in
figure 60.

To see why the arrow $s'_2 \mapright{h_3} s'_1$ exists, note first that the
marking on $s_3$ includes one of $c$--$e$ or $i$--$k$ in $U$, because $s_3
\mapright{h_3} s_2$ is an arrow; we can rule out $i$ because either $h_4$ has a
marked end and no marked prong can be opposite an arc with a marked end
(property 5 in the definition of a marking), or the unique marked prong is
adjacent to $h_4$; it follows that the marking on $s'_2$ includes one of
$a$--$c,h,i$ in $W$, so $s'_2 \mapright{h_3} s'_1$ is an arrow. 

\midinsert
\centeredepsfbox{60.CaseIIbii.eps}
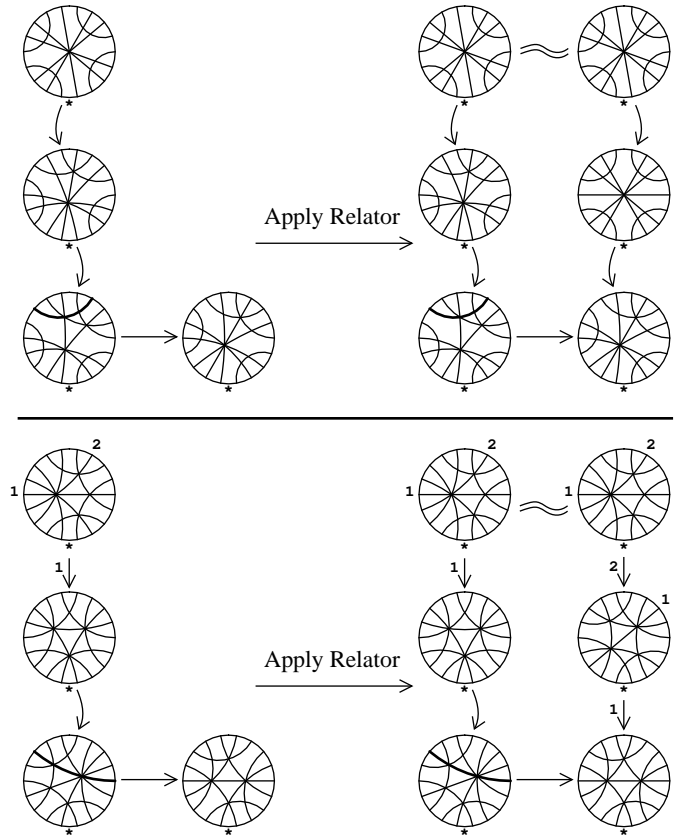
\botcaption{Figure 60} Examples of case IIbii
\endcaption\endinsert

\subsubhead Case IIbiii: $h_4$ is opposite the marked prong of $s_3$
\endsubsubhead Then $h_4$ is the chord inserted under $s_4 \to s_3$, so
$ps_4=W$. Locate the chord $h_3$ in $s_4$, and compute the arrow
$(s_4 = s'_2) \mapright{h_3} s'_1$. An example is given in figure 61.

To see why the arrow $(s_4 = s'_2) \mapright{h_3} s'_1$ exists, note that the
marking on $s_3$ must include one of $c$--$e$ or $i$--$k$ in $U$, since $s_3
\mapright{h_3} s_2$ is an arrow. However, $c$, $e$, and $k$ may be eliminated
by the requirement that the unique marked prong is opposite $h_4$. Also, $i$
may be eliminated, for if $i$ is included then $b$ and $g$ are not included, so
the marking on $s_2$ includes $i$ but not $b$ and $e$ in $V_2$, and by
uniqueness of the marked prong $a$,$c$,$d$, and $f$ are also not included, but
this violates the requirement that $s_2 \mapright{h_4} s_1$ is an arrow. Thus,
the marking on $s_3$ includes one of $d$ or $j$ in $U$. It follows that the
marking on $s_4$ includes one of $b$ or $h$ in $W$, so $(s_4 = s'_2)
\mapright{h_3} s'_1$ is an arrow. 

\midinsert
\centeredepsfbox{61.CaseIIbiii.eps}
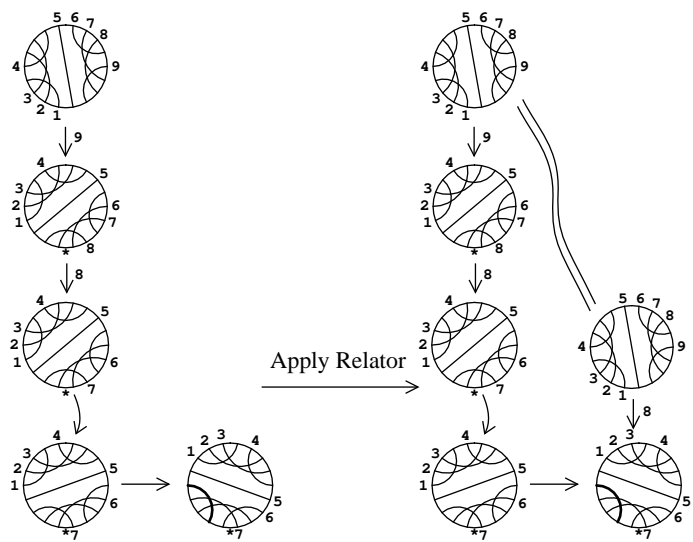
\botcaption{Figure 61} Example of case IIbiii
\endcaption\endinsert

\subsubhead Case IIc: $h_0$ has a marked end or end adjacent to the marked prong
of $s_2$ \endsubsubhead Compute the arrow $(s_2 = s'_4) \mapright{h_0} s'_3$,
whose inserted chord is $h_3$; note that $ps'_3=U$. Now locate the chord $h_4$
in $s'_3$, and compute the arrow $s'_3 \mapright{h_4} s'_2$, so $ps'_2 = W$.
Finally, locate the chord $h_3$ in $s'_2$ and compute the arrow $s'_2
\mapright{h_3} s'_1$. An example is given in figure 62. Note that occurences of
cases IIbiii and IIc are orientation reversals of each other; c.f.\ figures 61,
62.

To see why the arrows $s'_3 \mapright{h_4} s'_2$ and $s'_2 \mapright{h_3} s'_1$
exist, first note that the marking on $s_1$ must include one of $a$ or $b$ in
$V_1$, because $h_1$ is inserted under the arrow $s_2 \to s_1$; however, $b$ is
eliminated since $s_1 \mapright{h_0} s'_1$ is a bad elementary move, so $a$ is
included. It follows that the marking on $s_2$ includes one of $a$, $b$, or $c$
in $V_2$; let the numerical value of this marking be $n$. By hypothesis, the
marking on $s_2$ includes some end of $h_0$ or prong adjacent to $h_0$ in $V_2$;
let the numerical value of this marking be $m$. Then $m>n$, for after the arrow
$s_2\to s_1$ there is no marked end of $h_0$, nor marked prong adjacent to $h_0$,
because $s_1 \mapright{h_0} s'_1$ is a bad elementary move. If $a$ or $c$ is
marked with $n$ in $s_2$ we obtain a contradiction, since the prong marking is
greater than all end markings. Thus, $b$ is marked with $n$ in $s_2$. It follows
that after the arrow $(s_2 = s'_4) \mapright{h_0} s'_3$, the marking on
$s'_3$ includes $b$ in $U$. It then follows that there is an arrow $s'_3
\mapright{h_4} s'_2$, and that the marking on $s'_2$ includes $a$ in $W$.
Finally, it follows that there is an arrow $s'_2 \mapright{h_3} s'_1$. 

\midinsert
\centeredepsfbox{62.CaseIIc.eps}
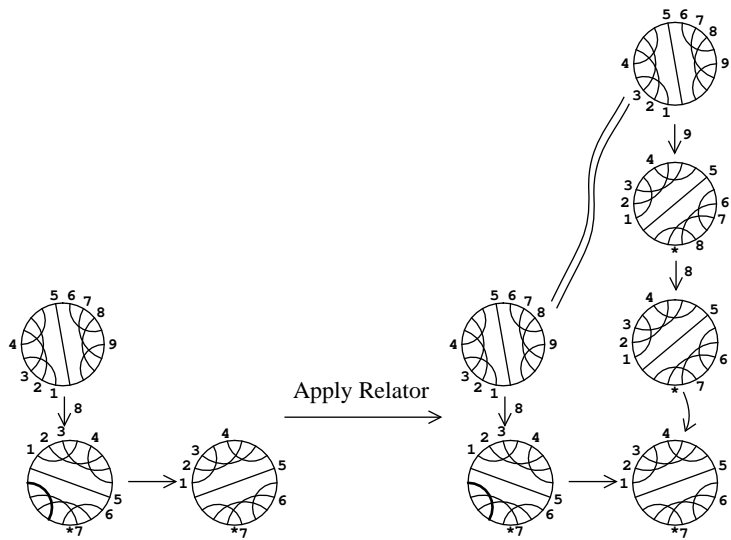
\botcaption{Figure 62} Example of case IIc
\endcaption\endinsert

This finishes the subroutine \tit{Do a bad elementary move}. A few comments:

Comment 1: Strictly speaking, we have gone around a relator in two ways to
obtain markings on $ps'_1 = V'_1$, and we should check that these two markings
are identical. It is obvious that the markings are identical outside of the
support of the relator, and by checking cases one may see that the markings are
identical in the support of the relator; alternatively, use the results of
\cite{M}.

Comment 2: When the arrows on the left side of the relator are all parity
arrows, then the relator must be case Ia, IIa, IIbi, or IIbii. Moreover, the
arrows on the right side must also be parity arrows. Moreover, in cases I, IIa,
and IIbi all of the arrows in the relator have the same parity; in case IIbii,
the arrows on the left side have the parities $dd'$ for some choice of $d \ne d'
\in \{L,R\}$, and the arrows on the right side have parities $d'd$. The
diagrams for each of these cases, figures 54,58,59,60, each show an example
where the arrows are all parity arrows. This observation can be used to make
some computational shortcuts: once it has been determined that the arrows on
the left hand side of the relator are all parity arrows, and once the top chord
diagram on the right hand side of the relator has been computed, then the
parities of the arrows on the right hand side are determined, and from this
information the arrows may be computed.

This completes the description of the algorithm for computing normal forms.

Finally, observe from the description of the algorithm that $\cal L_0$ is an
asynchronous automatic structure, because the input and output normal forms under
any run of the subroutine \tit{Do one move} are asynchronous fellow travellers. 

\subhead Examples of doing a bad elementary move \endsubhead

Figures 63,64 show some examples of applying the subroutine \tit{Do a bad
elementary move}. Figure 63 finishes the computation of the normal form	 for the
example word given in figure 41, using two more applications of
\tit{Do a bad elementary move}. The left side of the figure shows $v^4$ as
computed in figure 49, and the bottom of the figure shows the last two
generators $D_4 \to D_5 \to D_6$ from figure 41. Then $v^5$ is computed from
$v^4$ using an elementary move--relabelling relator, followed by a run of the
subroutine \tit{Do a bad elementary move}, using relators of types IIa, IIc; then
$v^6$ is computed using relators of types IIa, IIbi. Some of the relators used
in figures 63 and 64 were described in figures 54--62, and the numerals written
on these relators refer to the relevant figure. 

\topinsert
\vbox spread 12pt{
\hbox{
 \hss
	\vbox{
		\hsize=3.0truein
		\vfill
		\centeredepsfbox{63.AlgorithmExample.eps}
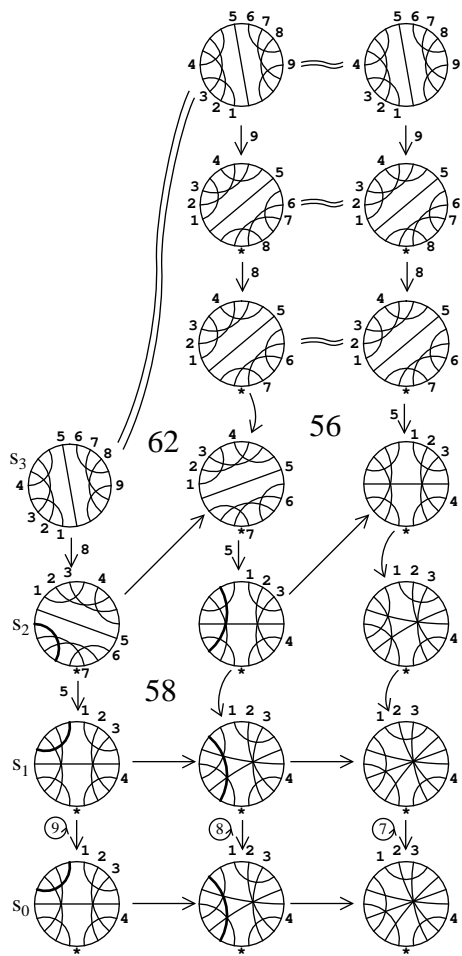
		\centerline{\it Figure 63}
	}
	\hfill\vrule\hfill
	\vbox{
		\hsize=3.0truein
		\centeredepsfbox{64.AlgorithmExample2.eps}
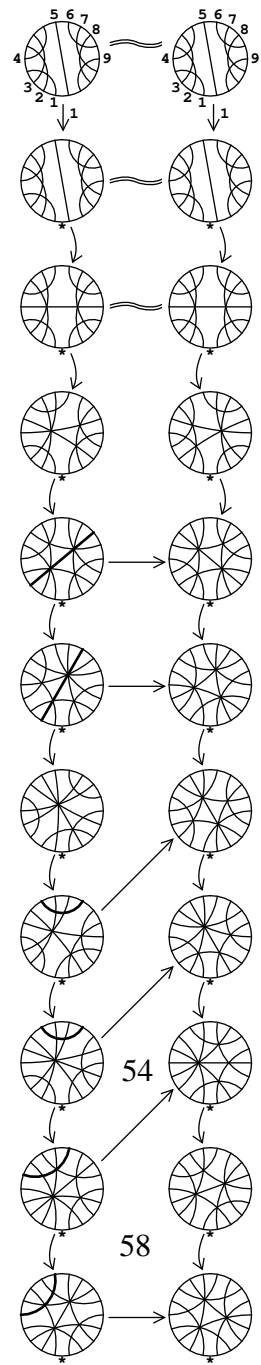
		\centerline{\it Figure 64}
	}
 \hss
}
}
\endinsert

Figure 64 shows another application of \tit{Do a bad elementary move}. This is a
more typical example than figure 63: given the fact that an arbitrary element of
$\cal L_0$ has at most $\kappa = 12g-6$ nonparity arrows, if the word is very
long the relators do not interact very often with nonparity arrows.

\head IV. Interlude: The suffix uniqueness property for an automatic structure
\endhead

In this section, we describe some properties of $\cal L_0$, as a motivation for
our proof in the next section that the computation of normal forms in $\cal
L_0$ runs in quadratic time.

In \cite{ECHLPT} it is proved that an automatic group has a quadratic time
algorithm for the word problem. We define a property of automatic structures
called \em{suffix uniqueness}, and using this property we give another
description of a quadratic time algorithm for the word problem. Our algorithm is
more efficient than the one described in \cite{ECHLPT}, as can be seen by
comparing the proofs. For asynchronous automatic structures, the method of
\cite{ECHLPT} yields an exponential time algorithm; the property of suffix
uniqueness applies equally well to asynchronous automatic structures, and in
this case we obtain an exponential time algorithm as well.

Let $G$ be a groupoid with finite generating set $\cal A$, and let $\cal L$ be a
synchronous or asynchronous automatic structure over $G$, with word acceptor
$\cal M$ and fellow traveller constant $K$. Let $B_K$ be the set of elements in
$G$ represented by words of length $\le K$ in the generating set $\cal A$;
elements of $B_K$ do not have to start at the base point of $\cal L$. Let ${}_n
\cal L$ be the set of suffixes of $\cal L$ of length between $1$ and $n$.
Assuming that $\cal M$ has no inaccessible states, then ${}_n \cal L$ is the set
of all nontrivial paths in $\cal M$ of length at most $n$ ending at an accept
state. We say that $\cal L$ satisfies \em{suffix uniqueness} if it is prefix
closed, no two normal forms represent the same element of $G$, and there exists
an integer $n \ge 1$, a finite subset $\cal S \subset {}_n\cal L \cross B_K$, and
a function $F \from \cal S \to {}_k \cal L$, with the following property:
\roster
\item"" If $w,w' \in \cal L$ are $K$-fellow travellers, then setting $g =
\overline w^\inverse w' \in B_K$, there exists a unique suffix $s$ of $w$ such
that $(s,g) \in \cal S$, and $F(s,g)$ is a suffix of $w'$. Moreover, writing $w
= \hat w s$ and $w' = \hat w' F(s,g)$, then $\hat w$ and $\hat w'$ are
$K$-fellow travellers.
\endroster

For example, from the description of the subroutine \tit{Do one move} it
follows that the asynchronous automatic structure $\cal L_0$ for the
groupoid $\mcgd$ satisfies suffix uniqueness, with fellow traveller constant
$K=1$ and maximal suffix length $n=3$.

From the suffix uniqueness property, we obtain an algorithm for the word problem
as follows. Given an arbitrary word $w = w_1 \cdots w_M$, by induction compute
the normal form $v^m$ representing $w(m) = w_1 \cdots w_m$. To do this, suppose
$v = v^{m-1}$ is computed. We must compute the normal form $v' = v^m$
representing $w^{m-1} w_m$. We compute $v'$ by induction, producing longer and
longer suffixes of $v'$. Set $a_0 = v$ and $a'_0 = v'$. Since $a_0$ and $a'_0$
are $K$-fellow travellers, then we may factor $a_0$ uniquely as $v=a_1 s_0$ so
that $(s_0,w_m) \in \cal S$, and then set $s'_0 = F(s,w_m)$, so $s'$ is a suffix
of $v'$ and $v' = a'_1 s'_0$ for some $a'_1$. It follows that $a_1$ and $a'_1$
are $K$-fellow travellers, so we may continue by induction. Since $a_i$ is
decreasing in length, eventually we compute $v' = s'_J \ldots s'_1 s'_0$.

Note that this algorithm is exactly the same as the algorithm described in the
previous section, for computing normal forms in $\cal L_0$ representing
elements of $\mcgd$. 

The computation time of this algorithm may be estimated as follows. Since
$\cal S$ is finite, then the number of steps $J$ in the computation of $v'$ is
bounded by a linear function of the length of $v$. If the structure $\cal L$ is
asynchronous, then $\Length(v^m)$ is growing exponentially, and we have an
exponential time algorithm for the computation of $v^M$. If the structure $\cal
L$ is synchronous, then $\Length(v^m)$ is growing linearly, and we have a
quadratic time algorithm for computing $v^M$. 

From this argument, the most we can conclude is that the algorithm described in
\section III runs in exponential time. However, we can perhaps do better using
the following ideas.

Let $\cal L$, $\cal L'$ be asynchronous automatic structures on a groupoid
$G$ with generating set $\cal A$. We say that $\cal L'$ is a
\em{factorization} of $\cal L$ if, for each $v \in \cal L$ and $v' \in
\cal L'$ such that $\overline v = \overline v'$, there exists a sequence
$0=n_0 < n_1 < \cdots < n_J = \Length(v)$ with steps of bounded length
such that $\overline v(n_j) = \overline v'(j)$ for $j=0,\ldots,J$. Suppose
moreover that $\cal L$ and $\cal L'$ both satisfy suffix uniqueness, and that
$\cal L'$ is an automatic structure. Then the above described algorithm for
the word problem, \em{using the asynchronous structure} $\cal L$, runs in
quadratic time, improving the a priori fact that the algorithm runs in
exponential time. The reason is that lengths of normal forms in $\cal L'$ grow
linearly, and the factors have bounded length, therefore lengths of normal
forms in $\cal L$ grow linearly, hence the algorithm runs in quadratic time. 

In the next section we use this technique for showing that normal forms in
$\cal L_0$ can be computed in quadratic time, by finding an automatic structure
$\cal L_1$ for $\mcgd$ that is a factorization of $\cal L_0$.

\head V. Dehn twists, synchronous normal forms, and quadratic computation time
\endhead

The key to understanding the synchronous normal forms is to see how Dehn twists
arise in the asynchronous normal forms $\cal L_0$. This is described in the
\em{Dehn twist lemma} of \cite{M}, which we review here. A word in $\cal L_0$
will be factored into subwords which represent either Dehn twists or fractions of
Dehn twists; this leads to a language $\cal L_1$ which is a factorization of
$\cal L_0$ as in the last section. The properties of this factorization are used
in \cite{M} to prove that $\cal L_1$ is an automatic structure for $\mcgd$. For
our present purposes, we use the language $\cal L_1$ to prove that the
algorithm described above, for computing normal forms in $\cal L_0$, runs in
quadratic time:

\proclaim{Theorem: Quadratic computation time} Given a word $v$ of relabelling
moves and labelled elementary moves, with $\Length(v) = K$, the algorithm
computes the normal form of $v$ using at most $(12g-12)K^2$ relators. 
\endproclaim

Comment 1: As we shall see, the number $12g-12$ is the maximum length of a factor
in the Dehn twist factorization.

Comment 2: At any stage of the algorithm, the time needed to apply the next
relator is bounded by a constant, hence the algorithm computes the normal form
of $v$ in quadratic time.

Comment 3: For genus 2 we need at most $12K^2$ relators. In any
given run of the algorithm, the author is able to apply the required relator
using at most 2 minutes of time, leading to a computation time of at most
$25K^2$ minutes (experience shows that this is a very conservative estimate).

\subhead Dehn twist blocks \endsubhead

Consider an ideal triangulation $\delta$ and a prong of $\delta$ marked with a
$*$. Then $(\delta,*)$ is the special case of a marked ideal triangulation, with
no end markings; we call this a \em{prong marked ideal triangulation} or \prmit.
We use $\Delta$ to denote $(\delta,*)$. Choose a parity $d \in \{L,R\}$; we use
$\lnot d$ to denote the opposite parity. Consider the arc $h^{\lnot d}$ of
$\delta$, equipped with a transverse orientation pointing into the marked
triangle. The arc $h^{\lnot d}$ forms a simple closed curve in $S$, whose
regular neighborhood $N$ is an annulus, and the transverse orientation points
towards one of the boundary components of $N$, a simple closed curve we denote
$\gamma = \gamma(\Delta,d)$. Let $\tau = \tau(\Delta,d)$ be the Dehn twist
of parity $d$ around $\gamma$; see figure 65. Our parity convention for Dehn
twists is that a Right Dehn twist is a positive one, i.e.\ on an oriented
annulus $A = \reals \cross [0,1] / (x,y)\equiv (x+1,y)$ forming a regular
neighborhood of $\gamma$, the Dehn twist is given by the linear map $(x,y)
\mapsto (x+y,y)$ which takes a vertical segment to a segment of positive slope,
i.e.\ a segment that slopes up and to the Right. 

\midinsert
\centeredepsfbox{65.DehnTwist.eps}
\botcaption{Figure 65} The Dehn twist $\tau(\Delta,R)$
\endcaption\endinsert

The following lemma is basically the first part of the \tit{Dehn twist lemma} of
\cite{M}:

\proclaim{Dehn twist lemma (part I)} Let $\,\cdots \mapright{d} \Delta_2
\mapright{d}\Delta_1\mapright{d}\Delta_0$ be a (finite or infinite) sequence of
$d$-elementary moves ending with $\Delta =\Delta_0$. Let $\tau =
\tau(\Delta,d)$ be the Dehn twist defined above. Then there exists a
constant $K=K(\Delta,d)$, depending only on the combinatorial type of $\Delta$,
such that $\tau(\Delta_i) = \Delta_{i+K}$ for all $i \ge 0$.
\endproclaim

The proof is sketched below.

The sequence $\Delta_K \mapright{d} \cdots \mapright{d} \Delta_0$ is called a
\em{Dehn twist sequence} of parity $d$. The number $K$ is called the \em{Dehn
twist length}. Note that after taking combinatorial types, then $\{ \Delta_K \}
\mapright{d} \cdots \mapright{d} \{ \Delta_0 \}$ is a closed path in $\cal M_0$,
lying entirely in $\cal M_0^1$: the path is closed because $\Delta_0$ and
$\Delta_K = \tau(\Delta_0)$ have the same combinatorial type, and it is in level
1 because there are no end markings. We call this a \em{Dehn twist block} of
parity $d$ in $\cal M_0^1$; later, after putting in end markings, we shall
define Dehn twist blocks in higher levels. 

Some examples are given in figures 66--69. Figure 66 shows a Left Dehn twist
block on a torus, of length 1. All Left Dehn twist blocks in $\cal M_0^1$ on a
torus are orientation preserving conjugate to this one, and all Right Dehn twist
blocks are orientation reversing conjugate. 

\midinsert
\centeredepsfbox{66.TorusDehnTwistBlock.eps}
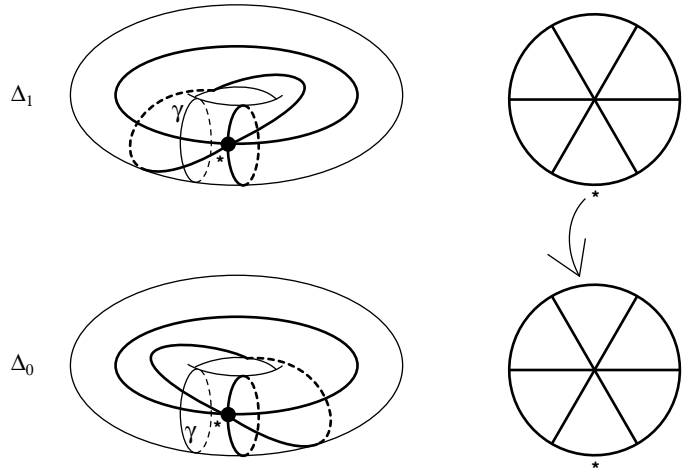
\botcaption{Figure 66} A Left Dehn twist block on a torus: $\Delta_1$ is
obtained from $\Delta_0$ by a Left Dehn twist about the curve $\gamma$.
\endcaption\endinsert

Figure 67 shows a Right Dehn twist block of length 1, on a surface of genus 2,
obtained by putting a marked prong into the second elementary move of figure
5. This example exhibits a characteristic property of Dehn twist sequences of
length 1 when the genus is at least 2: the Dehn twist length is 1 if and only if
$h^{\lnot d}$ is a boundary arc of a 1-handle piece, and the $*$ is inside the
1-handle piece. Figure 68 shows an example of a Right Dehn twist block of length
2.

\midinsert
\centeredepsfbox{67.Genus2TwistBlock1.eps}
\botcaption{Figure 67} A Right Dehn twist block of length 1 on a surface of
genus 2
\endcaption\endinsert

\midinsert
\centeredepsfbox{68.Genus2TwistBlock2.eps}
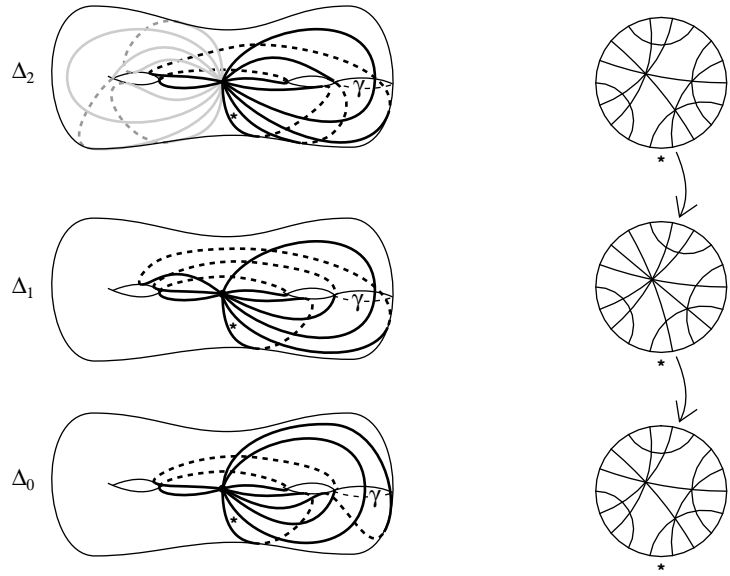
\botcaption{Figure 68} A Right Dehn twist block of length 2 on a surface of
genus 2. The shaded arcs of $\Delta_2$ should be included in $\Delta_1$ and
$\Delta_0$ as well, but are omitted from the diagram for clarity.
\endcaption\endinsert

In figure 69, note that a single elementary move returns to the original state,
completing a simple closed loop $w$ in $\cal M_0$ and thereby defining a mapping
class $\Phi$, but $\Phi$ is not a Dehn twist. The separating closed curve
$\gamma$ is invariant under $\Phi$; on the left component the restricted mapping
class is the identity, and on the right component the restricted class has some
finite order $k$ whose value we leave as an exercise. The mapping class $\Phi^k$
is a Dehn twist $\tau$, and the non-simple closed loop $w^k$ in $\cal M_0$ is a
Dehn twist block of length $k$. 

\midinsert
\centeredepsfbox{69.WhatIsTheLength.eps}
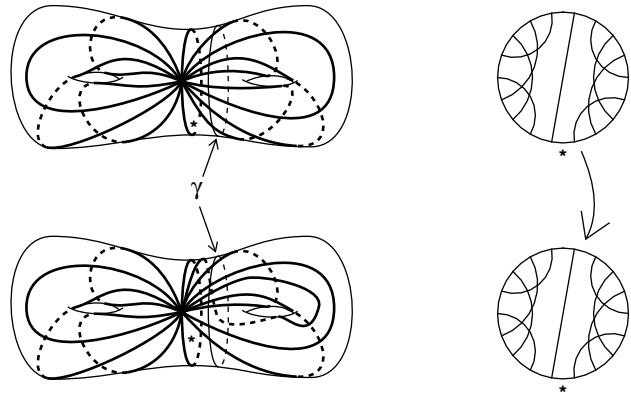
\botcaption{Figure 69} What is the Dehn twist length? (Hint: see figure 24c)
\endcaption\endinsert

Now we give the formula for Dehn twist lengths $K(\Delta,d)$, and we sketch the
proof of the Dehn twist lemma. The curve $\gamma$ cuts off certain half-arcs of
$\delta$, namely those half-arcs in the annulus bounded by $\gamma$ and
$h^{\lnot d}$. These half-arcs determine a subset of $\Ends(\delta)$ denoted
$\Ends^* = \Ends^*(\Delta,d)$. In figure 65 these half-arcs are labelled with a
$*$. In figures 66 and 67, $\bigl| \Ends^* \bigr| = 2$. In figure 68, $\bigl|
\Ends^* \bigr| = 3$. And in figure 69, $\bigl| \Ends^* \bigr| = 8$. Note that the
arc $h^\Opp$ always has at least one end in $\Ends^*$. Define $\hat{\Ends}^* =
\hat{\Ends}^*(\Delta,d) = \Ends^* - \Ends(h^\Opp)$. Let $K(\Delta,d) = \bigl|
\hat{\Ends}^* \bigr|$. Thus, $K(\Delta,d) = |\Ends^*| - 2$ or $|\Ends^*| - 1$
depending on whether or not $h^\Opp$ has one or two ends in $\Ends^*$. In
figures 66, 67, 68 only one end of $h^\Opp$ is in $\Ends{}^*$, so $K(\Delta,d) =
1,1,2$ respectively. In figure 69 both ends of $h^\Opp$ are in $\Ends{}^*$, so
$K(\Delta,d) = 6$. 

The number $K(\Delta,d)$ may be computed from the chord diagram of $\Delta$ as
follows. Recall that the ``marked triangle'' is the triangle having the marked
prong as a corner. Locate the chord corresponding to $h^{\lnot d}$. The endpoints
of this chord separate the remaining chord ends into two subsets; the subset
containing $e^d$ corresponds to $\Ends{}^*$. Now count the number of elements in
$\Ends{}^*$, subtract 1 if the marked triangle is twisted (because then
only one end of $h^\Opp$ is in $\Ends{}^*$), and subtract 2 if the marked
triangle is untwisted (because both ends of $h^\Opp$ are in $\Ends{}^*$);
the result is $K(\Delta,d)$; see figure 70. In figures 66-69, only figure 69
has an untwisted marked triangle, hence only in that case is 2 subtracted to
compute $K(\Delta,d)$; in the other cases 1 is subtracted. Other examples of
computing $K(\Delta,d)$ are given in figure 71.

\midinsert
\centeredepsfbox{70.Subtracting1or2.eps}
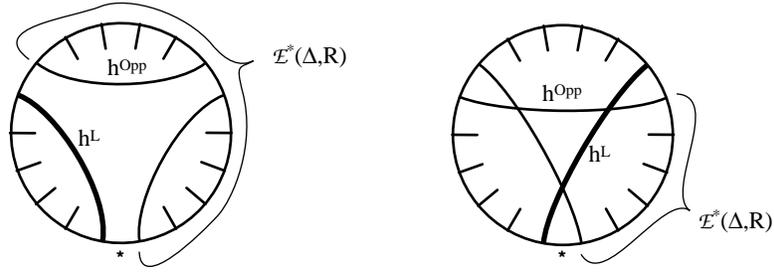
\botcaption{Figure 70} If the marked triangle is untwisted then both ends of
$h^\Opp$ are in $\Ends{}^*$; but if it is twisted only one end of $h^\Opp$ is in
$\Ends{}^*$.
\endcaption\endinsert

\midinsert
\centeredepsfbox{71.ComputingKDelta,d.eps}
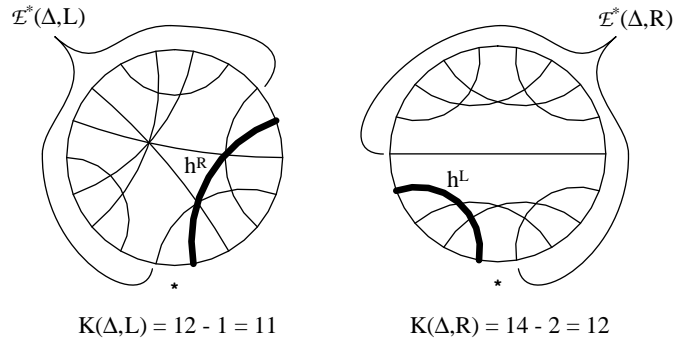
\botcaption{Figure 71} Examples of computing $K(\Delta,d)$
\endcaption\endinsert

The key observation in proving the Dehn twist lemma is that for any prong marked
ideal triangulation $\Delta$ and any $d \in \{L,R\}$, the arc $h^\Opp$ is
obtained up to isotopy from $h^d$ by letting the Dehn twist
$\tau(\Delta,d)^\inverse$ act on a half-arc representing $e^d$, where the
half-arc is chosen to intersect $\gamma$ exactly once. This is illustrated in
figure 72, which shows separately the cases where the marked triangle is twisted
and untwisted. Note in the twisted case that $h^d$ and $h^\Opp$ each have a
unique end in $\Ends{}^*$, and the twist about $\gamma$ takes $h^d$ to
$h^\Opp$. But in the untwisted case where $h^d$ and $h^\Opp$ have both ends
in $\Ends{}^*$, the twist \em{does not} take $h^d$ to $h^\Opp$; by allowing the
twist to act \em{only} on the end $e^d$ of $h^d$, we thereby obtain $h^\Opp$.

\midinsert
\centeredepsfbox{72.ProofOfDehnTwistLemma.eps}
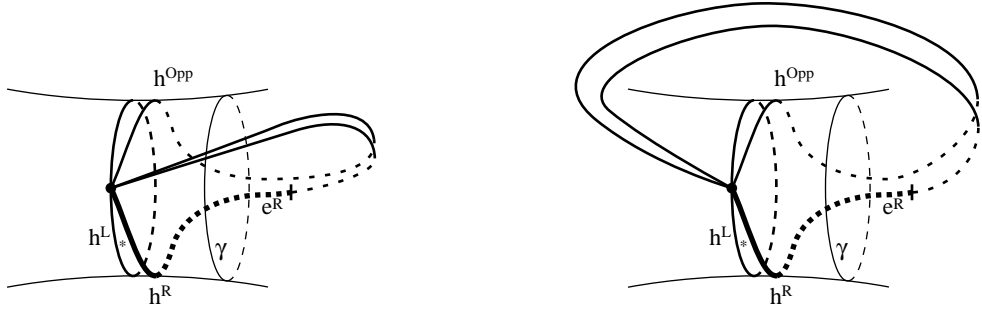
\botcaption{Figure 72} Up to isotopy, $h^\Opp$ is obtained from $h^R$ by
letting $\tau(\Delta,R)^\inverse$, the Left Dehn twist about $\gamma$, act on the
half-arc representing $e^R$.
\endcaption\endinsert

This observation is applied as follows. Let $\hat\Delta_i$ be obtained from
$\Delta_i$ by removing $h_i^\Opp$. After removing $h_i^\Opp$, the marked prong
is now located in a complementary 4-gon of $\hat\Delta_i$. Note that $\Delta_i$
may be recovered from $\hat\Delta_i$ by triangulating this 4-gon using an arc
opposite the marked prong. Consider the good elementary move $\Delta_{i+1} \to
\Delta_i$. Note that $\hat\Delta_i$ is obtained from $\hat\Delta_{i+1}$ by
inserting $h_{i+1}^\Opp$ and then removing $h_{i+1}^d$. Using the key
observation, it follows that $\hat\Delta_i$ is obtained from $\hat\Delta_{i+1}$
by letting $\tau^\inverse$ act on the half-arc representing $e_{i+1}^d$.
Travelling along the sequence $\hat\Delta_K, \ldots, \hat\Delta_1, \hat\Delta_0$,
then $\tau^\inverse$ acts in turn on a representative half-arc of each end in
$\hat{\Ends}_K^*$. Since these half-arcs represent all the points in
$\hat\Delta_K \intersect \gamma$, it follows that $\tau^\inverse(\hat\Delta_K) =
\hat\Delta_0$, hence $\tau^\inverse(\Delta_K) = \Delta_0$ so $\tau(\Delta_0) =
\Delta_K$. In order to get the full periodicity statement $\tau(\Delta_i)
=\Delta_{i+K}$, note that $h_i^{\lnot d} = h_{i+1}^{\lnot d}$ so
$\tau(\Delta_i,d) = \tau(\Delta_{i+1},d)$, and also $K(\Delta_i,d) = \bigl |
\hat{\Ends}^*_i \bigr| = \bigl| \hat{\Ends}^*_{i+1} \bigr| = K(\Delta_{i+1},d)$.
This finishes the proof of the Dehn twist lemma.

We can now determine the range of possible Dehn twist lengths. Obviously
$K(\Delta,d) \ge 1$. To determine when equality is acheived, note that $|
\Ends^*(\Delta,d)| \ge 2$ with the minimum acheived if and only if $h^{\lnot
d}$ is the boundary of a 1-handle piece and the $*$ is inside the 1-handle
piece, in which case the marked triangle is inside the 1-handle piece and
therefore twisted; this is the only way the Dehn twist length can be 1, because
if $|\Ends^*(\Delta,d)| = 3$ then the marked triangle is still twisted so the
Dehn twist length is 2. To find the maximum Dehn twist length, recall that
$\Ends(\delta) = 12g-6$. Now $h^{\lnot d}$ cuts off one subset of
$\Ends(\delta)$ of size at least 2, and $h^{\lnot d}$ itself has 2 ends, hence
$\bigl| \Ends^*(\Delta,d) \bigr| \le 12g - 6 - 2 - 2 = 12g-10$. This size is
acheived if and only if $h^{\lnot d}$ is on the boundary of a 1-handle piece
and the $*$ is outside the 1-handle piece, in which case the marked triangle is
untwisted and so the Dehn twist length is $12g - 10 - 2 = 12g-12$. Also,
when $|\Ends^*(\Delta,d)| = 12g-11$ then the Dehn twist length is at most
$12g-11-1=12g-12$. Therefore we have an optimal upper bound of $12g-12$ for the
Dehn twist length. An example of a maximal length Dehn twist sequence on a
surface of genus 2 is given in figure 72. It is an exercise to show that this is
the unique maximal length Right Dehn twist block in genus 2.

\midinsert
\centeredepsfbox{73.Length12TwistBlock.eps}
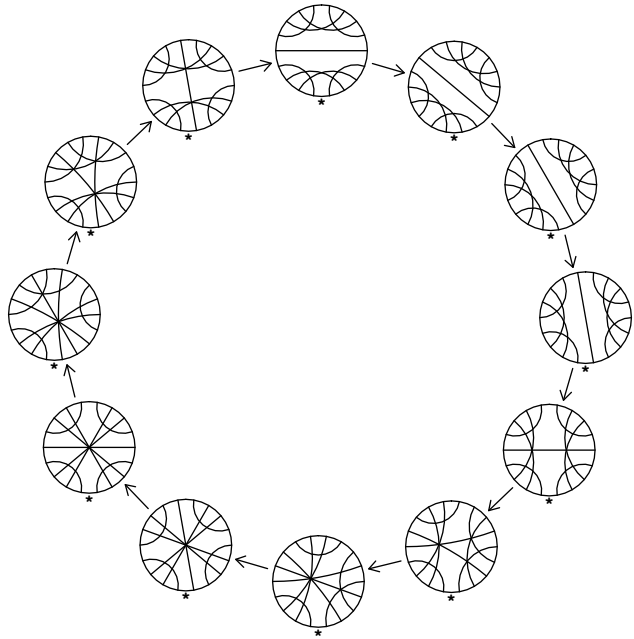
\botcaption{Figure 73} A Right Dehn twist block in genus 2, of maximal length
12
\endcaption\endinsert

\subhead An automaton for synchronous normal forms \endsubhead
First we define a new generating set $\cal A_1$ for $\mcgd$, over which the new
automaton $\cal M_1$ is defined. The set $\cal A_1$ is obtained from $\cal A_0$
by adding Dehn twist generators and fractions thereof.

In the last section we defined Dehn twist blocks in $\cal M_0^1$. Define a
\em{Dehn twist generator} to be the path in $X$ obtained by projecting a Dehn
twist block from $\cal M_0^1$ to $X$. Every Dehn twist generator is a closed
curve in $X$, so it defines a group element in $\mcg$. Define a \em{fractional
Dehn twist generator} to be any subword of a Dehn twist generator; this may not
be a closed curve, and so may not define a group element.

Define a new alphabet $\cal A_1$ to be the set of all relabelling generators,
Dehn twist generators, and fractional Dehn twist generators. Note that every
labelled elementary move is either a full or fractional Dehn twist generator, so
$\cal A_0 \subset \cal A_1$. There is a map from $\cal A_1$ to $\mcgd$, taking
each generator to its homotopy class.

Consider a full or fractional Dehn twist generator $w = w_1 \cdots w_n$. Let
$s_0 \mapright{w_1} s_1 \mapright{w_2} \cdots \mapright{w_n} s_n$ be any path in
$\cal M$ that lifts $w$ and stays among the accept states. This path is called a
full or fractional \em{Dehn twist block} if the subpath
$s_1 \to \cdots \to s_n$ stays in a single level of $\cal M$; we allow $s_0 \to
s_1$ to drop between levels.

Chord diagrams of Dehn twist blocks in arbitrary levels are understood as
follows. Suppose that $s_0 \mapright{d} s_1 \mapright{d} \cdots \mapright{d}
s_n$ is a full Dehn twist block in level 1, of parity $d$. Thus, none of the
chord diagrams have labelled ends, and each $s_i \to s_{i+1}$ is a parity $d$
arrow. Dehn twist blocks in higher levels are obtained by introducing end
markings, as follows. Recall the set of chord ends $\Ends{}^* =
\Ends{}^*(s_0,d)$: the two ends of
$h_0^{\lnot d}$ divide the remaining chord ends of $s_0$ into two subsets, one
of which contains $e^d$, that subset being $\Ends{}^*(s_0,d)$. If end markings
$1,\ldots,k-1$ are introduced to form a new state $s'_0$, and if none of the
labelled ends are in $\cal E^* \union \OppEnd(\cal E^*)$, then we obtain a Dehn
twist block $s'_0 \mapright{d} s'_1 \mapright{d} \cdots \mapright{d} s'_n$
staying entirely in level $k$; figure 74 shows an example. Then if additional end
markings
$k,\ldots,l-1$ are inserted so that the end $e^d$ is marked with $k$, and if a
different marked prong is then chosen, we obtain a Dehn twist block $s'_0
\mapright{d} s'_1 \mapright{d} \cdots
\mapright{d} s'_n$ where the first move drops from level $l$ to level $k$ and
the rest of the block stays in level $k$. Figure 75 shows an example, adding end
markings to the example from figure 74.

\midinsert
\centeredepsfbox{74.HigherLevelTwistBlock.eps}
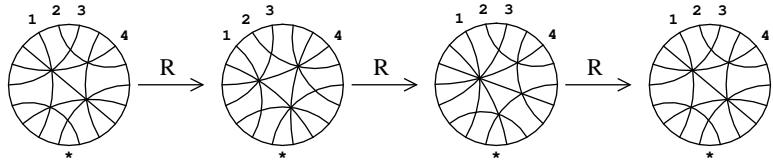
\botcaption{Figure 74} A Right Dehn twist block in level 5
\endcaption\endinsert

\midinsert
\centeredepsfbox{75.DroppingLevelTwistBlock.eps}
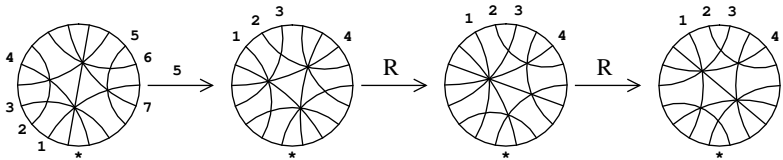
\botcaption{Figure 75} A Right Dehn twist block which drops from level 8 to
level 5
\endcaption\endinsert

Chord diagrams of fractional Dehn twist blocks are constructed similarly,
except that the restrictions on end markings are somewhat weaker. Suppose $s_0
\to s_1\to\cdots \to s_n$ is a fractional Dehn twist block in level 1, of parity
$d$. Instead of worrying about all of $\Ends{}^*$, only worry about those ends in
$\Ends{}^*$ which will eventually become $e^d$ for one of the states
$s_0,\ldots,s_{n-1}$, i.e.\ those ends which lie on a chord that will eventually
be removed in performing one of the elementary moves in the block; let that set
be denoted $\Ends{}^\#$. To get a fractional Dehn twist block that stays in level
$k$, we may mark chord ends of $s_0$ with $1,\ldots,k-1$ as long as the marked
ends are not in $\Ends{}^\# \union \OppEnd(\Ends{}^\#)$. To get a fractional
Dehn twist block that drops from level $l$ to level $k$, add more end markings
$k,\ldots,l-1$ so that $e^d$ is marked with $k$, and then move the marked prong
if desired. An example of a fractional Dehn twist block in level 6 is given in
figure 76, adding an end marking to the first two moves in figure 74. This
example cannot be extended to any longer full or fractional Dehn twist block,
because in the final state of the block the end $e^R$ in the final state is
marked with a 5, so the move on $h^R$ drops down to level 5.

\midinsert
\centeredepsfbox{76.FractionalTwistBlock.eps}
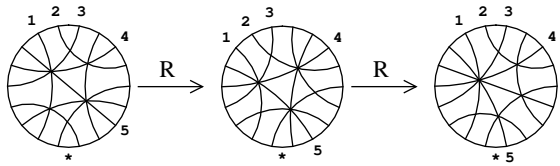
\botcaption{Figure 76} A Right fractional Dehn twist block in level 6, which
cannot be extended to a longer block
\endcaption\endinsert

Now we are ready to define the automaton $\cal M_1$. Its state set is (almost)
the same as the state set of $\cal M_0$. In addition to the old arrows of $\cal
M_0$, we add new arrows representing Dehn twist blocks and fractional Dehn
twist blocks, which jump over the corresponding path in $\cal M_0$. In order to
preserve uniqueness, we require that fractional Dehn twists can occur only at the
very beginning of a parity block, so the states of $\cal M_1$ must remember
whether a full Dehn twist has just occured.

To define the states of $\cal M_1$: the failure states, and
inconsistent accept states are the same as for $\cal M_0$. For every consistent
accept state $s$ of $\cal M_0$, we define three accept states in $\cal M_1$,
namely $(s,O)$, $(s,L)$, and $(s,R)$, whose meanings are as follows. In state
$(s,d)$ the previous letter was a full Dehn twist of parity $d$, and in state
$(s,O)$ the previous letter was not a full Dehn twist. If $s_0$ is the start
state of $\cal M_0$ then $(s_0,O)$ is the start state of $\cal M_1$.

Now we define arrows of $\cal M_1$. The arrows coming out of failure states
and inconsistent states all lead to failure states as before. For every
relabelling arrow $s \to s'$ of $\cal M_0$, noting that $s$ is a consistent
state and $s'$ is inconsistent, we define three relabelling arrows
$(s,O) \to s'$, $(s,L) \to s'$, and $(s,R) \to s'$, all named with the same
relabelling generator.

Consider now a consistent accept state $s_0$ of $\cal M_0$. Consider also a full
or fractional Dehn twist block $s_0 \to s_1\to \ldots \to s_n$ of parity $d$,
ending at $s_n$ in level $k$, and let $w$ be the full or fractional Dehn twist
generator to which this block projects. If $w$ is full, construct arrows from
the states $(s_0,O), (s_0,R), (s_0,L)$ to the state $(s_n,d)$, all named with
$w$. If $w$ is fractional and $s_0$ is in level $k$, construct arrows from the
states $(s_0,O), (s_0,\lnot d)$ to
$(s_n,O)$, both named with $w$; the arrow from $(s_0,d)$ named with $w$ leads
to the appropriate failure state. If $w$ is fractional and $s_0$ is not in level
$k$, construct arrows from the states $(s_0,O), (s_0,L), (s_O,R)$ to $(s_n,O)$,
all named with $w$. All other arrows which have not been specifically
constructed here should lead to the appropriate failure state. 

The effect of this construction is that a fractional Dehn twist block cannot
follow a full Dehn twist block of the same parity, unless the fractional block
drops down to a lower level. 

The language $\cal L_1$ accepted by $\cal M_1$ is related to the language $\cal
L_0$ in the following manner. Given any word $w \in \cal L_0$, recall that $w$
is factored into uncombing blocks, and each uncombing block is in turn factored
into parity blocks. Now look at a parity block. Using the Dehn twist lemma,
that block can be factored in a unique manner as a fractional Dehn twist block,
followed by some number of full Dehn twist blocks. Doing this factorization for
each parity block in $\cal L_0$, one obtains a word $w'$ in the generators $\cal
A_1$, such that $w$ and $w'$ represent the same element of $\mcgd$, and $w' \in
\cal L_1$. We say that $w'$ is the \em{Dehn twist factorization} of $w$.

\subhead Synchronization \endsubhead

In this section we review the results of \cite{M} that are used to prove that
$\cal L_1$ satisfies the fellow traveller property. We shall use these results
to prove directly that our algorithm for the word problem runs in quadratic
time. 

Given a word $w := w_1\cdots w_k$, recall the notation for a prefix subword $w(t)
:= w_1\cdots w_t$. We also use notation for an infix subword $w[i,j] := w_{i+1}
\cdots w_j$. Note that if $a<b<c$ then $w[a,b] w[b,c] = w[a,c]$ and $w(a)
w[a,b] = w(b)$.

Consider two normal forms $v,w \in \cal L_0$, such that $\alpha = \overline
v^\inverse \overline w$ is an elementary move generator. In applying the
subroutine \tit{Do one move} to the word $v\alpha$, the algorithm applies some
number of relators to obtain $w$, and we want to estimate that number. The most
interesting case is when the algorithm classifies the move $\alpha$ as a bad
elementary move, in which case the estimate will arise by studying how Dehn twist
boundaries interact with the operations of the algorithm. Understanding this
interaction is also the key to proving the synchronous fellow traveller property
for $\cal L_1$.

Let $v',w' \in \cal L_1$ be the Dehn twist factorizations of $v,w$, with
$\Length(v') = M$ and $\Length(w') = N$. Thus, we may write $v = v[s_0,s_1]
\composed \cdots \composed v[s_{M-1},s_M]$ where $v'_m = v[s_{m-1},s_m]$, and
similarly $w'_n = w[t_{n-1},t_n]$. We think of the parameter values
$s_0,\ldots,s_M$ and $t_0,\ldots,t_N$ as ``Dehn twist boundaries''. We regard
$v'_m$ and $w'_n$ as individual letters of $\cal A_1$, but they may also be
regarded as words in $\cal A_0$ and as such we may speak of their subwords.

To start, we study how Dehn twist boundaries interact under the situation where
$\alpha$ is a good or inverse good elementary move, or more generally when one of
$v,w$ is a prefix subword of the other. Supposing $v$ is a prefix subword of
$w$, then there exists $n \le N-1$ such that $t_n \le s_M < t_{n+1}$, and we say
that \em{$w'$ extends $v'$ by $N-n$ Dehn twist units}, namely $w[s_M,t_{n+1}],
w'_{n+2}, \ldots, w'_N$. All but the first of these units are letters
of $w'$; the first unit is a suffix of $w'_{n+1}$, possibly the whole word, but
either way it forms a letter in the alphabet $\cal A_1$. To compare Dehn twist
boundaries in this situation:

\proclaim{Good proposition} Suppose $w'$ extends $v'$ by $K$ Dehn twist units.
Then $w'(N-t)$ extends $v'(M-t)$ by at most $K$ Dehn twist units, for all $t \ge
0$.
\endproclaim

\demo{Proof} The boundaries of the Dehn twist factorizations of $v$ and $w$ are
identical up until the last parity block in $v$; suppose that $v'(A)$ and $w'(A)$
end at the beginning of that parity block. In that parity block, the initial
fractional Dehn twist factors of $v$ and $w$ may have different lengths, but for
the rest of the parity block the full Dehn twist factors have the same length,
hence if $A < B \le M$ then one of $v'(B)$ or $w'(B)$ extends the other by at
most one Dehn twist unit. It follows that if $t_n \le s_M < t_{n+1}$, then
$t_{n-1} \le s_{M-1} < t_{n+1}$, so $w'(N-1)$ extends $v'(M-1)$ by at most $K$
Dehn twist units. Now continue by induction. \qed\enddemo

Now suppose that $\alpha$ is a bad elementary move. The algorithm \tit{Do a bad
elementary move} applies a sequence of relators, producing a sequence of bad
elementary moves connecting shorter and shorter initial subwords of $v$ and $w$,
until reaching identical initial subwords. More precisely, there exists $P \ge
1$ and sequences $0 \le i_0 < i_1 < \cdots < i_P$, $0 \le j_0 < j_1 < \cdots <
j_P$ with the following properties:
\roster
\item"(1)" If $0 < p \le P$ then $v(i_p)$ and $w(j_p)$ differ by a bad elementary
move; we say that $v(i_p)$ and $w(j_p)$ are \em{matching \bem\ ends} (\bem\ is
the acronym for ``bad elementary move'').
\item"(2)" If $1 < p \le P$ then there exists $(a,b) \in \{(1,1), (1,2), (2,1)\}$
such that $i_{p-1} + a = i_p$, and $j_{p-1} + b = j_p$.
\item"(3)" There exists $(a,b) \in \{(1,1),(1,2),(2,1),(2,2)\}$ such that $i_0
+ a = i_1$ and $j_0 + b = j_1$. Furthermore:
\itemitem{(3a)} If $(a,b) = (2,2)$ then $v(i_0) = w(j_0)$.
\itemitem{(3b)} If $(a,b) \in \{(1,1),(1,2),(2,1)\}$ then one of $v(i_0), w(j_0)$
extends the other by a single good elementary move. 
\item"(4)" In particular, $\left|i_0 - j_0\right| \le 1$. We refer to the
intervals $[i_0,i_1]$ and $[j_0,j_1]$ as the \em{irregular regions}.
\endroster
Remark: in \cite{M} we also say, in case (3b), that $v(i_0)$ and $w(j_0)$ are
\em{matching \bem\ ends} (despite the fact that their difference is a good or
inverse good elementary move).

The relation between the Dehn twist factorizations of $v$ and $w$ is given in
the following, which although not stated explicitly in \cite{M} is proved
implicitly:

\proclaim{Bad proposition} Suppose $\alpha = \overline v^\inverse \overline w$
is a bad elementary move. Recall the notation $s_0,\ldots,s_M$ and
$t_0,\ldots,t_N$ for the Dehn twist boundaries of $v,w$, and note that $v(s_M)$
and $w(t_N)$ are matching \bem\ ends. Then there exists a constant $A \ge 0$
such that
$v(s_{M-a})$ and $w(t_{N-a})$ are matching \bem\ ends for $0 \le a < A$. Moreover,
one of the following happens:
\roster
\item"(1)" One of $v(s_{M-A})$ and $w(t_{N-A})$ extends the other by one Dehn
twist unit.
\endroster
\noindent 
or
\roster
\item"(2)" $v(s_{M-A})$ and $w(t_{N-A})$ differ by 2 elementary moves (one bad
and one good), and one of $v(s_{M-A-1})$ and $w(t_{N-A-1})$ extends the other by
at most two Dehn twist units.
\endroster\endproclaim

The point of this proposition is that as you move backwards along $v$ and $w$,
moving synchronously one Dehn twist block per step, then the corresponding Dehn
twist boundaries will be matching \bem\ ends, hence differing by a single
generator in $\cal A_0$. This continues until you reach the irregular regions,
at which time the difference can become as large as two Dehn twist units. 

\demo{Sketch of proof} We have defined three progressively finer factorizations
of normal forms in $\cal L_0$: the uncombing block factorization, the parity
block factorization, and the Dehn twist factorization. Corresponding to each of
these is a proposition in \cite{M} which describes the interaction of the
factorization with \bem\ ends: \tit{\bems respect combing blocks}, \tit{\bems
respect parity blocks}, and the second part of the \tit{Dehn twist lemma}. We
invoke these in the proof.

Now look at the final Dehn twist blocks $v[s_{M-1},s_M]$ and $w[t_{N-1},t_N]$,
and go case by case through the different possibilities.

If $v[s_{M-1},s_M]$ is a full uncombing block, then by \tit{\bems respect
combing blocks} it follows that $w[t_{N-1},t_N]$ is also a full combing block,
and $v(s_{M-1}), w(t_{N-1})$ differ by a single elementary move, either good,
bad, or inverse good. If it is good or inverse good, then evidently one of
$v'_{M-1}$ or $w'_{N-1}$ extends the other by one Dehn twist unit, proving item
(1) of the \tit{Bad proposition}. It it is bad, then the proof continues by
induction.

If $v[s_{M-1},s_M]$ is a full Dehn twist block, then by the second part of the
\tit{Dehn twist lemma} it follows that $w[t_{N-1},t_N]$ is also a full Dehn
twist block, and $v(s_{M-1}), w(t_{N-1})$ differ by a single bad elementary
move. The proof now continues by induction. 

The remaining case is where $v[s_{M-1},s_M]$, $w[t_{N-1},t_N]$ are fractional
Dehn twist blocks which are not full uncombing blocks. In this case,
either $v(s_{M-1})$ and $w(t_{N-1})$ differ by a bad elementary move, or both are
in the irregular regions; this follows from \tit{\bems respect parity blocks}.
When they differ by a bad elementary move, continue by induction as before. 

When both $s_{M-1}$ and $t_{N-1}$ are in the irregular regions, then the proof
of synchronization in \cite{M} analyzes carefully where the Dehn twist boundaries
may occur. Roughly speaking, since the two sides of the relation in the
irregular regions are quite short, they cannot throw off the synchronization by
too much. The conclusions of the argument from \cite{M} are as follows, proving
item (2) of the \tit{Bad proposition}: the irregular region is of type $(2,2)$ as
in (3a) above; the relation which applies is always of type IIbii; the words
$v'(M-1)$ and $w'(N-1)$ differ by two elementary moves (one bad and one good);
and one of $v'(M-2)$, $w'(N-2)$ extends the other by at most two Dehn twist
units. This argument is summarized in figure 18 of \cite{M}, the last figure of
section III.2. There are three cases to the argument, and in figure 77 we
present examples for each of the three cases, paralleling the schematic pictures
in figure 18 of \cite{M}. These cases are distinguished as follows. We assume
that all the arrows in the relator IIbii are parity arrows, with parities LR on
the left side and RL on the right side. The arrows below the relator on the left
and right have the same parity, say L; these arrows are part of the final
letters of $v'_M$ and $w'_N$. The three cases are distinguished by whether the
arrow above the relator is no parity, Left parity, or Right parity. If no
parity, the example in figure 77a is typical: $v'(M-2)$ extends $w'(M-2)$ by at
most two Dehn twist units, each unit being a single elementary move; it could
happen that the no parity arrow pictured is an entire Dehn twist block, in which
case the Dehn twist boundary
$t_{N-1}$ would be one arrow lower in figure 77a, and $v'(M-2)$ would extend
$w'(M-2)$ by a single Dehn twist unit. If Left parity, the example in figure 77b
is typical: $v'(M-2)$ extends $w'(M-2)$ by one Dehn twist unit, consisting of a
single elementary move. If Right parity, the example in figure 77c is typical:
$v'(M-2)$ extends $w'(M-2)$ by at most two Dehn twist units, one being a single
elementary move and the second being a fractional Dehn twist generator; the
example in figure 77c shows the second unit being a single elementary move, but
the letter $w'(M-2) = w[t_{N-2},t_{N-1}]$ can be any full or fractional Dehn
twist generator, in which case the second unit can be an arbitrarily long
fractional Dehn twist generator. The general argument given in \cite{M} shows
that figures 77a--c are typical: the irregular regions affect the Dehn twist
boundaries in one of the three ways exemplified in these figures.

\midinsert
\centeredepsfbox{77.IrregularRegions.eps}
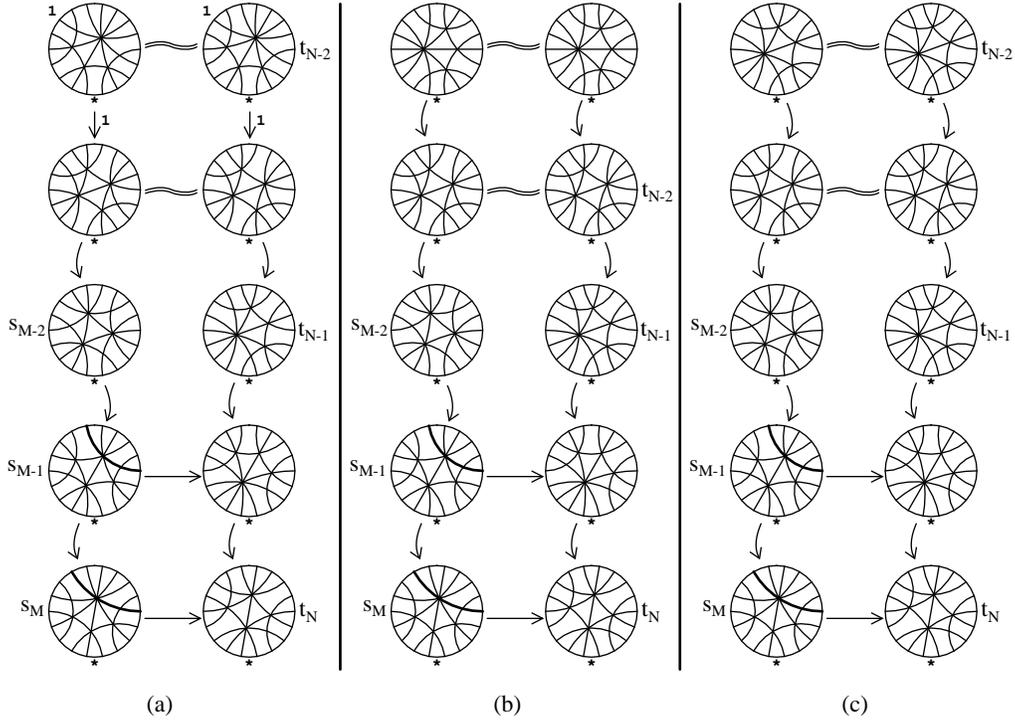
\botcaption{Figure 77} Dehn twist boundaries in the irregular regions
\endcaption\endinsert

\qed\enddemo

The \tit{Good and Bad propositions} may be used to prove the synchronous fellow
traveller property for $\cal L_1$, as follows. Consider $v',w' \in \cal L_1$,
with $M = \Length(v')$ and $N = \Length(w')$. Let $d_i$ denote distance
measured with the generating set $\cal A_i$. Noting that each letter in $\cal
A_1$ is a word of length at most $12g-12$ in the letters of $\cal A_0$, to prove
the fellow traveller property it suffices to assume that $d_0(\overline
v',\overline w') = 1$ and prove that $d_1(\overline v'(t),\overline w'(t)
\le 4$ for all $t$. The case where $\overline v',\overline w'$ differ by a
relabelling generator is easy. 

Suppose $\overline v',\overline w'$ differ by an elementary move generator.
Applying the \tit{Good and bad propositions} we see that $d_1(\overline
v'(M-t),\overline w'(N-t))\le 2$ for all $t$. This shows that $|M-N| \le 2$, so
$d_1(\overline v'(t),\overline w'(t)) \le 4$ for all $t$.

From this argument we obtain the following important fact which is needed in
estimating computation time:

\proclaim{Lemma: Length grows additively} Given $v,w \in \cal L_0$ such that
$d_0(\overline v,\overline w) \le 1$, if $v',w' \in \cal L_1$ are the Dehn
twist factorizations, then $\bigl|\Length(v') -\Length(w')\bigr| \le 2$.
\qed\endproclaim

\subhead Proof of theorem: Quadratic computation time \endsubhead
Consider $v,w \in \cal L_0$ with $d_0(\overline v,\overline w) \le 1$, and
$\alpha =\overline v^\inverse \overline w$. Let the Dehn twist factorizations be
$v',w'\in \cal L_1$. 

We claim that the number of relators used by the subroutine \tit{Do one move} to
compute $w$ from $v\alpha$ is at most $(12g-12)\Length(v')$. If the generator
$\alpha =\overline v^\inverse \overline w$ is a relabelling generator, then only
1 relator is used. If $\alpha$ is a good or inverse good elementary move, then at
most 2 relators are used. 

Suppose that $\alpha$ is a bad elementary move, and apply the \tit{Bad
proposition}. For $1 \le a < A$, consider the relators that are applied by the
algorithm to compute the (possibly fractional) Dehn twist block $w[t_{N-a},
t_{N-a+1}]$ from the Dehn twist block $v[s_{M-a}, s_{M-a+1}]$. These relators
are all of types Ia, IIa, or IIbi, each relator touching one of the good
elementary moves in
$v[s_{M-a}, s_{M-a+1}]$ and no two relators touching the same one, but there are
at most $12g-12$ elementary moves since this is a Dehn twist block or fraction
thereof. Thus, at most $12g-12$ relators are used. Also, the algorithm uses at
most $12g-12$ relators to compute $w[t_{N-A}, t_{N-A+1}]$ from $v[s_{M-A},
s_{M-A+1}]$: if these are full Dehn twist blocks then it follows as before; and
if these are fractional Dehn twist blocks then except for the top relator, each
relator touches some elementary move in $v[s_{M-A}, s_{M-A+1}]$ and no two
relators touch the same one, but there are at most $12g-13$ elementary moves in
a fractional Dehn twist block, and adding one more for the top relator gives
$12g-12$. This proves the claim.

Taking this claim together with the lemma \tit{Length grows additively}, we
reach the conclusion that for any word $v$ of length $K$ in the
generators $\cal A_0$, the number of relators used by the algorithm to compute
the normal form of $v$ is at most $(12g-12)[1+3+5+\cdots+(2K-1)] = (12g-12)K^2$.

Remark: The constant $12g-12$ can be improved slightly, by noticing that in
constructing the relators touching a Dehn twist block, if the block has full
length $12g-12$ then there must be at least one relator of type IIbi which
touches two moves in the block, so the number of relators adjacent to a Dehn
twist block is at most $12g-13$. Therefore the number of relators needed by
the algorithm is at most $(12g-13)K^2$.


\Refs\nofrills{Bibliography}

\widestnumber\key{ECHLPT}



\ref\key ECHLPT 
\by D. Epstein, J. Cannon, D. Holt, S. Levy, M. Paterson, W. Thurston 
\book Word processing in groups 
\publ Jones \& Bartlett 
\yr 1992 
\endref

\ref\key Har
\by J. Harer
\paper The virtual cohomological dimension of the mapping class group of an
oriented surface 
\pages 157--176
\yr 1986
\vol 84
\jour Invent. Math.
\endref

\ref\key Hat
\by A. Hatcher
\paper On triangulations of surfaces
\jour Topology Appl.
\vol 40
\issue 2
\yr 1991
\pages 189--194
\endref

\ref\key M
\by L. Mosher
\paper Mapping class groups are automatic
\yr 1993
\paperinfo Preprint
\endref

\ref\key STT
\by D. D. Sleator, R. E. Tarjan, W. P. Thurston
\paper Rotation distance, triangulations, and hyperbolic geometry
\jour J. Amer. Math. Soc.
\vol 1
\yr 1988
\pages 647--681
\endref
 


\endRefs
\enddocument